\documentclass[12pt]{article}

\usepackage[top=1.25in,bottom=1.25in,left=0.5in,right=0.5in]{geometry}
\usepackage{natbib}

\usepackage{tikz}
\usepackage{graphicx}
\usepackage{amssymb}
\usepackage{amsmath}
\usepackage{dsfont}
\usepackage{framed}
\usepackage{color}
\usepackage{caption}
\usepackage{float}
\usepackage{multirow}
\usepackage{rotating}
\usepackage{lscape}
\usepackage{booktabs}
\usepackage{subfigure}
\usepackage{adjustbox}
\usepackage{arydshln}
\usepackage{diagbox}
\usepackage{colortbl}
\usepackage[title,titletoc]{appendix}
\newtheorem{thm}{Theorem}[section]

\newtheorem{dfn}[thm]{Definition}
\newtheorem{cor}[thm]{Corollary}
\newtheorem{pro}[thm]{Proposition}
\newtheorem{lem}[thm]{Lemma}
\newtheorem{rem}[thm]{Remark}
\newtheorem{asm}[thm]{Assumption}
\usepackage{enumitem}
\def\argmin{\mathop{\rm arg\,min}}%

\def\ECSwitch{%
\clearpage
\normalsize
\pagestyle{ECheadings}%
\ECHowTheorems
\ECHowEquations
\ECHowSectionsAppendix
\setcounter{figure}{0}%
\renewcommand\thefigure{EC.\@arabic\c@figure}%
\setcounter{table}{0}%
\renewcommand\thetable{EC.\@arabic\c@table}%
\setcounter{page}{1}\def\thepage{ec\arabic{page}}%
\hspace*{1em}\relax}

\begin{document}

\title{\vspace{-0.5in}\Large
Information Relaxation and A Duality-Driven Algorithm for Stochastic Dynamic Programs
}
\author{
\ Nan Chen\thanks{Department of Systems Engineering and Engineering Management, The Chinese University of Hong Kong, Shatin, N. T., Hong Kong. Email: nchen@se.cuhk.edu.hk}
\and
\ Xiang Ma\thanks{Department of Systems Engineering and Engineering Management, The Chinese University of Hong Kong, Shatin, N. T., Hong Kong. Email: xma@se.cuhk.edu.hk}
\and
\ Yanchu Liu\thanks{Department of Finance, Lingnan (University) College, Sun Yat-sen University, Guangzhou, China. Email: liuych26@mail.sysu.edu.cn}
\and
\ Wei Yu\thanks{World Quant (Singapore), 1 Wallich Street, \#20-01 Guoco Tower, Singapore. Email: wei.yu@worldquant.com}
}

% \date{July 8, 2019}
\date{This Version: July 24, 2020\\
First Version: July 8, 2019}
\maketitle
\vspace{-0.2in}

\begin{abstract}
We use the technique of information relaxation to develop a duality-driven iterative approach to obtaining and improving confidence interval estimates for the true value of finite-horizon stochastic
dynamic programming problems. We show that the sequence of dual value estimates yielded from the proposed approach in principle monotonically converges to the true value function in a finite
number of dual iterations. Aiming to overcome the curse of dimensionality in various applications, we also introduce a regression-based Monte Carlo algorithm for implementation. The new
approach can be used not only to assess the quality of heuristic policies, but also to improve them if we find that their duality gap is large. We obtain the convergence rate of our Monte Carlo method in
terms of the amounts of both basis functions and the sampled states. Finally, we demonstrate the effectiveness of our method in an optimal order execution problem with market friction and in an
inventory management problem in the presence of lost sale and lead time. Both examples are well known in the literature to be difficult to solve for optimality. The experiments show that our method
can significantly improve the heuristics suggested in the literature and obtain new policies with a satisfactory performance guarantee.
\end{abstract}


\section{Introduction}
\label{sec:intro}

Stochastic dynamic programming (SDP) provides a powerful framework for modeling and solving decision-making problems under a random environment in which uncertainty is resolved
and actions are taken sequentially over time. Recently it also has become increasingly important to help us understand the general principle behind reinforcement learning, a rapidly developing
area of artificial intelligence. The Bellman backward recursion fully characterizes the structure of the optimal policies of an SDP problem. However, hampered by the curse of dimensionality,
it is practically infeasible to implement this principle of optimality to derive the solutions for many high dimensional applications.  Hence, people often have to settle for a suboptimal control policy that
strikes a reasonable balance between convenient implementation and adequate performance. This practice naturally gives rise to the following two research questions:
\begin{enumerate}
\item How can we assess the quality of a given control policy?
\item If we know the performance of a policy is not satisfactory, do we have a systematic way to improve it?
\end{enumerate}

Motivated by these two questions, especially the second one, we develop in this paper a duality-driven iterative approach to obtaining and improving confidence interval estimates for the
true value of an SDP problem with finite time horizon. This new approach stems from information relaxation and the corresponding dual formulation in the SDP literature.
Take a cost minimization problem as an example. Within the dual framework laid out in \cite{b}, we relax the admissible constraint that requires policies to be dependent only upon the
information up to the moment when a decision is made, and meanwhile impose a penalty in the problem's objective function that punishes any violations of the admissible constraint.
This two-step construction results in a lower bound on the optimal expected cost.

The above duality bounds enable us to assess the performance of a candidate policy. Fixing the policy we are interested in assessing, we can use standard simulation techniques to estimate the
expected costs under this policy (refer to, for example, \cite{p} for other related statistical learning approaches for policy evaluation). Note that every policy is suboptimal and thus produces a
value higher than the optimal cost. If the difference, referred to as the duality gap hereafter, between the expected value of this policy and the aforementioned lower bound from the dual
formulation is tight, we can assert that the policy must be close to the optimality. A variety of applications of this duality based policy assessment can be found, just to name a few, in
\cite{lms} and \cite{lai} for natural gas storage valuation, \cite{bs2}, \cite{hw}, and \cite{hiw} for dynamic portfolio investment, \cite{b} and \cite{bs} for inventory management, \cite{got}
for multi-vehicle routing, \cite{bs} for revenue management, \cite{kl} for robust multi-armed bandits, \cite{das} for an integrated optimization problem of procurement, processing, and trade of commodities, \cite{bbc} for stochastic scheduling problems, \cite{bb} for stochastic knapsack problems,
stochastic scheduling on parallel machines, and sequential search problems, and most recently, \cite{bs3} for dynamic selection problems.

Complementing the applications of the SDP duality in policy assessments, the primary focus of our work is how to improve a candidate policy if we find that its duality gap is not small.
The paper makes two contributions to the literature on SDP duality. First, we propose a new duality-driven dynamic programming (DDP) algorithm that is capable of
iteratively improving the estimates of the dual value to an SDP problem. In each iteration, the algorithm utilizes the dual values from the last iteration as inputs to construct the penalty and
then outputs new dual values for the next round. We manage to show that the sequence of lower bound estimates that result from the proposed algorithm monotonically converges from below
to the true value function of an SDP problem with a cost minimization objective. More importantly, for problems with a finite time horizon, we also prove that such convergence will be accomplished
in a finite number of dual iterations and the optimal control can thereby be obtained on the basis of the dual value function that is output at the termination of the DDP algorithm. With these
important theoretical underpinnings, the new algorithm systematizes the improvement of a policy with a large duality gap, which addresses the second issue imposed at the beginning of the paper
that remains largely unanswered in the SDP duality literature. We demonstrate this convergence result by applying the DDP algorithm to the linear-quadratic control (LQC) problem, one of the most
fundamental problems in control theory. Corroborating the above theoretical discovery, the calculation reveals that, from a suboptimal policy, our DDP algorithm can yield the optimal linear policy within
just two dual iterations.

The second contribution of this paper is that we present a high-dimensional numerical implementation approach for DDP and develop its related performance guarantee. To overcome the curse
of dimensionality in the high-dimensional setup, we combine the regression architecture with Monte Carlo simulation to extrapolate the dual estimates observed on the sampled states to the entire state
space for approximating dual functions in each iteration of the DDP algorithm. The dual bound yielded from this algorithm can help us build up effective confidence interval estimates on the
value of the SDP problem, from which we can determine the optimality of the improved policy. Though the approach shares some common features with the existing simulation and approximation
methods in the study of approximate dynamic programming (see, e.g., \cite{bt}, \cite{ls}, \cite{tv1,tv2}, \cite{p}), the special structure of the dual formulation distinguishes it from the others in several
key aspects:
\begin{itemize}
\item Compared with the Monte Carlo duality in American option pricing (see, e.g., \cite{r1}, \cite{hk}, \cite{ab}, \cite{cg}, and \cite{dfm3}), one additional layer of complexity in dealing with a general
dynamic program is that the policies taken by the decision maker will affect the evolution of the underlying system. This leads us to face the challenging tradeoff between exploration and exploitation
when we try to numerically implement the DDP algorithm; see the counterexample in Appendix \ref{app:ee}. To avoid the exploration pitfall, we introduce a device called a state sampler
into our Monte Carlo approach and analyze its role in determining the convergence of the method.
\item To determine the dual value in each iteration, the DDP algorithm requires solving an optimization problem before taking expectation. Along one sample path of randomness, such an optimization
problem is deterministic. This salient characteristic is in stark contrast to the classical value iteration algorithm widely used in dynamic programming where one has to solve stochastic programs to
optimize an expected value. As shown in the discussion on the LQC problem (Sec. \ref{sec:lqc}) and the numerical examples (Sec. \ref{sec:num}), the vast research base of deterministic
optimization enables us to have a high degree of flexibility in choosing effective numerical procedures for our DDP algorithm.
\item Another advantage of solving optimization inside expectation is that it allows us to deploy parallel computing to accelerate the execution of the DDP algorithm. In particular, we can
simulate different groups of sample paths in parallel processors and solve the corresponding optimization programs simultaneously; then we can take the average across all the outcomes
collected from the parallel processors to compute the dual values. The parallelization grants scalability to the DDP algorithm.
\end{itemize}

To develop a performance guarantee for the above regression-based simulation approach, we characterize its rate of convergence to the true value in terms of the amounts of both basis
functions for the purpose of function approximation and the sampled states on which the dual values are estimated. Our analysis reveals an intriguing trade-off between model complexity and
simulation efforts. More specifically, the number of sampled states should be proportionally sufficient relative to the number of basis functions; otherwise, the effect of model overfitting may cause the
outcome from the DDP algorithm to diverge, rather than converge, even if both amounts tend to infinity. The paper quantifies a relative growth order between the numbers of the sampled states
and basis functions as a sufficient condition to warrant the convergence.

We demonstrate the effectiveness of our DDP algorithm with two numerical examples. One is about portfolio execution (a variant of \cite{bl}) and the other is about inventory management
(\cite{zpa,zpb}). Both examples are widely known in the literature to be intractable due to the constraints imposed on the policies and the complex high-dimensional dynamics. Using the above
DDP algorithm, we significantly improve a variety of conventional heuristics suggested in the literature, such as lookahead and linear programming approximation, to yield new policies with
satisfactory performance. It is worthwhile mentioning that, aiming at the convex structure in these examples, we apply difference-of-convex (DC) programming to solve the inner
optimization problem in their dual formulation. The tightness of the resulted confidence intervals strongly indicates this programming technique works very effectively for convex control problems.

As noted earlier, the paper extends and complements the literature on information relaxation and SDP dualities initiated by \cite{b}. Along this research line, \cite{bs} consider dynamic programs that
have a convex structure and use the first-order linear approximations of value functions to construct gradient penalties that can provide tight bounds. \cite{bh} and \cite{yz} generalize the information
relaxation approach for calculating performance bounds for infinite horizon Markov decision processes and continuous-time controls, respectively.  \cite{dfm} compare the duality in the perfect
information relaxation (called martingale duality in their paper) with the approximate linear programming approach in the literature (e.g., \cite{ss}, \cite{dv,dv2}). They find that the former one
can produce tighter lower bounds on the optimal cost-to-go function of a Markov decision problem. More recently, \cite{hr} derive the information relaxation bounds to Markov decision processes
with partial observations.

To the best of our knowledge, the idea of information relaxation based duality can be dated back to \cite{rw}, who show the possibility of associating with the non-anticipative
requirement on the solution of a multi-stage stochastic program a Lagrange multiplier that satisfies a martingale property. \cite{d1,d2} and \cite{dz} also find that introducing appropriate Lagrange
multiplier terms in the objective function of an LQC problem and solving the corresponding pathwise optimization problem will lead to the optimal controls for the original problem. Later, \cite{r2}
represents the value function of a discrete-time controlled Markov process in a dual Lagrangian form with the help of measure-change arguments and the perfect information relaxations.

A lot of interesting theoretical results, such as weak and strong dualities under various setups, have been established by the aforementioned papers. People especially find that the
dual value should be identical to the true value of the original SDP problem for an optimally chosen penalty --- the strong duality relation. However, solving for this optimal penalty is not easy. Thus
the existing literature typically heuristically selects ``good" martingale penalty functions and numerically examine its quality. Contributing to this literature, the DDP algorithm presents a systematic
approach to iteratively construct the optimal duality.

As a special case of the general SDP problem,  \cite{r1}, \cite{hk}, \cite{ab}, and \cite{dfm3} investigate the dual representation of American option pricing and more generally the optimal
stopping problem. In particular, \cite{cg} discuss how to improve the dual bounds on the option prices iteratively. However, what differentiates the case of American option pricing or more
broadly optimal stopping from a general SDP problem is that the state transition probabilities in the former case generally do not depend on the exercising actions taken by the option holder.
In this sense, our paper extends the study of \cite{cg} to a general setup of dynamic programming.

The remainder of the paper is organized as follows. In Section \ref{sec:dual}, we review the basic duality results developed by \cite{b}. We develop the theory underpinning the DDP algorithm
in Section \ref{sec:subsolution} and illustrate how it works using the LQC problem as an example. Section \ref{sec:mc} is devoted to the regression-based Monte Carlo simulation implementation and
the related convergence analysis. Section \ref{sec:num} presents two numerical experiments. All the proofs and some supplementary discussions are deferred to the AppendixAppendix.

\section{The Dual Formulation of an SDP Problem}
\label{sec:dual}

To fix the idea, we consider a generic finite-horizon discrete-time SDP problem in a
probability space $(\Omega, \mathcal{F}, \mathbb{P})$. Suppose that a planner makes sequential control decisions on a system
over a $T$-period time horizon indexed by $t=0,1,...,T$. At the beginning of each time period $t$,  given the system state
$x_t \in \mathbb{R}^n$, she takes an action $a_t \in A_{t} \subseteq \mathbb{R}^{m}$, where $A_{t}$ is the set of all feasible actions
at that moment. A random vector $\xi_{t}: \Omega \rightarrow \mathbb{R}^d$ will materialize during the period. To make the problem Markovian,
we assume that all $\xi_{t}$'s are independent. The purpose of this assumption is only for notational simplicity. Most of the subsequent results still
hold when we generalize the discussion to non-Markovian cases in which the probability distribution of $\xi_{t}$ may depend on the whole trajectories of
$\{\xi_0, ..., \xi_{t-1}\}$ and $\{x_{0}, ..., x_{t}\}$. The planner then incurs a cost amounting to $r_{t}$ that may be dependent on $x_{t}$, $a_{t}$, and $\xi_{t}$. The system
evolves to a new state according to the following recursive dynamic
\begin{eqnarray}
\label{dynamic}
x_{t+1}=f_t(x_t,a_t,\xi_t)
\end{eqnarray}
and the next round of decision making starts. Here $f_{t}$, $t=0, 1, ..., T-1$, is a function from $\mathbb{R}^n \times A_{t} \times \mathbb{R}^d$
to $\mathbb{R}^n$, mapping the current state, the selected action, and the realized randomness to another state.
%Figure \ref{fig:timeline} depicts
%the timeline of decision making in a representative period.
%\begin{figure}[htp]
%\centering
%\includegraphics[width=2.5in]{timeline.pdf}
%\caption{Timeline of the system control for a representative period.}
%\label{fig:timeline}
%\end{figure}
The planner attempts to minimize the expected aggregate costs
\begin{eqnarray}
\label{reward}
\mathbb{E}\left[\sum_{t=0}^{T-1}r_{t}(x_t, a_{t}, \xi_{t}) + r_{T}(x_{T})\Big| x_{0}\right]
\end{eqnarray}
in this process by taking proper actions, where $r_{T}(x_{T})$ stands for the terminal cost received at the end of the planning horizon.

We call $\alpha=(\alpha_{0}, \dots, \alpha_{T-1})$ a \emph{policy} if each argument $\alpha_{t}$ of it is a function from $\Omega$ to $A_{t}$, $t=0, \cdots, T-1$. In other words,
a policy prescribes the rule of action selection for the planner for each possible outcome $\omega$ in $\Omega$ in each period. To reflect the information constraint that the planner
faces, assume that she cannot peek into the future of the system dynamics. Hence, the decision that she makes in period $t$ relies only on what is known about the past trajectory
of the system at the beginning of the period. More formally, letting $\mathcal{F}_t=\sigma(x_{0}, ..., x_{t})$ be the $\sigma$-algebra generated by the information about the system
states up to time $t$, we require the planner's policy to be \emph{admissible} in the sense that $\alpha_{t}$ is $\mathcal{F}_t$-measurable for all $0 \le t \le T-1$.
Denote $\mathbb{F}=(\mathcal{F}_{0}, \mathcal{F}_{1}, \dots, \mathcal{F}_{T-1})$ with $\mathcal{F}_{0}=\{\emptyset, \Omega\}$. The objective of the decision maker can then be
formulated as optimizing
\begin{eqnarray}
\label{control}
V_{0}(x)=\inf_{\alpha \in \mathcal{A}_{\mathbb{F}}}\mathbb{E}
\left[\sum_{t=0}^{T-1}r_{t}(x_t, \alpha_{t}, \xi_{t}) + r_{T}(x_{T})\Big| x_{0}=x\right],
\end{eqnarray}
where $\mathcal{A}_{\mathbb{F}}$ denotes the collection of all admissible policies with respect to the information filtration $\mathbb{F}$.

It is well known that we may invoke the \emph{principle of dynamic programming} (or the Bellman equation) to solve the above SDP problem (\ref{control}). Let $V_t(x)$ be the cost-to-go function of the system from time $t$ onward; that is,
\begin{eqnarray}
\label{value-to-go}
V_t(x)=\inf_{\alpha \in
\mathcal{A}_{\mathbb{F}}|t}\mathbb{E}
\left[ \sum_{s=t}^{T-1}r_{s}(x_s, \alpha_{s}, \xi_{s}) + r_{T}(x_{T})\Big| x_t=x \right],
\end{eqnarray}
where
\[
\mathcal{A}_{\mathbb{F}}|t=\Big\{\alpha=(\alpha_{t}, \dots, \alpha_{T-1}):\ \alpha_{s}\ \textrm{is $\mathcal F_{s}$-measurable for all $t \le s \le T-1$}\Big\}.
\]
The Bellman equation dictates that we can determine the value of $V_{t}$ in a backward fashion:
\begin{eqnarray}
\label{bellman1}
V_{T}(x)&=&r_{T}(x);\\
\label{bellman2}
V_t(x)&=&\inf_{a_t \in A_t}\mathbb{E}\left[r_t(x, a_{t}, \xi_{t})+V_{t+1}(f_t(x,a_t,\xi_t))\right]
\end{eqnarray}
for all $t=0, \cdots, T-1$ and $x \in \mathbb{R}^n$. The expectation in (\ref{bellman2}) is taken with respect to the probability distribution of $\xi_{t}$. Furthermore, if
$a^{*}_{t}=\alpha^{*}_{t}(x)$ minimizes the right hand side of (\ref{bellman2}) for each $x$ and $t$, the policy $\alpha^*=(\alpha^*_{0}, \dots, \alpha^*_{T-1})$ is optimal.

However, the curse of dimensionality prevents us from directly utilizing the Bellman equations (\ref{bellman1}-\ref{bellman2}) to solve the SDP problem because the
computational complexity that this procedure incurs grows exponentially as the dimensionality of the state, randomness, and action spaces increase; see, e.g., Sections
1.2 and 4.1 in \cite{p} for detailed discussions on this issue. In light of this difficulty, people often have to settle for a computationally tractable approximate (thus, suboptimal)
policy of adequate performance. This gives rise to a natural question about how to assess such approximate policies without knowing where the optimality is.
As noted in the introduction, the dual formulation proposed in \cite{b} presents a systematic approach by which we can measure the quality of a suboptimal policy, or in other words,
how close it is to the optimal one.

The key ingredients of their duality are the concept of \emph{information relaxation} and a related \emph{penalty}. For the purpose of this paper, we only consider the case of
perfect relaxation and refer readers to their paper for a rigorous development of the dual theory under a general framework. Intuitively, if we relax the requirement of information
admissibility on policies by allowing the decision maker to take actions after she observes the entire realization of randomness $(\xi_{1},\cdots, \xi_{T})$, we should be able to obtain
a lower bound to the true cost value $V_{0}$. More precisely, by Jensen's inequality, we have
\begin{eqnarray}
\label{explain}
\mathbb{E}\left[\inf_{a \in A}\left(\sum_{s=0}^{T-1}r_{s}(x_s, a_{s}, \xi_{s}) + r_{T}(x_{T})\right)\Big| x_0=x\right] \le V_{0}(x)
\end{eqnarray}
for all $x$. Note that the minimizer of the optimization inside the expectation on the left hand side of (\ref{explain}) is not admissible in the original problem because
it may depend on the whole trajectory of $(\xi_{1},\cdots, \xi_{T})$.

\cite{b} further points out that we can achieve equality in (\ref{explain}) if properly penalizing the objective function inside the expectation. Corresponding to the
above perfect relaxation, one possible penalty can be constructed as follows. Let $W=(W_{1}(\cdot), \dots, W_{T}(\cdot))$ be any sequence of functions such that each argument
$W_{t}: \mathbb{R}^n \rightarrow \mathbb{R}$ maps the system state to real numbers. Given an action sequence $a=(a_{0}, \dots, a_{T-1}) \in A:=A_{0}\times \cdots \times A_{T-1}$
and a sequence of randomness $\xi=(\xi_{0}, \dots, \xi_{T-1})$, we can use Eq. (\ref{dynamic}) to recursively generate a trajectory of system states $(x_{1},\cdots, x_{T})$.
Along it, define a penalty function such as
\begin{eqnarray}
\label{penalty_def}
z(a, \xi)=\sum_{t=0}^{T-1}\left\{\mathbb{E}[r_{t}(x_{t}, a_{t}, \xi_{t})+W_{t+1}(f_{t}(x_{t}, a_{t}, \xi_{t}))]-(r_{t}(x_{t}, a_{t}, \xi_{t})+W_{t+1}(f_{t}(x_{t}, a_{t}, \xi_{t})))\right\},
\end{eqnarray}
where the expectation inside the sum is taken with respect to the distribution of $\xi_{t}$. Then, \cite{b} show that
\begin{eqnarray}
\label{dual} V_{0}(x)=\sup_{W}\mathbb{E}\left[ \inf_{a \in A}\left(\sum_{s=0}^{T-1}r_{s}(x_s, a_{s}, \xi_{s}) + r_{T}(x_{T})+z(a, \xi)\right)\Big| x_{0}=x\right].
\end{eqnarray}

The strong duality relationship (\ref{dual}) paves a useful way to assessing the quality of a specific admissible policy $\alpha$. First, we may evaluate the policy by calculating
\[
\overline{V}_{t}(x)=\mathbb{E}\left[ \sum_{s=t}^{T-1}r_{s}(x_s, \alpha_{s}, \xi_{s}) + r_{T}(x_{T})\Big| x_t=x \right]\ \textrm{for all $0 \le t \le T$}.
\]
Surely $\overline{V}_{t}(x) \ge V_{t}(x)$ for any $x$ because of the sub-optimality of $\alpha$. Then, we replace the generic $W$ in (\ref{penalty_def}) by $\overline{V}_{t}$ to
construct a penalty $z$ and compute the associated dual value
\[
\underline{V}_{0}(x)=\mathbb{E}\left[ \inf_{a \in A}\left(\sum_{s=0}^{T-1}r_{s}(x_s, a_{s}, \xi_{s}) + r_{T}(x_{T})+z(a, \xi)\right)\Big| x_{0}=x\right].
\]
From (\ref{dual}), we have $V_{0}(x) \ge \underline{V}_{0}(x)$, which implies
\[
0 \le \overline{V}_{0}(x)-V_{0}(x) \le \overline{V}_{0}(x)-\underline{V}_{0}(x).
\]
When the dual gap  $\overline{V}_{0}-\underline{V}_{0}$ is sufficiently tight, we can conclude that the performance of policy $\alpha$ must be very close to the optimality.
One can refer to those works mentioned in the introduction for various applications of the above duality-based policy assessment.

\section{DDP: A Duality-Driven Dynamic Programming Method}
\label{sec:subsolution}

Beyond the aforementioned policy assessment, the primary interest of the current paper is on the second research question posed in the introduction: can we develop a systematic
approach to improving the policy in hand if we find that its dual gap is not tight enough? In this section, we build up an iterative method on the basis of the SDP information duality
to achieve the goal of policy improvement.

\subsection{Subsolutions and Dual Value Iteration}
\label{sec:value_iteration}

Central to our investigation are the notion of \emph{subsolution} and, more importantly, its close relationship with the information duality.
\begin{dfn}[subsolution]
\label{dfn:subsolution}
A functional sequence $S=(S_0, S_1,..., S_{T})$ with $S_{t}: \mathbb{R}^{n} \rightarrow \mathbb{R}$, $0 \le t \le T$, is called a subsolution to the problem (\ref{control}) if it satisfies
\begin{eqnarray*}
S_t(x) \le \inf_{a_t \in A_t}\mathbb{E}\left[r_t(x, a_{t}, \xi_{t})+S_{t+1}(f_t(x,a_t,\xi_t))\right]
\end{eqnarray*}
for any $t=0,1,...,T-1$ and $x \in \mathbb{R}^{n}$ with the convention that $S_{T}(x) = r_{T}(x)$.
\end{dfn}
The concept of subsolutions to a generic SDP problem has been long known in the literature; one may see, for instance, Theorem 6.2.2 in \cite{pu} or Theorem 3.4.1 in \cite{p}.
It just generalizes the Bellman equation (cf. (\ref{bellman2})) by replacing the equality with an inequality. One well-known fact is that any subsolution provides a lower bound on
the true value of the primal problem (\ref{control}) (e.g., Theorem 6.2.2 in \cite{pu}). Using the subsolution requirement on each state as the constraints, \cite{dv} developed a
linear programming based approach to approximate solutions to the SDPs. Let $\mathcal{S}$ denote the collection of all the subsolutions to the problem (\ref{control}).

As one of the key underpinnings of our DDP algorithm, Proposition \ref{pro:super} points out that the dual operation actually offers us a way to construct subsolutions.
Introducing some operator notations here will help us present the main results in a compact way. Take any functional sequence $W=(W_{0}(\cdot), \dots, W_{T}(\cdot))$ and consider
the tail subproblem (\ref{value-to-go}) for each $t$, $0 \le t \le T-1$. Note that it is still an SDP problem. Hence, we can apply the corresponding dual formulation to it, namely,
construct the associated penalty
\begin{eqnarray}
\label{sub_penalty}
\resizebox{.95\hsize}{!}{$
z_{t}(a, \xi)=\sum_{s=t}^{T-1}\left\{\mathbb{E}[r_{s}(x_{s}, a_{s}, \xi_{s})+W_{s+1}(f_{s}(x_{s}, a_{s}, \xi_{s}))]-(r_{s}(x_{s}, a_{s}, \xi_{s})+W_{s+1}(f_{s}(x_{s}, a_{s}, \xi_{s})))\right\}, $}
\end{eqnarray}
by using the tail sequence of $W$, $(W_{t+1}(\cdot), \dots, W_{T}(\cdot))$, and obtain the dual function
\begin{eqnarray}
\label{sub_weak_dual}
W'_{t}(x):=\mathbb{E}\left[ \inf_{a \in A|t}\left(\sum_{s=t}^{T-1}r_{s}(x_s, a_{s}, \xi_{s}) + r_{T}(x_{T})+z_{t}(a, \xi)\right)\Big| x_{t}=x\right]
\end{eqnarray}
for each $t$, where $A|t=A_{t} \times \cdots \times A_{T-1}$. In this way, as implied by the duality theory discussed in the last section, we reach a sequence of lower bounds
$W'=(W'_{0}(\cdot), \cdots, W'_{T}(\cdot))$ to the true cost value of every tail problem. From now on, let $\mathcal{D}$ denote the dual operator defined through
(\ref{sub_penalty}-\ref{sub_weak_dual}) that can be viewed as acting on any functional sequence $W$ to produce another function sequence $\mathcal{D}W=((\mathcal{D}W)_{0}, \cdots,
(\mathcal{D}W)_{T})$, where $(\mathcal{D}W)_{t}(x)=W'_{t}(x)$ for $0 \le t \le T-1$ and $(\mathcal{D}W)_{T}(x)=r_{T}(x)$.

Examining the relationship among $(\mathcal{D}W)_{t}$ across all $t$'s, we have
\begin{pro}
\label{pro:super}
Let $W=(W_0,W_1,...,W_T)$ be any functional sequence. Then,
$\mathcal{D}W \in \mathcal{S}$, i.e., for all $t$,
\begin{eqnarray*}
\mathcal{D}W_t(x) \le \inf_{a_t \in A_t}\mathbb{E}\left[r_t(x, a_{t}, \xi_{t})+\mathcal{D}W_{t+1}(f_t(x,a_t,\xi_t))\right].
\end{eqnarray*}
\end{pro}
This proposition reveals that the information relaxation based duality and the subsolutions are closely related. In particular, the former presents a systematic way of constructing the latter.
This finding is new to the existing literature to the best of our knowledge. Moreover, Proposition \ref{pro:super} indicates that, if we repeatedly apply the the operator $\mathcal{D}$ on $W$, i.e.,
letting $\mathcal{D}^{n}W=\mathcal{D}(\mathcal{D}^{n-1}W)$ for all $n \ge 1$, we can obtain a sequence of subsolutions $\{\mathcal{D}^nW,\ n \ge 1\}$.

Now we are ready to present Theorem \ref{thm:convergence}, one of the main results of the paper. In it, we show that the above dual value sequence increasingly converges to the true cost-to-go
function of the primal problem (\ref{control}).
\begin{thm}
\label{thm:convergence}
(i) The subsolution sequence $\{\mathcal{D}^{n}W, n \ge 1\}$ is increasing in $n$ in the sense that $(\mathcal{D}^{n+1}W)_{t}(x) \ge (\mathcal{D}^{n}W)_{t}(x)$
for all $n\ge 1$, $0 \le t \le T$, and $x \in \mathbb{R}^{n}$;\\
(ii) if, for some $n$, $(\mathcal{D}^{n+1}W)_{t}(x)=(\mathcal{D}^{n}W)_{t}(x)$ for all $t$ and $x$, then $\mathcal{D}^{n}W \equiv V$;\\
(iii) $\mathcal{D}^{T+1}W=V$.
\end{thm}
Recall that any subsolution is dominated by the true cost-to-go function. Hence, one implication of Part (i) of Theorem \ref{thm:convergence} is that $\mathcal{D}^{n}W \le \mathcal{D}^{n+1}W \le V$.
In other words, the subsolution sequence $\{\mathcal{D}^{n}W\}$ iteratively improves its quality of approximation as lower bounds on $V$. Two key facts underpin the proof of Part (i). First,
we need to show that, for any given subsolution, applying the dual operation on it will lead to a tighter lower bound on the true value function of the primal problem. It is worth noting that
similar results have been established in the setups of optimal stopping problems (\cite{cg}) and infinite-horizon Markov decision processes (\cite{dfm} and \cite{bh}). To prove Part (i),
we manage to extend the fact to a finite-horizon framework. The second fact, Proposition \ref{pro:super}, also plays an important role in the proof. It guarantees that $\mathcal{D}^{n}W$, as the output
of the last dual iteration, is still a subsolution. So, implied by the first fact, we can further apply $\mathcal{D}$ on it in the next iteration to yield more improvement. In other words, Proposition
\ref{pro:super} accomplishes the inductive step for us to carry out induction on the sequence of $\{\mathcal{D}^{n}W\}$ to show (i).

A more powerful conclusion stems from Parts (ii) and (iii) of the theorem. That is, the improvements in the sequence $\{\mathcal{D}^{n}W, n \ge 1\}$ will terminate in a finite number of iterations
and when it terminates, the optimal value of the primal problem has been achieved. From this, we propose  the following DDP algorithm to solve the problem (\ref{control}) in an
iterative manner:
\begin{framed}
{\small \begin{center}
\textbf{Table I: A Duality Driven Dynamic Programming (DDP) Algorithm}
\end{center}
\begin{itemize}[labelwidth=-0.25in,itemindent=-0.5in]
\item
\textbf{Step 0.} Initialization:
\begin{itemize}[leftmargin=0in]
\item \textbf{Step 0a.} Select an initial approximate value function sequence $W^{0}=(W^{0}_{0}, \cdots, W^{0}_{T})$. One way to do it, for instance, is to use a feasible policy $\alpha$ to compute its corresponding value
\begin{eqnarray*}
W^{0}_{t}(x):=\mathbb{E}\left[ \sum_{s=t}^{T-1}r_{s}(x_s, \alpha_{s}, \xi_{s}) + r_{T}(x_{T})\Big| x_t=x \right]
\end{eqnarray*}
for all $x \in \mathbb{R}^{n}$ and $0 \le t \le T-1$. Let $\underline{V}^{0}=W^{0}$.
\item \textbf{Step 0b.} Set $n=1$.
\end{itemize}
\item
\textbf{Step 1.} Construct subsolutions using the dual operator $\mathcal{D}$:
\begin{itemize}[leftmargin=0in]
\item \textbf{Step 1a.} For $\underline{V}^{n-1}$, define a penalty function sequence $z^{n}=(z^{n}_{0}, \cdots, z^{n}_{T})$ such that $z^{n}_{T}(a, \xi)=0$ and for any given $0 \le t \le T-1$,
{
\begin{eqnarray}
\label{penalty}
\resizebox{.85\hsize}{!}{$
z^{n}_{t}(a, \xi)=\sum_{s=t}^{T-1}\left\{\mathbb{E}[r_{s}(x_{s}, a_{s}, \xi_{s})+\underline{V}^{n-1}_{s+1}(f_{s}(x_{s}, a_{s}, \xi_{s}))]-(r_{s}(x_{s}, a_{s}, \xi_{s})+\underline{V}^{n-1}_{s+1}(f_{s}(x_{s}, a_{s}, \xi_{s})))\right\}
 $}
\end{eqnarray}
}
with $a=(a_{0}, \cdots, a_{T-1}) \in A$ and $\xi=(\xi_{0}, \cdots, \xi_{T-1})$.
\item \textbf{Step 1b.} For all state $x$ and time $t$, determine the value of the following lower bound
\begin{eqnarray}
\label{algo_1}
\underline{V}^{n}_{t}(x)=\mathbb{E}\left[ \inf_{a \in A|t}\left(\sum_{s=t}^{T-1}r_{s}(x_s, a_{s}, \xi_{s}) + r_{T}(x_{T})+z^{n}_{t}(a, \xi)\right)\Big| x_{t}=x\right].
\end{eqnarray}
\end{itemize}
\item
\textbf{Step 2.}
If $\underline{V}^{n}(x) \neq \underline{V}^{n-1}(x)$ for some $x$, let $n=n+1$ and go to Step 1.
\end{itemize}}
\end{framed}

Though the DDP algorithm focuses on updating the dual value, it can be used to improve control policies as well. For a suboptimal policy $\alpha$, we may run Step0a to evaluate it and
initiate the algorithm with its policy value. Suppose that the algorithm terminates at the $n$th iteration. Replace the value function $V_{t+1}$ in the one-step Bellman equation (\ref{bellman2})
with $\underline{V}^{n}_{t+1}$ and solve
\begin{eqnarray}
\label{policy}
\alpha^{n}_{t}(x)=\argmin_{a_{t} \in A_{t}}\mathbb{E}\left[r_t(x, a_{t}, \xi_{t})+\underline{V}^{n}_{t+1}(f_t(x,a_t,\xi_t))\right]
\end{eqnarray}
for a new policy $\alpha^{n}_{t}(\cdot)$ at time $t$, $0 \le t \le T-1$. According to Part (ii) of Theorem \ref{thm:convergence}, we should have achieved the optimality, i.e., $\underline{V}^n=V$.
The optimality of $\alpha^n=(\alpha^{n}_{0}, \cdots, \alpha^{n}_{T-1})$ ensues.

\subsection{An Illustration: Linear-Quadratic Control}
\label{sec:lqc}

Below we will use the classical LQC problem to demonstrate the effectiveness of policy improvement of the algorithm. In this case, the DDP algorithm can yield
the optimal policy after just \emph{two} iterations of the dual operation, no matter how long the time horizon of the problem is. By (\ref{algo_1}), the key steps in each dual iteration involve
solving the inner optimization problem and determining the outer expectation value. One caveat is that, unlike the LQC example in which closed-form expressions for both are available,
it is in general difficult to explicitly carry out these two types of computation, especially in high-dimensional problems. To address this issue, we shall explore in Section \ref{sec:mc} how to
resort to some numerical techniques, such as Monte Carlo simulation and the related approximation architectures, to implement the DDP algorithm effectively. One error bound is also
developed therein (cf. Theorem \ref{thm:error}) to deliver the performance guarantee. Theorem \ref{thm:convergence}, despite its theoretical nature, still serves as an
important cornerstone for us to obtain such numerical performance guarantees.

The LQC problem has received a lot of attention in control theory because of its tractability. It is widely applied in automatic control of a motion or a
process to formulate how to regulate a system to stay close to the origin. The closed-form solution to the problem is well known in the literature. The intention
of this subsection is definitely not to repeat these known results. Instead, we want to corroborate the result of the last subsection by showing its
policy-improving effect. Following the standard setup of a LQC problem, consider a system whose dynamic equation is given by
\begin{eqnarray}
\label{dynamic_LQ}
x_{t+1}=D_{t}x_t+B_{t}a_t+\xi_t,\quad t=0, \cdots, T-1.
\end{eqnarray}
When it runs, it will incur a cost of
\begin{eqnarray}
\label{cost_LQ}
\sum_{t=0}^{T-1}\left(x^{tr}_t Q_{t} x_t +a^{tr}_{t} R_{t}a_{t}\right) + x^{tr}_T Q_{T} x_T.
\end{eqnarray}
In these expressions, $D_{t} \in \mathbb{R}^{n \times n}$, $B_{t} \in \mathbb{R}^{n \times m}$, $Q_{t} \in \mathbb{R}^{n \times n}$, and $R_{t} \in \mathbb{R}^{m \times m}$,
are all given. The matrices $Q_{t}$ are positive semidefinite symmetric and the matrices $R_{t}$ are positive definite symmetric. There is no constraint on the controls $a_{t}$, i.e.,
we may take any vector in $\mathbb{R}^{m}$ as its value. Each $\xi_{t}$ has zero mean and a finite second moment. Assume that the decision maker has perfect information of the state
$x$ over the course of system evolution.

From the above description, it is not difficult to see that this problem is just a special case of (\ref{dynamic}-\ref{reward}) by taking a linear form for the evolution function $f_{t}$ and a
quadratic form for the cost $r_{t}$. Its optimal policy is explicitly known in the literature (see, e.g., Sec. 4.1 in \cite{ber}) to be of the following linear form: $\alpha^*_{t}(x)=-L_{t}x$,
for $t=0, \cdots, T-1$. Accordingly, the optimal cost function equals
\begin{eqnarray}
\label{truev_LQ}
V^*_{t}(x)=x^{tr}K_{t}x+\sum_{s=t}^{T-1}\mathbb{E}\left[\xi^{tr}_{s}K_{s+1}\xi_{s}\right].
\end{eqnarray}
Here, both matrices $L_{t} \in \mathbb{R}^{m \times n}$ and $K_{t}\in \mathbb{R}^{n \times n}$ are explicitly computable. Detailed discussions are deferred to Electronic
Companion \ref{LQ_iteration}.

Applying the DDP algorithm to the LQC problem, we have
\begin{pro}
\label{pro:lq}
Fix a matrix $P_{t} \in \mathbb{R}^{m \times n}$ and a vector $E_{t} \in \mathbb{R}^{m \times n}$ for each $t$. Consider a policy of the linear form
\begin{eqnarray}
\label{LQ_a0}
\alpha_{t}(x)=P_t x+E_t,\quad 0 \le t \le T-1.
\end{eqnarray}
If we start the DDP algorithm with this policy, then it will terminate after two iterations at $\underline{V}^{2} \equiv V^*$.
\end{pro}
Corroborating the results in Theorem \ref{thm:convergence}, Proposition \ref{pro:lq} shows that our DDP algorithm warrants the convergence to the true cost function of the LQC problem
in two iterations. There are several studies in the literature related to the application of the information relaxation technique in LQC. \cite{dz} postulate a linear form for the optimal penalty
and thereby present a new proof of the LQC optimal control theorem based on the dual formulation. \cite{hl2012} develop two types of approaches to constructing optimal penalties
for an LQC problem. However, both of their constructions require some prior knowledge on the optimal value function of LQC. Compared with these studies, our DDP algorithm provides
a more mechanical way to find the optimal penalty with little prior knowledge required.

The proof of Proposition \ref{pro:lq} is contained in Appendix \ref{LQ_iteration}. This example highlights one advantage of working with the
duality-driven method in the computational aspect. That is, to compute the dual value, we just need to solve a deterministic optimization problem inside the expectation for which
there is a vast research base that we can draw on for help. In particular, the proof of Proposition \ref{pro:lq} shows that the minimization problem leading to the duality for the LQC problem
turns out to be a quadratic program, which is well known to be tractable in the optimization literature (see, e.g., \cite{n&w}, Chapter 16).

\section{Monte Carlo Implementation of the DDP Algorithm}
\label{sec:mc}

As noted at the beginning of Section \ref{sec:lqc}, the intrinsic difficulty of dealing with a general SDP problem lies in the fact that the inner optimization and the
outer conditional expectation in (\ref{algo_1}) often cannot be analytically solved. Below we propose the use of regression to estimate the duality $\underline{V}^{n}$ from
simulated states for the purpose of implementing the DDP algorithm via Monte Carlo simulation. A related convergence analysis is developed in Section \ref{sec:convergence}.

\subsection{Regression-based Algorithm}
\label{sec:regression}

In the first step of the algorithm, we need to generate a group of states on which the value of
\begin{eqnarray}
\label{inner}
 \inf_{a \in A|t}\left(\sum_{s=t}^{T-1}r_{s}(x_s, a_{s}, \xi_{s}) + r_{T}(x_{T})+z^{n}_{t}(a, \xi)\right)
\end{eqnarray}
will be estimated so that we can use regression to build up the approximation to the conditional expectation in (\ref{algo_1}). Many of the regression-based methods in
the literature on American option pricing (see, e.g., \cite{c}, \cite{ls}, \cite{tv1,tv2}) directly invoke the dynamic of the underlying asset to simulate states for continuation value
estimation. Note that the exercising decision for an American option has no impact on the underlying price dynamics. One additional layer of complexity encountered here in
a general SDP problem is that its state evolution hinges on the policy that we are using. As illustrated by the example in Appendix \ref{app:ee}, using a suboptimal
policy of the original problem (\ref{dynamic}) to generate the states that our DDP algorithm will visit later can possibly lead to being stuck in suboptimality, because with this policy
the algorithm may have no chance to access such states that contain useful information for us to improve the estimation.

To avoid this exploration pitfall, we suggest that the sequence of probability density functions $\{G_{1}, \cdots, G_{T}\}$ that are utilized for the purpose of state selection should
be independent of the current policy of the SDP problem. In particular, if the support sets of all the $G$'s contain the entire state space of the problem, these sampling
distributions enable us to reach any states in the space with nonzero chance. Imposing this ergodic requirement on $G$ as one of the sufficient conditions, we investigate in
Theorem \ref{thm:error} the asymptotic properties of the regression-based implementation of the DDP algorithm. Denote the number of simulated representative states by $L$
 hereafter. We independently draw $L$ states, $(x^{(1)}_{t}, \cdots, x^{(L)}_{t})$, from the distribution $G_{t}$ at each time period $t$, $t=1, \cdots, T$, where the superscript $(l)$,
 $1 \le l \le L$, indicates the $l$th sample at time $t$. The values of the dual functions will be estimated on these points.

We now turn to present the core step of the implementation, i.e., how to use Monte Carlo regression to obtain an approximation to $\underline{V}^{n}$ from the previous estimate of
$\underline{V}^{n-1}$ (cf. Step 1 in Table I). Let $\{\psi_{1}, \cdots, \psi_{M}\}$ denote a pre-specified set of basis functions, where each argument $\psi_m$ is a function mapping from
$\mathbb{R}^n$ to $\mathbb{R}$. Assume that the previous iteration has yielded that $\underline{V}^{n-1}$ can be approximated by
\begin{eqnarray}
\label{approx}
\underline{V}^{n-1}_{t}(x) \approx \underline{\widehat{\mathfrak{V}}}^{n-1}_{t}(x):=\sum_{m=1}^{M}\widehat{\beta}^{n-1}_{t, m}\psi_{m}(x)
\end{eqnarray}
for some constants $\widehat{\beta}^{n-1}_{t, m}$, $1 \le m \le M$ and $1 \le t \le T$. Following Step 1a in DDP, to construct a new penalty for the next round, we substitute the right hand side of (\ref{approx})
into (\ref{penalty}).  We then have the following approximation to $z^{n}_{t}(a, \xi)$:
\begin{eqnarray*}
&&{\mathfrak{z}}^{n}_{t}(a, \xi)\\
&&\resizebox{0.95\hsize}{!}{$ =\sum_{s=t}^{T-1}\left\{\mathbb{E}\left[r_{s}(x_{s}, a_{s}, \xi_{s})+\sum_{m=1}^{M}\widehat{\beta}^{n-1}_{s+1, m}\psi_{m}(f_{s}(x_{s}, a_{s}, \xi_{s}))\right]-\left(r_{s}(x_{s}, a_{s}, \xi_{s})
+\sum_{m=1}^{M}\widehat{\beta}^{n-1}_{s+1, m}\psi_{m}(f_{s}(x_{s}, a_{s}, \xi_{s}))\right)\right\} $}
\end{eqnarray*}
for any $a$ and $\xi$, where the expectation in the first term of ${\mathfrak{z}}^{n}_{t}$ is taken with respect to $\xi_s$.

In evaluating ${\mathfrak{z}}^{n}_{t}$, we need to compute $\mathbb{E}[\psi_{m}(f_{s}(x_{s}, a_{s}, \xi_{s}))]$. For many applications, especially when $\psi_{m}$ is a polynomial, $f_{s}$ is simple,
and the distribution of $\xi_{s}$ is analytically known, we can explicitly compute this expectation. For the cases in which its closed-form expression is not available, we may rely on
Monte Carlo simulation to generate samples from the distribution of $\xi_{s}$ and then use sample averages to approximately evaluate it. To expedite the computation in this step, we also
attempt an alternative simulation method, which is the low-discrepancy method from the quasi-Monte Carlo (QMC) literature, in the numerical experiments of the next section. Different from plain
Monte Carlo, this QMC approach deterministically chooses representative points for $\xi$. We find that QMC can deliver excellent approximation performance with a relatively smaller
number of simulation trials, consistent with the well known fact that the QMC converges faster than the ordinary Monte Carlo. One may refer to Chapter 5 of \cite{gl} for a
comprehensive coverage of this subject.

Once the value of ${\mathfrak{z}}^{n}_{t}(a, \xi)$ is determined, we proceed to build up the regression estimators for the conditional expectation (\ref{algo_1}) in Step 1b of the DDP algorithm.
To this end, we posit that (\ref{algo_1}) can be represented as a linear combination of the basis functions, i.e.,
\begin{eqnarray}
\label{reg}
\resizebox{.95\hsize}{!}{$
\mathbb{E}\left[J_{t,n}(\xi | t, x_{t})\right]:=\mathbb{E}\left[ \inf_{a \in A|t}\left(\sum_{s=t}^{T-1}r_{s}(x_s, a_{s}, \xi_{s}) + r_{T}(x_{T})+z^{n}_{t}(a, \xi)\right)\Big|x_{t}=x\right]=\sum_{m=1}^{M}\beta^{n}_{t, m}\psi_{m}(x), $}
\end{eqnarray}
at $x_{t}=x$ for any given $t$, $t=0, \cdots, T-1$. Here, $J_{t,n}(\xi | t, x_{t})$ is a shorthand notation for the minimization problem inside the expectation, whose value apparently depends on
the tail vector of random perturbation $\xi | t=(\xi_{t}, \cdots, \xi_{T-1})$ and the system state $x_{t}$ at time $t$. The standard least square arguments imply that the coefficient vector
$\beta^{n}_{t}=(\beta^{n}_{t, 1}, \cdots, \beta^{n}_{t, M})^{tr}$ in (\ref{reg})  should be given by
\begin{eqnarray}
\label{reg_theory}
\beta^{n}_{t}=(B^{t}_{\psi\psi})^{-1}B^{t, n}_{J\psi}:=(\mathbb{E}^{G}[\Psi_M(X_{t})\Psi_M(X_{t})^{tr}])^{-1}\mathbb{E}^{G\otimes \xi}[\Psi_M(X_{t})J_{t,n}(\xi | t, X_{t})].
\end{eqnarray}
In (\ref{reg_theory}), $B^{t}_{\psi\psi}$ is the indicated $M \times M$ matrix $\mathbb{E}^{G}[\Psi_M(X_{t})\Psi_M(X_{t})^{tr}]$ (assumed nonsingular) with
$\Psi_M(x)=(\psi_1(x),\cdots,\psi_M(x))^{tr}$. The superscript $G$ stresses that the expectation is defined on $G_{t}$, the distribution of $X_{t}$. Meanwhile,
$B^{t, n}_{J\psi}$ is the indicated vector of dimension $M$ computed from $\mathbb{E}^{G\otimes \xi}[\Psi(X_{t})\mathfrak{J}_{t,n}(\xi | t, X_{t})]$ with
$X_{t} \sim G_{t}$ and $\xi |t$ independently drawn from its own distribution.

Both $B^{t}_{\psi\psi}$ and $B^{t, n}_{J\psi}$ can be estimated on the basis of observations of pairs $(\xi | t, X_{t})$. More explicitly, starting from each point $x^{(l)}_{t}$, we independently
simulate one path of $\xi^{(l)}|t=(\xi^{(l), t}_{t}, \xi^{(l), t}_{t+1}, \cdots, \xi^{(l), t}_{T-1})$ from the distribution of $\xi$. Suppose for a moment that the value of $J_{t,n}(\xi | t, X_{t})$ can be
(approximately) computed at each pair $(\xi^{(l)}|t, x^{(l)}_{t})$ and denote that quantity by $\mathfrak{J}^{(l)}_{t, n}$. Let $\hat{B}^{t}_{\psi\psi}$ be an $M \times M$ matrix with the $(i, j)$-entry
\begin{eqnarray}
\label{bij}
 \frac{1}{L}\sum_{l=1}^{L}\psi_{i}(x^{(l)}_{t})\psi_{j}(x^{(l)}_{t})
\end{eqnarray}
and $\widehat{B}^{t, n}_{\mathfrak{J}\psi}$ be an $M$-vector with the $k$th entry
\begin{eqnarray}
\label{bj}
\frac{1}{L}\sum_{l=1}^{L}\mathfrak{J}^{(l)}_{t, n}\psi_{k}(x^{(l)}_{t,n}).
\end{eqnarray}
Then, an estimate to $\beta^{n}_{t}$ can be formed by $\hat{\beta}^{n}_{t}=(\hat{B}^{t}_{\psi\psi})^{-1}\hat{B}^{t, n}_{\mathfrak{J}\psi}$.  From it, we complete one iteration in our DDP algorithm
by building up a new approximate to the dual value function
\begin{eqnarray*}
\underline{V}^{n}_{t}(x) \approx \underline{\widehat{\mathfrak{V}}}^{n}_{t}(x):=\sum_{m=1}^{M}\widehat{\beta}^{n}_{t, m}\psi_{m}(x),
\end{eqnarray*}
where $\widehat{\beta}^{n}_{t, m}$ is the $m$th entry of $\hat{\beta}^{n}_{t}$.

To determine the value of $\mathfrak{J}^{(l)}_{t, n}$, we replace $z^{n}_{t}(a, \xi)$ by its approximation ${\mathfrak{z}}^{n}_{t}(a, \xi)$ in $J_{t,n}(\xi | t, X_{t})$ and solve the following
optimization problem:
\begin{eqnarray}
\label{j}
\mathfrak{J}^{(l)}_{t, n}:= \inf_{a \in A|t}\left(\sum_{s=t}^{T-1}r_{s}(x_s, a_{s}, \xi^{(l), t}_{s}) + r_{T}(x_{T})+\mathfrak{z}^{n}_{t}(a, \xi^{(l)}|t)\right)
 \end{eqnarray}
subject to the constraints $x_t=x$ and
\begin{eqnarray}
\label{constraint}
x_s=f_{s-1}(x_{s-1}, a_{s-1},\xi^{(l), t}_{s-1})\
\end{eqnarray}
for all $s=t+1, \cdots, T$. The outcome of the above optimization problem, denoted by $\mathfrak{J}^{(l)}_{t, n}$, will be used as one observation of $J_{t,n}(\xi | t, X_{t})$ at $(\xi^{(l)}|t, x^{(l)}_{t})$
to estimate $B^{t, n}_{J\psi}$.

It is worth mentioning that, given $\xi^{(l)}|t=(\xi^{(l), t}_{t}, \xi^{(l), t}_{t+1}, \cdots, \xi^{(l), t}_{T-1})$, the problem (\ref{j}-\ref{constraint}) is indeed a deterministic optimization
program. Compared with many of the SDP algorithms in which stochastic optimization is involved, the computation for the solution to (\ref{j}-\ref{constraint}) is less
demanding. Similar to the case of LQC problems, the vast research literature on deterministic optimization provides us various flexible and potent methodologies that we can draw on to
solve it. In particular, we develop in the next section an efficient numerical scheme based on the DC programming to solve this inner optimization problem for a broad class of control
problems. Simplifying the underlying probabilistic structure of an SDP problem to yield some computational advantages is a commonly used
strategy in approximate dynamic programming. For instance, the approach of certainty equivalent control replaces the stochastic disturbances with deterministic quantities so as to reduce the
SDP problem to a deterministic one; see, e.g. Chapter 2.3 of \cite{be}. Such simplification arises naturally in the duality formulation.

We encapsulate the implementation procedure discussed above in Table III in Appendix \ref{ec_cov_no}. As noted in the introduction, another computational advantage of the
algorithm is that we can deploy parallel computing to expedite it. Note that, for different representative state $x^{(l)}$, $1 \le l \le L$, simulation of the associated $\xi^{(l)}|t$ and the
subsequent inner optimization in Step 1b of Table III are independent. It is easy to parallelize the execution of these procedures at different $x^{(l)}$ using multiprocessor machines.
Finally, with the help of the approximate lower bound $\underline{\widehat{\mathfrak{V}}}^{n}$ obtained from our regression-based algorithm, we can also build up a
confidence interval estimate, which many approximate dynamic programming methods are short of, for the true cost-to-go value of the original problem. One may refer to the discussion
around Table IV in \ref{ec_cov_no} and the numerical examples in the next section for details in this regard.

\subsection{Convergence Analysis}
\label{sec:convergence}

Theorem \ref{thm:convergence} establishes that in principle the DDP method should lead a convergence to the true value of the SDP problem in finite rounds of iterations. In contrast,
its regression-based implementation, as discussed in Section \ref{sec:regression}, is apparently subject to the biases coming from three sources. First, the functional
approximations built upon the basis functions may be biased relative to the true dual function $\underline{V}^{n}$. Second, the states $\{x^{(l)}, 1 \le l \le L\}$ simulated at the
beginning of algorithm execution may not be sufficiently representative. Third, the solver of the optimization problem (\ref{j}-\ref{constraint}) may only be able to find its local
optimal solution. However, in comparison with the first two errors, the error that arises in solving the deterministic optimization problem is typically not significant if a proper optimizer is used,
as suggested by the numerical examples in Section \ref{sec:num}. Hence, we focus only on the characterization of how the performance of the DDP algorithm will be affected by those factors in Theorem \ref{thm:error}.

Without loss of generality, let us assume that the state space $\mathcal{X}$ of the original problem is compact in the subsequent convergence analysis.
Many numerical examples, including the ones in Section \ref{exm2}, satisfy this assumption. In addition, for those cases with unbounded state spaces, we can obtain
approximations with sufficient accuracy by truncating the spaces into compact ones; see, for instance, \cite{alt}, \cite{kd}, \cite{dp}, and \cite{sly} for more discussions in that direction.

Consider an infinite series of basis functions $\{\psi_{m}(x), m \ge 1\}$. Suppose that we take the first $M$ functions from this set to form a functional vector $\Psi_M(x)=(\psi_1(x),\cdots,\psi_M(x))^{tr}$
to perform the DDP algorithm. We intend to characterize how its outcome will converge to the true value as we increase both the number of representative states $L$ and the number of the basis
functions $M$. We need several other technical assumptions to proceed. First,
\begin{asm}
	\label{as5}
	There exists a measure $F$, whose support is $\mathcal{X}$, such that the basis function sequence $\{\psi_m(x), m \ge 1\} $ is orthonormal under this measure $F$; that is to say,
	\begin{eqnarray*}
\int_{\mathbb{R}^n} \psi_i(x)\psi_j(x) dF(x)=\left\{
\begin{array}{rcl}
0     &      & i\neq j,\\
1     &      & i=j.
\end{array}\right.
\end{eqnarray*}
\end{asm}
Note that this assumption is not restrictive at all because we may perform the celebrated Gram-Schmidt orthogonalization to construct an orthogonal basis from any given set of linearly independent
functions.

The second assumption is about the distributions $\{G_{t}, 1 \le t \le T\}$ that are used for sampling representative states.
\begin{asm}
	\label{as1}
	Each of the sampling distributions $G_{t}(x)$ is absolutely continuous with respect to the measure $F$ in Assumption \ref{as5}. Furthermore, the Radon-Nykodym derivative between
	these two measures $dG_{t}/dF(x)$ is bounded away from zero and infinity on $\mathcal{X}$. In other words, there exist strict positive constants $\epsilon$ and $D$ such that
	$\epsilon<dG_t/dF(x)<D$ for all $x \in \mathcal{X}$.
	\end{asm}
Essentially the purpose of Assumption \ref{as1} is to help us avoid the aforementioned exploration pitfall (cf. Appendix \ref{app:ee}). The positiveness of $dG/dF$ over the entire
state space $\mathcal{X}$ ensures that the state samplers introduced in the algorithm have non-zero probability to access any part of $\mathcal{X}$.

Finally, we assume
\begin{asm}
	\label{as2}
	There exists a constant $C$ such that, for any positive integer $M$, a functional vector consisting of the first $M$ basis functions in the set,
	$\Psi_M(x)=(\psi_1(x),\cdots,\psi_M(x))^{tr}$,  satisfies
	\[
	\sup_{x\in \mathcal{X}}\left(\sum_{m=1}^{M}\psi_m^2(x)\right)^{1/2}\le CM\ \ \textrm{and}\ \ \sup_{1\le m\le M}\mathbb{E}^G[\psi_m^2(X)]\le C,
	\]
\end{asm}
and
\begin{asm}
	\label{as3}
The optimal cost-to-go function $\{V_t(x)\}_{0\le t\le T}$ is bounded on the compact set $\mathcal{X}$.
\end{asm}
Indeed, one can show that Assumption \ref{as2} holds for many popular series used in the literature on the approximation theory, including Fourier series, spline series, and
local polynomial partition series; see \cite{bcck} for a detailed discussion. The boundedness of the value function $V$ in Assumption \ref{as3} is natural if we can establish
its continuity. \cite{hl} suggest some technical conditions under which a general SDP problem has continuous value functions.

Now we turn to present the main asymptotic result for our regression-based algorithm.  Let
\begin{eqnarray*}
	\Delta_M:=\max_{0 \le t \le T}\inf_{\gamma_t=(\gamma^1_{t}, \cdots, \gamma^M_{t}) \in \mathbb{R}^M}\|V_t-\Psi_M^{tr}\gamma_t\|_{\infty},
\end{eqnarray*}
where $\|\cdot\|_{\infty}$ is the $L_\infty$ norm such that $\|f\|_{\infty}=\sup_{x\in \mathcal{X}}|f(x)|$. The quantity $\Delta_{M}$ measures the least error
magnitude that we can achieve if we approximate the true value function of the SDP problem by linearly combining the $M$ basis functions. Recall that, absent both the simulation and the
approximation errors, the dual value sequence from the DDP method should converge to the true value in at most $T+1$ iterations as shown in Theorem \ref{thm:convergence}.
Correspondingly, we develop an upper bound on the bias of $\underline{\widehat{\mathfrak{V}}}^{T+1}$, the approximate dual value after $T+1$ rounds of iterations of the regression-based
algorithm, in the next theorem.
\begin{thm}
	\label{thm:error}
	Suppose that Assumptions \ref{as5} to \ref{as3} hold. Then, there exists a constant $C$,  independent of $L$ and $M$, such that
	\begin{eqnarray}
	\label{ec:thm}
	\mathbb{E}\Big[\big|\underline{\widehat{\mathfrak{V}}}^{T+1}_0(x)-V_0(x)\big|\Big]\le
	\left(1+2l_M+C\left(\frac{M^6}{L}\right)\right)^{T}\left[(1+l_M)\Delta_{M}+C\left(\frac{M^6}{L}\right)^{1/4}\right],
	\end{eqnarray}
	where $l_{M}$ is the corresponding Lebesgue constant of the basis functions $\{\psi_{m}(x), m \ge 1\}$ (cf. Definition \ref{dfn:lebesgue} in Appendix \ref{app:con})
\end{thm}

Theorem \ref{thm:error} clearly shows how the algorithm accuracy is determined by the choice of basis functions and the amount of simulation effort. Note that both $\Delta_M$ and
$l_{M}$ are the characteristics of the basis functions that we choose. In Remarks \ref{rem1} and \ref{rem2} below, we present the corresponding orders of $\Delta_{M}$ and $l_{M}$ with
respect to $M$ under a variety of commonly used basis functions. For instance, $l_{M}$ will be bounded by a constant and $\Delta_{M}$ decays in a power order of $M$ for some choices
of basis functions. Once the set of basis functions is chosen and $M$ is fixed, we need to pick up a sufficiently large $L$ to control the right-hand side of (\ref{ec:thm}). Theorem \ref{thm:error}
spells out explicitly that $L$ should grow faster than $O(M^6)$ in order to keep such error in check. It is well known in regression analysis that a model can be overfitted if the amount of observed
data is insufficient relative to the number of regressors. When $L$, the number of simulated states on which we estimate the dual values, is not adequate in our DDP algorithm, this overfitting
effect will cause a divergence for the DDP algorithm, as illustrated by the numerical examples in Section \ref{sec:num}. From the above discussion, we can see that Theorem \ref{thm:error} can
help us understand the asymptotic behavior of our DDP estimator when both $L$ and $M$ tends to infinity for a given basis function set and thereby provide us valuable guidances on the choice of
basis functions and such parameters as $L$ and $M$. The numerical examples in Section \ref{sec:num} also show that the relative ratio between $L$ and $M$ for the DDP algorithm to converge
could be lower under some specific cases. We leave the investigation on tighter error bounds to the future work.

\begin{rem}
\label{rem1}
The approximation theory has produced some bounds on the Lebesgue constant $l_{M}$ for a variety of basis function sets. Suppose that the density function of $G$ on
$\mathcal{X}$ is bounded away from zero and infinity. We can show that $l_{M}$ should be bounded by a constant $C$ for spline series, wavelet series and local polynomial
partition series, and $l_M\le C\log(M)$ for Chebyshev polynomial series and Fourier series. See, e.g., \cite{zyg}, \cite{h2003}, \cite{bcck}, and \cite{cc}.
\end{rem}

\begin{rem}
\label{rem2}
As for $\Delta_{M}$, some studies show that, if the true value function is $s$ times continuously differentiable, the approximation error of the spline or polynomial
regressors is bounded by
 \[
\Delta_{M} \le M^{-\kappa},
 \]
where $\kappa=s/d$ and $d$ stands for the dimensionality of the function. For the proofs of this property, one may refer to Section 7.6 of \cite{dl}, Section 5.3.2 of \cite{timan},
and Theorem 12.8 of \cite{sch}.
\end{rem}

\section{Numerical Experiments}
\label{sec:num}

In this section we shall apply the regression-based Monte Carlo DDP algorithm to solve two problems related to order execution and inventory management. Both are well known to be intractable
in the literature and only approximate methods are available so far. Our algorithm demonstrates great potential in effectively assessing and improving these heuristic policies towards optimality.

\subsection{Optimal Order Execution in the Presence of Market Frictions}
\label{sec:execution}
The first numerical example we consider in the paper is an optimal order execution problem in which a trader plans to transact a large block of equity over a fixed time framework
with minimum impact costs. It can be viewed as a variant of the models proposed in \cite{bl}, \cite{ac}, and \cite{hw}. Assume that there are $n$ different assets traded in the market,
and the trader aims to acquire $\bar{\mathbf{R}}=[\bar{R}_{1}, \cdots, \bar{R}_{n}]^{tr}$ shares in each of the assets in $T$ periods. The objective of the trader is to determine a trading schedule, i.e., how many shares to purchase in each period, denoted by $\{\mathbf{S}_{1}, \cdots, \mathbf{S}_{T}\}$, $\mathbf{S}_{t} \ge 0$, $t=1, 2, \cdots, T$,
to minimize the associated transaction cost.
Let  $\mathbf{R}_{t} \in \mathbb{R}^{n}$ denote the number of shares in each asset short of the target $\mathbf{\bar{R}}$ at time $t$.
Then, a feasible trading schedule should satisfy
\begin{eqnarray}
\label{exm1:target1}
&&\sum_{t=1}^{T}\mathbf{S}_{t}=\bar{\mathbf{R}},\quad \mathbf{S}_{t} \ge 0,\quad \mathbf{S}_{t} \in \mathbb{R}^{n},\\
\label{exm1:target2}
&&\mathbf{R}_{t+1}=\mathbf{R}_{t}-\mathbf{S}_t,\ \mathbf{R}_1=\bar{\mathbf{R}}, \ \textrm{for all $t=1, \cdots, T$}.
\end{eqnarray}

To complete the statement of the problem, we must specify the price dynamics. In particular, we use $\tilde{\mathbf{P}}_{t} \in \mathbb{R}^{n}$ and $\mathbf{P}_{t} \in \mathbb{R}^{n}$ to
represent the fundamental values and actual transaction prices of all assets at time $t$, respectively, and assume that $\tilde{\mathbf{P}}_{t}$ and $\mathbf{P}_{t}$ follow the evolution laws
such that
\begin{eqnarray}
\label{exm1:dyn1}
&&\tilde{\mathbf{P}}_{t}=\tilde{\mathbf{P}}_{t-1}+\mathbf{A}\mathbf{S}_t+\mathbf{B}\mathbf{X}_{t}+\boldsymbol{\epsilon}_t,\\
\label{exm1:dyn2}
&&\mathbf{P}_t=\tilde{\textbf{P}}_{t}+h(\mathbf{S}_t),
\end{eqnarray}
for all $t$, where ${\mathbf A}\in \mathbb{R}^{n \times n}$ is a positive definite matrix and $\mathbf{B} \in \mathbb{R}^{n \times m}$. Here $\{\boldsymbol{\epsilon}_t, t=1, \cdots, T\}$ is
a sequence of white noises with mean zero and covariance matrices $\Sigma_{\epsilon}$. As shown in (\ref{exm1:dyn1}-\ref{exm1:dyn2}), our model incorporates both permanent and
temporary price impacts of transaction activities. In  (\ref{exm1:dyn1}), the constant matrix $\mathbf{A}$ is used to capture the intensity of
the permanent impact: trading the amount of $\mathbf{S}_{t}$ changes the assets' fundamental values by $\mathbf{A}\mathbf{S}_t$ and this change will last persistently in the future via
the iterative relation of $\tilde{\mathbf{P}}$. Note that this permanent price impact takes a linear form, which is a commonly adopted modeling assumption in the literature; see \cite{bl}, \cite{ac},
\cite{hs}, and \cite{hw}, for example. \cite{hs2} and \cite{g} argue that including a nonlinear permanent price impact will introduce the possibility of arbitrage.

On the other hand, we introduce the function $h(\cdot): \mathbb{R}^n \rightarrow \mathbb{R}^n$ in (\ref{exm1:dyn2}) to reflect the trading-caused impacts that will not last into the next
period. The literature documents that this kind of temporary impact in the real-life market should be concave in trading quantities (cf. \cite{bfl}). However, an assumption of concavity
often makes the control problem intractable. To demonstrate that our algorithm still works well when analytical solutions are unavailable, we assume
$h(\mathbf{S}_t)=\mathbf{D}\sqrt{\mathbf{S}_{t}}$ in the experiment, where $\mathbf{D}$ is a constant coefficient matrix.

In addition, our model allows the trader to incorporate some predictive ``signals" to extract information about the stock's future movements for
improving the performance of her trade execution. The auxiliary process $\mathbf{X}_{t} \in \mathbb{R}^{m}$ in (\ref{exm1:dyn1}) serves this purpose.
There are several possibilities proposed in the literature for the choice of
such signals. For instance, \cite{bl} suggest that
$\mathbf{X}$ could be the return of a broader market index such as S\&P500, a factor commonly used in traditional asset pricing models such as CAPM, or the outputs of an ``alpha"
model from the trader's private stock-specific analysis that is not yet impounded into market prices. In the following experiment, we abstract out the true meaning of $\mathbf{X}$ and
assume it to follow a stationary AR(1):
\begin{eqnarray}
\label{exm1:dyn3}
&&\mathbf{X}_t=\mathbf{C} \mathbf{X}_{t-1}+\boldsymbol{\eta}_t,
\end{eqnarray}
where $\mathbf{C} \in \mathbb{R}^{m \times m}$ is a matrix with all of its eigenvalues less than unity in modulus, which determines the ``decay" speed of the information, and the random
noises $\boldsymbol{\eta}_{t} \sim N(0, \Sigma_{\eta})$ are Gaussian white, independent of $\boldsymbol{\epsilon}_t$. It is worthwhile to point out that the particular form of (\ref{exm1:dyn3})
is not essential for our algorithm to work. We have tried some other specifications in the experiments for $\mathbf{X}_{t}$ and found that does not affect the effectiveness of the method. \cite{gp2013, gp2016} use the same dynamic to model the return-predicting factor in the investigation of portfolio policy when trading is costly and security returns are predictable by signals.

As aforementioned, the trader's problem is to minimize
\begin{eqnarray}
\label{exm1:objective}
\min_{\{\mathbf{S}_{t}, 1 \le t \le T\}}\mathbb{E}\left[\sum_{t=1}^{T}\mathbf{P}^{tr}_{t}\mathbf{S}_{t}\right],
\end{eqnarray}
where $\mathbf{P}^{tr}_{t}\mathbf{S}_{t}$ is how much the trader actually pays in period $t$. In \ref{app:num}, we show that this objective is indeed equivalent to
\begin{eqnarray}
\label{exm1:objective_new}
\min_{\{\mathbf{S}_{t}, 1 \le t \le T\}}\mathbb{E}\left[\sum_{t=1}^{T}\mathbf{S}_{t}^{tr}h(\mathbf{S}_{t})+\sum_{t=0}^{T-1}(\tilde{\mathbf{P}}_{t+1}-\tilde{\mathbf{P}}_{t})^{tr}\mathbf{R}_{t+1}\right].
\end{eqnarray}
The new representation (\ref{exm1:objective_new}) clearly points out two sources that are contributing to the ultimate transaction costs of the trader. The first term corresponds to the temporary
impact cost that the trader needs to pay in the process of purchasing $\bar{\textbf{R}}$ shares of assets due to the presence of  $h(\mathbf{S})$. The second term consists of the changes
in the fundamental value of the assets because of the permanent price impact that her trading activities will generate. It is easy to see that the above SDP problem has a convex objective
function. Hence, the optimal policy of the problem uniquely exists. It should be a function of state variables $\mathbf{X}$ and $\mathbf{R}$.

\medskip

\noindent Next we use the DDP algorithm to solve the minimization problem (\ref{exm1:objective_new}) with the constraints (\ref{exm1:target1}-\ref{exm1:dyn3}). The nonnegative constraint
$\mathbf{S}_t \ge 0$ turns out to be the most difficult one to deal with. As suggested by \cite{bl}, imposing it will introduce a partition structure to the optimal policy, and more seriously, the
number of partitioned regions increases combinatorially with the time horizon $T$. This renders solving the problem through the Bellman equation computationally infeasible; see also
\cite{bmdp} for more discussion on this issue in the context of a general constrained linear quadratic system. Aiming at some applications in market microstructure, \cite{ckw}
develop a partitioning algorithm for linear-quadratic Markov decision processes with linear inequality constraints. Their method recursively constructs polyhedral regions in which the optimal
value function and policy have analytical quadratic and linear forms, respectively. Note that the complexity of their method is still exponential in $T$ (cf. Notes 5 and 6 of their paper). Moreover,
it cannot be applied here because the existence of the concave temporary impact $h(\cdot)$ makes our model no longer a linear-quadratic problem.

A variety of heuristic approaches can help us derive approximate solutions to this problem. The following numerical experiments show that our DDP algorithm can be used not only for evaluating
the performance but, more importantly, to effectively improve them. Below is a summary of the heuristics that we consider in this paper.

\smallskip

\begin{itemize}
\item From a tractable simplification of the problem.  It is straightforward to see that, if we ignore the no-sales constraint $\mathbf{S}_{t} \ge 0$ and the temporary price impact $h(\mathbf{S}_{t})$, the problem (\ref{exm1:objective_new}) with the constraints (\ref{exm1:dyn1}-\ref{exm1:dyn3}) indeed degenerates to a standard LQC. The computation in
\ref{app:num} shows that the optimal policy of this simplified problem is analytically known:
\begin{eqnarray}
\label{exm1:policy1}
\tilde{\mathbf{S}}_t(\mathbf{X}_t, \mathbf{R}_t)=\left(\mathbf{I}-\frac{1}{2}\mathbf{Q}_{t+1}^{-1}\mathbf{A}^{tr}\right)\mathbf{R}_t+\left(\frac{1}{2}\mathbf{Q}_{t+1}^{-1}\mathbf{K}_{t+1}\mathbf{C}\right)\mathbf{X}_t,
\end{eqnarray}
where $\mathbf{Q}_{t} \in \mathbb{R}^{n \times n}$ and $\mathbf{K}_{t} \in \mathbb{R}^{n \times m}$ are two matrices that can be determined by the matrix equations in
(\ref{app:exm1}-\ref{app:exm2}). This linear policy emphasizes the importance of trading on signals in the process of meeting the execution target, as the current information level
$\mathbf{X}_t$ affects the amount of trading volume $\tilde{\mathbf{S}}_t$. However, such $\tilde{\mathbf{S}}_t$ is not feasible to the original problem because it may take
negative values when $\mathbf{X}_t$ is negative. To restore a feasible policy, we may project  $\tilde{\mathbf{S}}_t$ into the region $[0, \bar{\mathbf{R}}]$ by letting
\begin{eqnarray}
\label{exm1:policy}
\mathbf{S}^{LQ}_t(\mathbf{X}_t, \mathbf{R}_t)=\min\left(\max\left(\tilde{\mathbf{S}}_t(\mathbf{X}_t, \mathbf{R}_t),0\right),\mathbf{R}_t\right).
\end{eqnarray}
\item From linear program approximation. Introduced by \cite{ss} and further developed by \cite{dv, dv2} and \cite{dfm2, dfm}, the linear programming based approach provides us an attractive
way to construct approximate solutions to the dynamic programs. Consider a collection of basis functions $\{\psi_{1}, \cdots, \psi_{K}\}$ and use the following regression to approximate the
optimal value functions at every time $t=1, \cdots, T$:
\begin{eqnarray*}
V_{t}(\mathbf{X}_t, \mathbf{R}_t) \approx \sum_{k=1}^{K}\theta_{k, t}\psi_{k}(\mathbf{X}_t, \mathbf{R}_t),
\end{eqnarray*}
where $\theta_{t}=(\theta_{1, t}, \cdots, \theta_{K, t})$ is the regression coefficient to be determined. As noted in the discussion around Definition \ref{dfn:subsolution}, the true value function
of an SDP must be the largest subsolution. Select some representative states $\{(\mathbf{X}^i, \mathbf{R}^i): i=1, \cdots, I\}$. Let $c_{t, i}$ be a positive constant for all $t=1, \cdots, T$ and
$i=1, \cdots, I$. We may recast this fact as a linear program for the problem (\ref{exm1:objective_new}):
\begin{eqnarray*}
\max_{\theta_{t}: t=1, \cdots, T} \sum_{i=1}^{I} c_{t, i} \sum_{k=1}^K \theta_{k,t} \psi_k(\mathbf{X}^i, \mathbf{R}^i)
\end{eqnarray*}
subject to
\begin{eqnarray}
  \label{lp_approx}
  \resizebox{.9\hsize}{!}{$
  \sum_{k=1}^K \theta_{k,t} \psi_k(\mathbf{X}^i, \mathbf{R}^i) \le \min_{\mathbf{S}_{t} \ge 0} \mathbb{E}\left[\mathbf{S}_{t}^{tr}h(\mathbf{S}_{t}) + (\tilde{\mathbf{P}}_{t+1}
  -\tilde{\mathbf{P}}_{t})^{tr}\mathbf{R}_{t+1} + \sum_{k=1}^K \theta_{k,t+1} \psi_k(\mathbf{X}_{t+1}, \mathbf{R}_{t+1})\Big|(\mathbf{X}_t, \mathbf{R}_t)=(\mathbf{X}^i, \mathbf{R}^i)\right].
  $}
\end{eqnarray}
Note that the constraint (\ref{lp_approx}) is just a rephrasing of Definition \ref{dfn:subsolution} and it is a linear inequality with respect to the regression coefficient $\theta_t$.
\item Lookahead. One-step and multistep lookahead constitute another class of commonly used approaches to produce approximate solutions to dynamic programs.
We replace $V_{t+1}$ in the Bellman equation (\ref{bellman2}) with any of its approximation $\tilde{V}_{t+1}$. Then, the minimization
\begin{eqnarray}
\label{lookahead}
\resizebox{0.9\hsize}{!}{$
\mathbf{S}^{LO}_{t}(\mathbf{X}_t, \mathbf{R}_t)=\argmin_{\mathbf{S}_{t} \ge 0}\mathbb{E}\left[\mathbf{S}_{t}^{tr}h(\mathbf{S}_{t}) + (\tilde{\mathbf{P}}_{t+1}
-\tilde{\mathbf{P}}_{t})^{tr}\mathbf{R}_{t+1} +\tilde{V}_{t+1}(\mathbf{X}_{t+1}, \mathbf{R}_{t+1})\Big| (\mathbf{X}_t, \mathbf{R}_t)\right] $}
\end{eqnarray}
defines the one-step lookahead policy at state $(\mathbf{X}_t, \mathbf{R}_t)$. For the purpose of illustration, we make use of the value function of the simplified problem discussed in the
first bullet as $\tilde{V}$. To derive the multistep lookahead, we can minimize the cost of the first $L>1$ steps with the future cost approximated by a function $\tilde{V}_{t+L}$.
\item Backward dynamic programming. To overcome the curse of dimensionality encountered in utilizing the equations (\ref{bellman1}-\ref{bellman2}), one may use the basis functions to
obtain low-dimensional regression representations of the value functions and repeatedly substitute them into the one-step Bellman's equation (\ref{bellman2}) to
produce approximate solutions to the problem in a backward fashsion. The regression coefficients can be estimated by using the least square method on some representative states that are
fixed beforehand.
\end{itemize}

\smallskip

Starting from any of these heuristics, our DDP algorithm demonstrates a strong ability to construct improved approximations for all of them. Table \ref{exm1:t1} displays the related
convergence results. Here we consider  a case with three assets and a signal vector of two variables, i.e. $\mathbf{R}=[R_{1}, R_{2}, R_{3}]^{tr}$ and $\mathbf{X}=[X_1, X_2]^{tr}$.
To deal with this 5-dim problem, we use the following set of basis functions in the experiment:
\begin{eqnarray}
\label{exm1:basis1}
&&\Big\{1, (X_{i})_{i=1, 2},\ (R_k)_{k=1,2,3},\ (X_{i}X_{j})_{1\le i,j \le 2},\ (R_{k}R_{l})_{1\le k,l\le 3},\ (R_{k}X_{j})_{1\le k\le 3,\ 1\le j\le 2},\nonumber\\
&&\ (R_{k}\sqrt{R_{k}})_{1\le k\le 3},\ (X_{i}^3)_{1\le i\le m},\ (R_{i}^3)_{1\le i\le n},\ (X_{i}^4)_{1\le i\le m},\ (R_{i}^4)_{1\le i\le n}\Big\}.
\end{eqnarray}
The abbreviation $(X_{i})_{i=1, 2}$, for example, represents that both functions $X_{1}$ and $X_{2}$ are included. The other notations should be understood in
the same way. We also include a constant, represented by 1 in the set (\ref{exm1:basis1}), in the regressors. In the interest of space, the values of all the model
parameters are reported in Appendix \ref{app:num}.

As noted in Section \ref{sec:mc}, we need a state selector $G$ to generate a number of representative pairs of $(\mathbf{X}, \mathbf{R})$ in the state space at each period $t$
so that we can run regressions to extrapolate the dual values observed on these pairs. Note that the signal process $\mathbf{X}_{t}$ has an autonomous dynamic (\ref{exm1:dyn3}),
independent of the control policies taken by the trader. We thereby use its marginal distribution in the experiment to simulate samples for $\mathbf{X}$. Meanwhile, since the sample
trajectory of $\mathbf{R}_{t}$ resides in $[0, 10^5]^3$ under any trading scheme, we take the uniform distribution in this cube to sample $\mathbf{R}$. To speed up the overall calculation
when evaluating the dual values, we parallelize the simulation of the state pairs $(\mathbf{X}, \mathbf{R})$ and random noises $(\boldsymbol{\epsilon}_{t}, \boldsymbol{\eta}_{t}, t=1, \cdots, 20)$
to multicore CPUs (32 cores in our experiments) and solve the corresponding optimization programs simultaneously.

\begin{table}[htbp]
	\centering
	\renewcommand\arraystretch{1}
	\begin{tabular}{cccccccc}
		\toprule
		Approximation & Iteration & Dual Values & (SE)  && Primal Values & (SE)   & Gap  \\
		\midrule
		& 0     &       &              && 325.49  & (0.49)  & 18.88\% \\
		\cdashline{6-8}
	        & 1     & 264.05 & (2.71)   &\diagbox[dir=SW,height=2em]{}{}& -  & -  & - \\
		\cdashline{3-4}
		LQ & 2     & 267.12 & (0.42)   && -  & -  & - \\
		& 3     & 269.26 & (0.25)   && 272.82  & (0.53)  & 1.28\% \\
		\cdashline{6-8}
		& 4     & 269.33 & (0.25)   &\diagbox[dir=SW,height=2em]{}{}&       &       &  \\
		\cdashline{3-4}
		\bottomrule
		\toprule
		Approximation & Iteration & Dual Values & (SE)  && Primal Values & (SE)   & Gap  \\
		\midrule
		& 0     &       &              && 433.28  & (8.21)  & 40.38\% \\
		\cdashline{6-8}
		& 1     & 258.32 & (0.63)   &\diagbox[dir=SW,height=2em]{}{}& -  & -  & - \\
		\cdashline{3-4}
		& 2     & 266.38 & (0.62)   && -  & -  & - \\
		Linear & 3     & 268.76 & (0.32)   && -  & -  & - \\
		& 4     & 269.63 & (0.26)   && 272.44  & (0.51)  & 1.06\% \\
		\cdashline{6-8}
		& 5     & 269.55 & (0.25)   &\diagbox[dir=SW,height=2em]{}{}&       &       &  \\
		\cdashline{3-4}
		\bottomrule
		\toprule
		Approximation & Iteration & Dual Values & (SE)  && Primal Values & (SE)   & Gap  \\
		\midrule
		& 0     &       &              && 324.80  & (0.49)  & 18.64\% \\
		\cdashline{6-8}
		& 1     & 264.23 & (2.59)   &\diagbox[dir=SW,height=2em]{}{}& -  & -  & - \\
		\cdashline{3-4}
		Lookahead & 2     & 266.99 & (0.44)   && -  & -  & - \\
		& 3     & 269.10 & (0.25)   && 272.38  & (0.50)  & 1.12\% \\
		\cdashline{6-8}
		& 4     & 269.33 & (0.25)   &\diagbox[dir=SW,height=2em]{}{}&       &       &  \\
		\cdashline{3-4}
		\bottomrule
		\toprule
		Approximation & Iteration & Dual Values & (SE)  && Primal Values & (SE)   & Gap  \\
		\midrule
		& 0     &       &              && 285.55  & (0.50)  & 8.21\% \\
		\cdashline{6-8}
		Backward & 1     & 262.08 & (2.80)   &\diagbox[dir=SW,height=2em]{}{}& -  & -  & - \\
		\cdashline{3-4}
		& 2     & 269.53 & (0.24)   && 272.70  & (0.50)  &1.28\% \\
		\cdashline{6-8}
		& 3     &  269.56& (0.23)   &\diagbox[dir=SW,height=2em]{}{}&       &       &  \\
		\cdashline{3-4}
		\bottomrule				\end{tabular}%
    \captionsetup{font=scriptsize}
	\caption{
  The convergence results of the DDP algorithm in the example of order execution. The four approximation methods are used to construct the initial policies as the inputs to
	the DDP algorithm. We denote them by LQ, Linear, Lookahead, and Backward, respectively, in the table. We simulate $K=1\times10^4$ sample paths of random noises
	$(\boldsymbol{\epsilon}_{t}, \boldsymbol{\eta}_{t}, t=1, \cdots, 20)$ to estimate their corresponding values, which are reported in the cell of Primal Values in Iteration 0 of every subparts
	of the table.  The standard error of this policy estimation is shown in the column ``(SE)". The entry in Row ``Iteration 1" and Column ``Dual Values" displays the dual value associated with
	each approximate policy. We sample $L=1.5 \times 10^4$ pairs of $(\mathbf{X}, \mathbf{R})$ in each time step from the distribution $G$ mentioned in the body text to compute the dual
	values in each iteration. The same distribution $G$ is also used to generate representative states for the methods of linear programming approximation and backward dynamic programming.
	In LP, we simulate 300 state pairs and thereby solve a linear program with 300 constraints. In Backward DP, we simulate $1.5 \times 10^4$ states for carrying out the least square estimation.
	The numbers in the parentheses in the column next to ``Dual Values" are the standard errors of the dual estimations. The percentage gaps in the last column of the table are computed according
	to the ratio of $(\textrm{Primal}-\textrm{Dual})/\textrm{Primal}$. The default parameters used in the experiments are $\lambda=10$ and $\delta=1$. The values of other parameters are reported in
	\ref{app:num}. All the computation experiments are conducted on a PC equipped with an Intel Xeon 32-core 2.93 GHz CPU and 12.0 GB of RAM. The computation environment is Windows 7
	and MATLAB R2017a and parallel pool. The average computational time is 1416.1s per iteration.}
	\label{exm1:t1}%
\end{table}%

Table \ref{exm1:t1} consists of four subparts. Each of them reports the respective convergence results for the four approximate heuristics.  We first assess the performance of each approximate policy
by evaluating its corresponding average transaction costs along $K=1 \times 10^4$ simulated paths of random noises $(\boldsymbol{\epsilon}_{t}, \boldsymbol{\eta}_{t}, t=1, \cdots, 20)$.  The
outcome is reported in the first row of each subpart. Meanwhile, we compute the dual value associated with each heuristic policy in the second row of the column ``Dual Values". All the approximations
have significant duality gaps, which show that the performance of all the policies are not satisfactory.

Consistent with the theoretical convergence results in the last sections, the dual values increase as we run more iterations of the DDP algorithm, no matter which approximate heuristic we start
with. These dual values thereby provide a sequence of increasingly tighter lower bounds to the true value of the problem. The algorithm terminates after several rounds of iterations when it produces no
essential changes on the dual value. More precisely, the termination criterion is that the dual value in the penultimate iteration falls within the 95\% confidence interval of the dual value in the terminal
iteration. We then apply the direct policy evaluation scheme (cf. Table IV in \ref{ec_cov_no}) to estimate the value of the policy obtained through our DDP algorithm. As shown by the last row of each
subpart, the dual gap of the improved policy shrinks down to around 1\%, strongly suggesting that the new policy is very close to the optimality. In addition, we find that, irrespective of the initial
approximation that we start with, all the final outcomes that the DDP algorithm converges to are identical. Denote hereafter the policy we obtain by $\mathbf{S}^{DDP}$.

In this experiment, we use the DC programming to solve the inner optimization problem in the dual formulation. The consideration underlying this choice is that the penalty $z$, one part of the
objective function of the inner optimization problem (cf. (\ref{inner}) and (\ref{j})), takes a very special form of functional difference. After decomposing the objective function to the difference of
two convex functions in $\mathbf{S}_t$, we rely on sequential convex relaxation to transform the optimization job down to solving a sequence of convex programs. The literature
has established the property of global convergence for this approach; that is, starting from any given initial point, the sequence generated by it converges to a solution to DC programs
that satisfies the Karush-Kuhn-Tucker condition; see \cite{yr}, \cite{lp2005}, \cite{ls2009}, \cite{l}, \cite{lp2018}, and \cite{bv}. A brief introduction on the DC programming is also provided in
Appendix \ref{app:dc}.

While in theory it is possible that the above sequential convex programming may only lead to local optimal solutions for the inner optimization problem, we need to stress that the numerical
evidence shows that does not affect the optimality of $\mathbf{S}^{DDP}$ reported in Table \ref{exm1:t1}. To see this, we develop a sanity check in the following remark.
\begin{rem}
At the termination of the DDP algorithm, we expand the output dual value function $\underline{\mathfrak{V}}_{t}$ to its first order, i.e., for $t=0. 1, \cdots, T-1$,
\begin{eqnarray}
\label{sanity}
\underline{\mathfrak{V}}_{t}(\mathbf{X}, \mathbf{R}) \approx \underline{\mathfrak{V}}_{t}(\mathbf{X}^0, \mathbf{R}^0) + \nabla_{\mathbf{x}}\underline{\mathfrak{V}}_{t}(\mathbf{X}^0, \mathbf{R}^0)(\mathbf{X}-\mathbf{X}^0)+ \nabla_{\mathbf{R}}\underline{\mathfrak{V}}_{t}(\mathbf{X}^0, \mathbf{R}^0)(\mathbf{R}-\mathbf{R}^0)
\end{eqnarray}
where $(\mathbf{X}^0, \mathbf{R}^0)$ is a state pair that we fix in advance, and $\nabla_{\mathbf{x}}$ and $\nabla_{\mathbf{R}}$ are the gradients with respect to $\mathbf{X}$ and
$\mathbf{R}$, respectively. Note that $\underline{\mathfrak{V}}_{t}$ is indeed a linear combination of the basis functions in (\ref{exm1:basis1}). Hence, the function on the right hand side of
(\ref{sanity}) is explicitly known and it is linear in the variable $\mathbf{R}$. In Table \ref{exm1:t1} (cf. the dual value in the last row of each subpart), we use  $\underline{\mathfrak{V}}_{t}$ to
construct a duality to assess the quality of $\mathbf{S}^{DDP}$. Alternatively we may substitute the linear function on the right-hand side of (\ref{sanity}) into (\ref{penalty_def}) to construct
another penalty. Note that the inner optimization problem in the resulting dual formulation will become a convex program, which is globally solvable. So we do not need to worry about the issue
of local solutions for this new duality.  It turns out that the dual value we obtain in this way is 269.47 with a standard deviation 0.25, which is very close to the ultimate dual values reported in
Table \ref{exm1:t1}. This strongly suggests that the sequential convex programming procedure can effectively lead us to find a policy with adequate performance.
\end{rem}

\begin{figure}[htbp]
		\centering
		\includegraphics[width=2.7in]{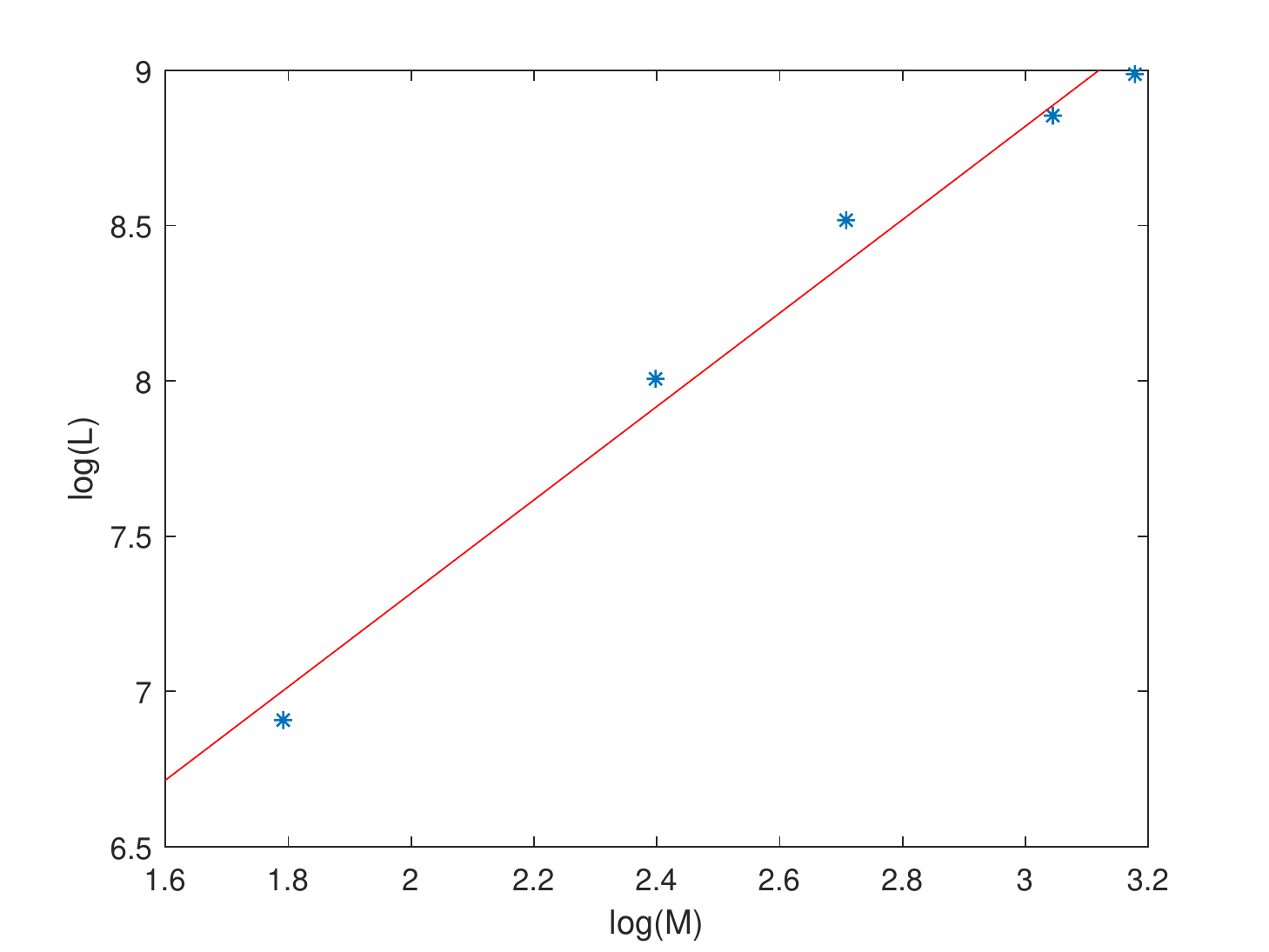}
		\caption{The regression result of the critical $\log(L)$ against $\log(M)$. As shown in Table \ref{ec:ml}, $L$ has to be as large as
		$[1000,3000,5000,7000,8000]$, respectively, for $M=[6, 11, 15, 21, 24]$ to ensure the convergence of the DDP algorithm.
		The red straight line is the linear extrapolation between $\log(L)$ and $\log(M)$. We have approximately $\log(L)=1.5\log(M)+4.3$. }
		\label{ec:ml_fig}
\end{figure}

Recall that Theorem \ref{thm:error} reveals a crucial trade-off between the model complexity and the sampling adequacy facing us in the implementation of the DDP algorithm; that is, given the
number of basis functions $M$, we need a sufficiently large number of samples $L$ to ensure the convergence of the DDP algorithm. Both Table \ref{ec:ml} and Figure \ref{ec:ml_fig} corroborate
this conclusion. In Table \ref{ec:ml}, we can easily see that, for a fixed basis function set, there exists a minimum $L$ for the DDP algorithm to converge. Moreover, as the number of basis functions
used in the approximation increases, this critical $L$ tends to become larger. Figure \ref{ec:ml_fig} empirically examines how fast this minimum $L$ grows with $M$ using the log-log plot. The slope
suggests that the number of representative states $L$ should be at least as large as $O(M^{3/2})$ to ensure the convergence of the DDP algorithm. Note that this rate is much smaller
than the theoretical rate established in Theorem \ref{thm:error}. We leave the research on tightening the bound to future work.
\begin{table}[htbp]
  \centering
  \renewcommand\arraystretch{1.5}
    \begin{tabular}{cccccccccccccccc}
    \toprule
     & \multicolumn{15}{c}{$L$}\\
     \cmidrule{2-16}
        $M$ & 500   & & 1000  & 2000  &&3000  & 4000  && 5000  & 6000  && 7000  && 8000  & 9000 \\
        \midrule
          	\cdashline{2-2}
    6     & $\times$ &\multirow{1}{0.2cm}{ \diagbox[height=1.8em]{}{}} & $\surd$& $\surd$& &$\surd$&         $\surd$  && $\surd$      &$\surd$       &&     $\surd$  && $\surd$      & $\surd$ \\
          	\cdashline{4-5}
    11    &    $\times$   && $\times$& $\times$   &\multirow{1}{0.2cm}{ \diagbox[height=1.8em]{}{}} &$\surd$ &$\surd$&& $\surd$ &  $\surd$     &&  $\surd$     && $\surd$      & $\surd$ \\
          \cdashline{7-8}
    15    &   $\times$   &&    $\times$   &    $\times$   &&$\times$& $\times$&\multirow{1}{0.2cm}{ \diagbox[height=1.8em]{}{}} &$\surd$& $\surd$ && $\surd$ && $\surd$      & $\surd$ \\
               \cdashline{10-11}
    21    &   $\times$  &&    $\times$   &    $\times$   &&     $\times$ &     $\times$  && $\times$&$\times$&\multirow{1}{0.2cm}{ \diagbox[height=1.8em]{}{}} &$\surd$& & $\surd$ &  $\surd$ \\
         \cdashline{13-13}
    24    &   $\times$  &&    $\times$   &    $\times$   &     &$\times$ &     $\times$  && $\times$&$\times$&& $\times$&\multirow{1}{0.2cm}{ \diagbox[height=1.8em]{}{}} &  $\surd$ &  $\surd$ \\
              \cdashline{15-16}
    \bottomrule
    \end{tabular}%
    % \captionsetup{font=scriptsize}
    \captionsetup{font=small,justification=centering}
  \caption{The convergence performance of the DDP algorithm under different choices of $L$ and $M$. If it converges, we input $\surd$ in the corresponding entry; otherwise, we use $\times$. The basis function sets we choose for each row are $\{1, (X_{i})_{i=1, 2},\ (R_k)_{k=1,2,3} \}$, $\{1, (X_{i})_{i=1, 2},\ (R_k)_{k=1,2,3},\ (X_{i}^2)_{1\le i \le 2},\ (R_{k}^2)_{1\le k\le 3}\}$, $\{1, (X_{i})_{i=1, 2},\ (R_k)_{k=1,2,3},\ (X_{i}X_{j})_{1\le i,j \le 2},\ (R_{k}R_{l})_{1\le k,l\le 3}\}$, $\{1, (X_{i})_{i=1, 2},\ (R_k)_{k=1,2,3},\ (X_{i}X_{j})_{1\le i,j \le 2},\ (R_{k}R_{l})_{1\le k,l\le 3},\ (R_{k}X_{j})_{1\le k\le 3,\ 1\le j\le 2}\}$ and $\{1, (X_{i})_{i=1, 2},\ (R_k)_{k=1,2,3},\ (X_{i}X_{j})_{1\le i,j \le 2},\ (R_{k}R_{l})_{1\le k,l\le 3},\ (R_{k}X_{j})_{1\le k\le 3,\ 1\le j\le 2},(R_{k}\sqrt{R_k})_{1\le k\le 3}\}$, respectively. }
  \label{ec:ml}%
\end{table}%

\medskip

\begin{figure}[htbp]
	\centering
	\subfigure{
		\includegraphics[width=5.1cm]{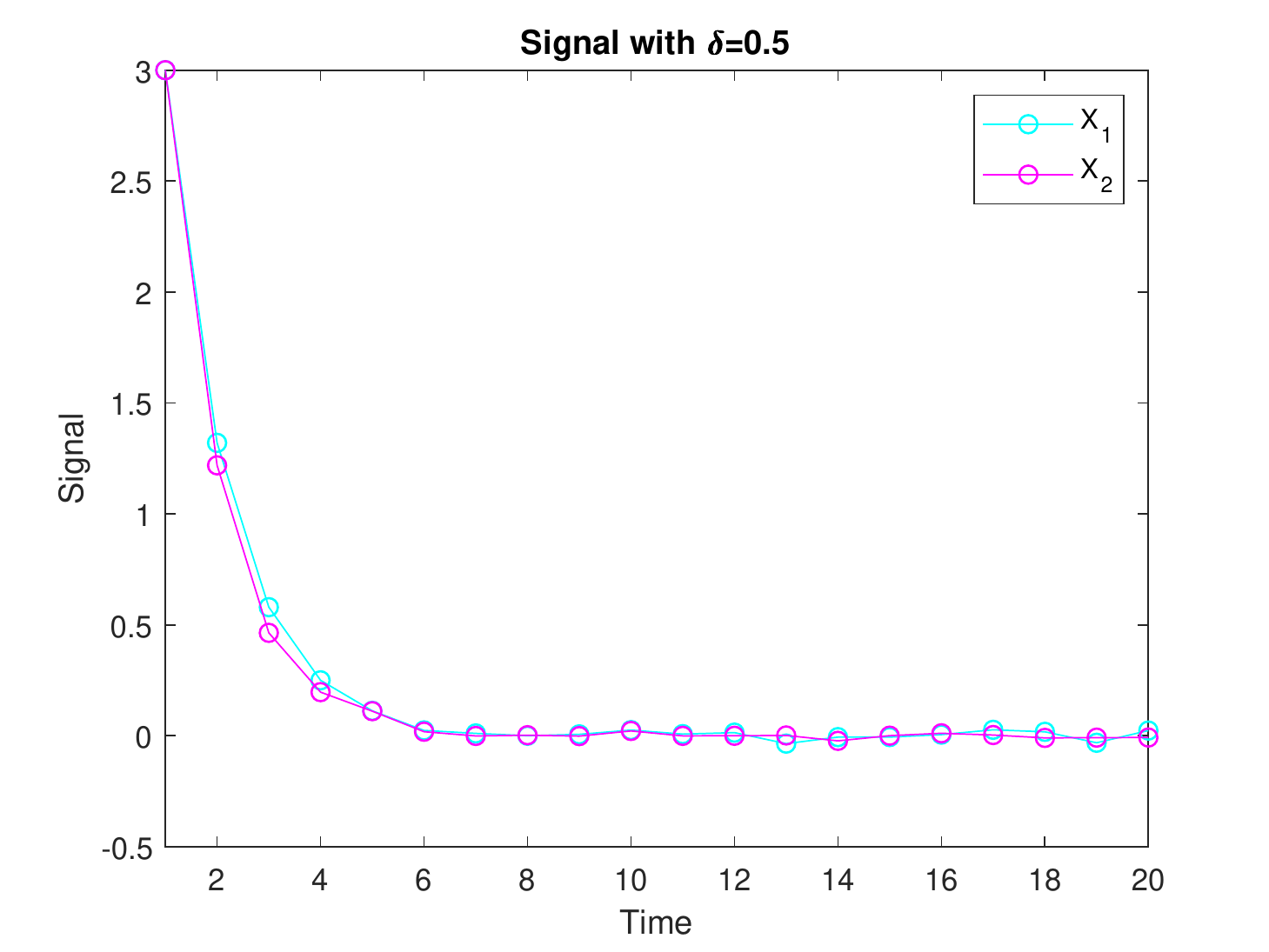}
	}
	\subfigure{
		\includegraphics[width=5.1cm]{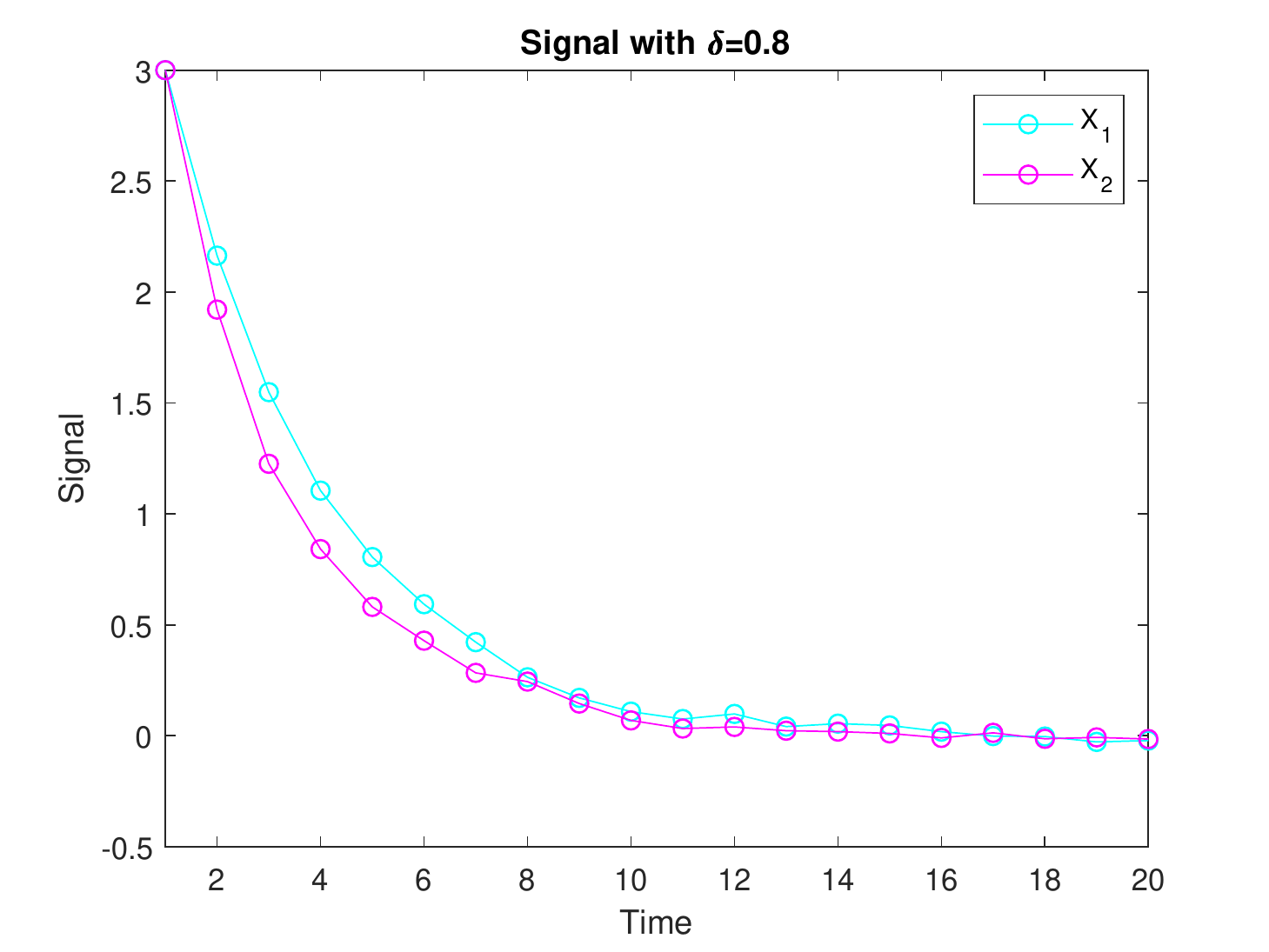}
	}
	\subfigure{
		\includegraphics[width=5.1cm]{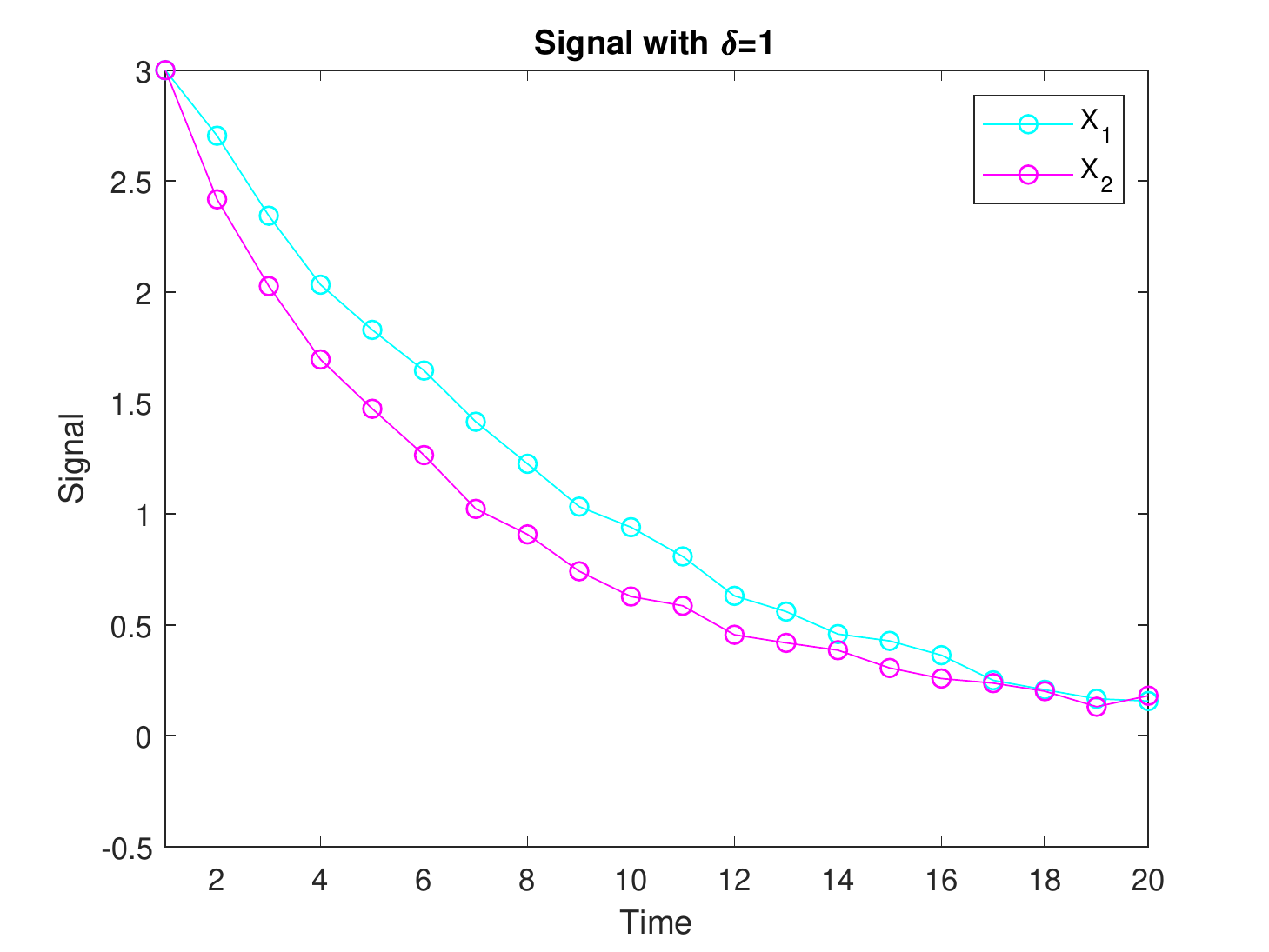}
	}
	\subfigure{
		\includegraphics[width=5.1cm]{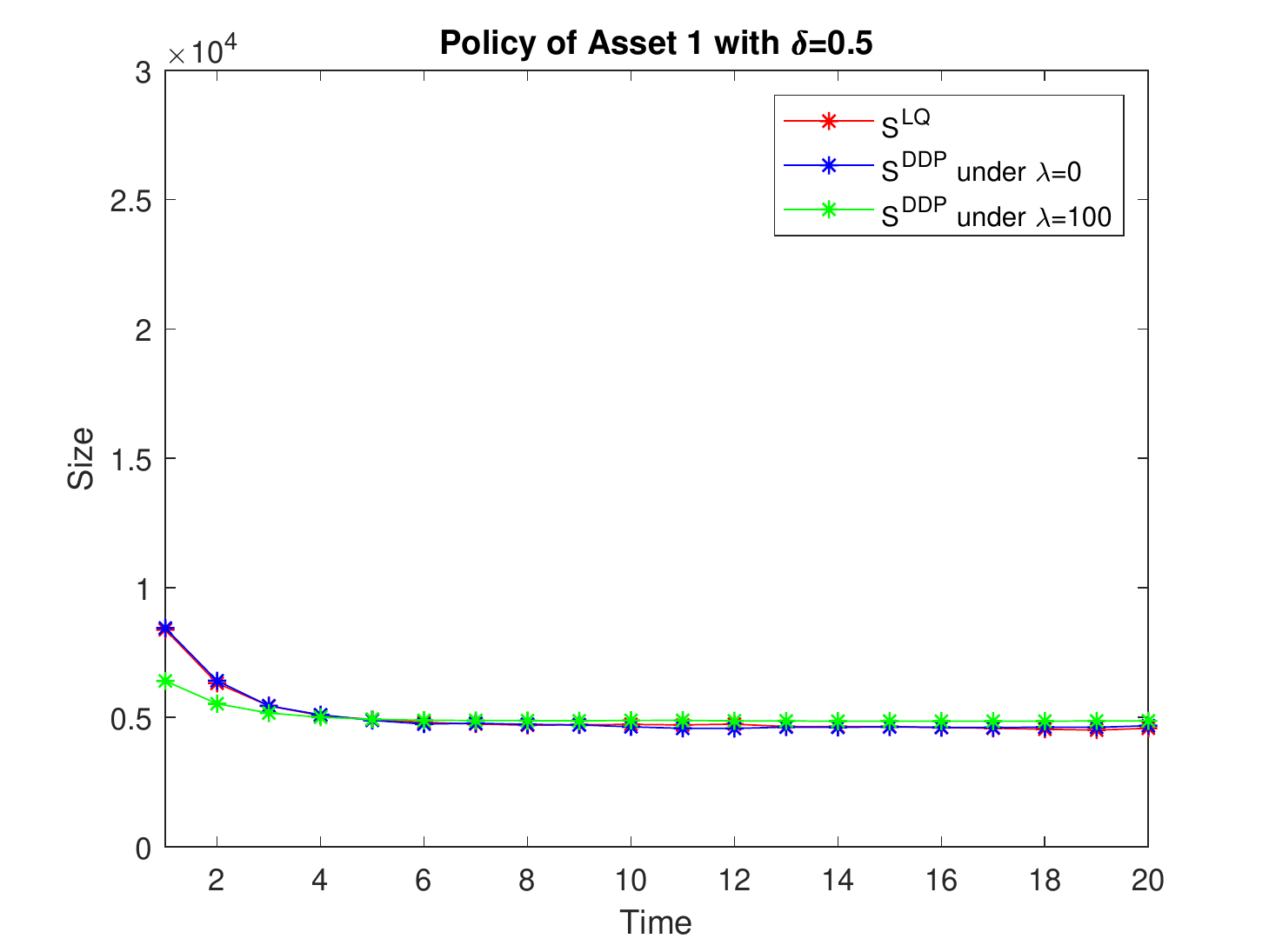}
	}
	\subfigure{
		\includegraphics[width=5.1cm]{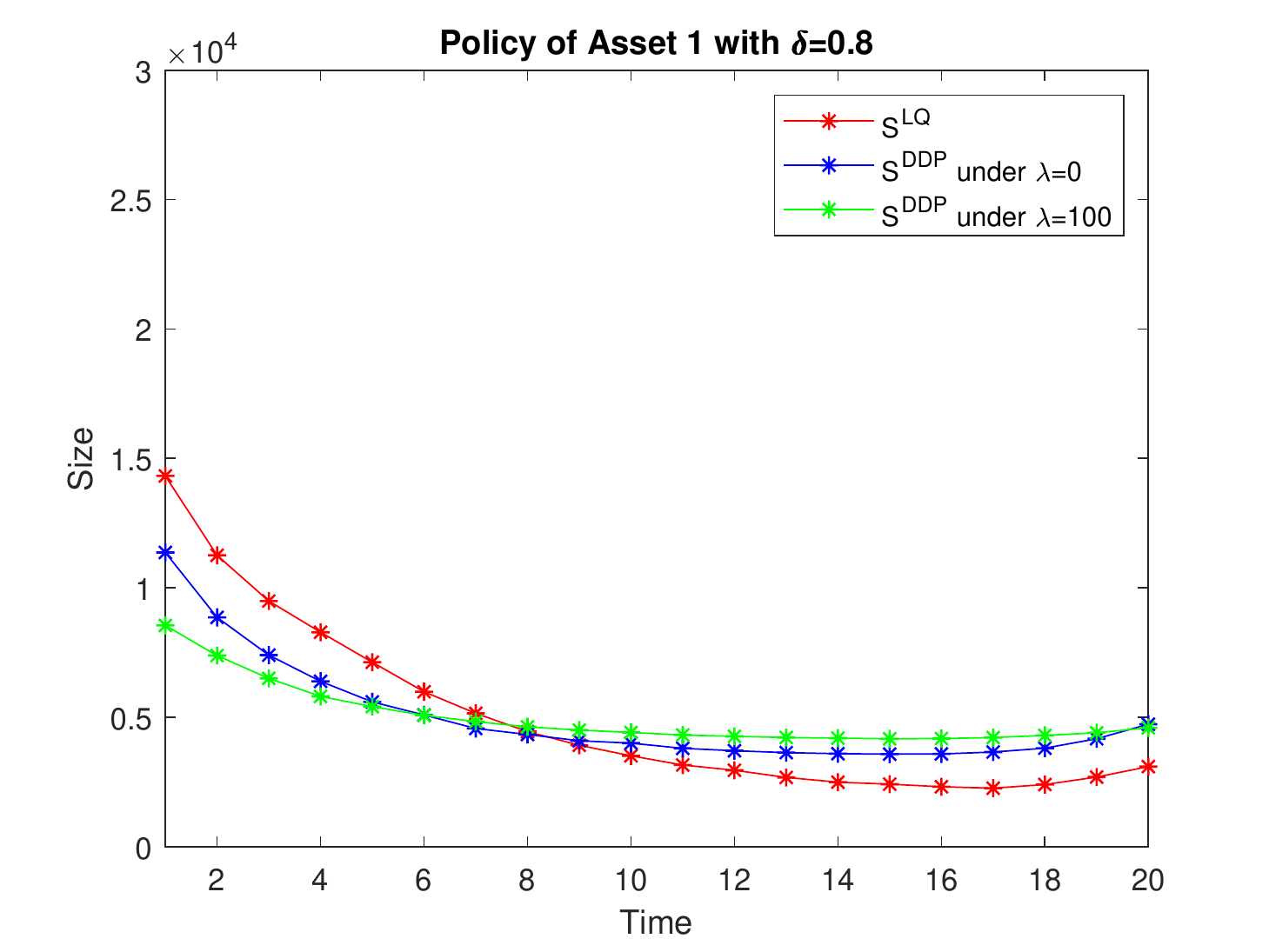}
	}
	\subfigure{
		\includegraphics[width=5.1cm]{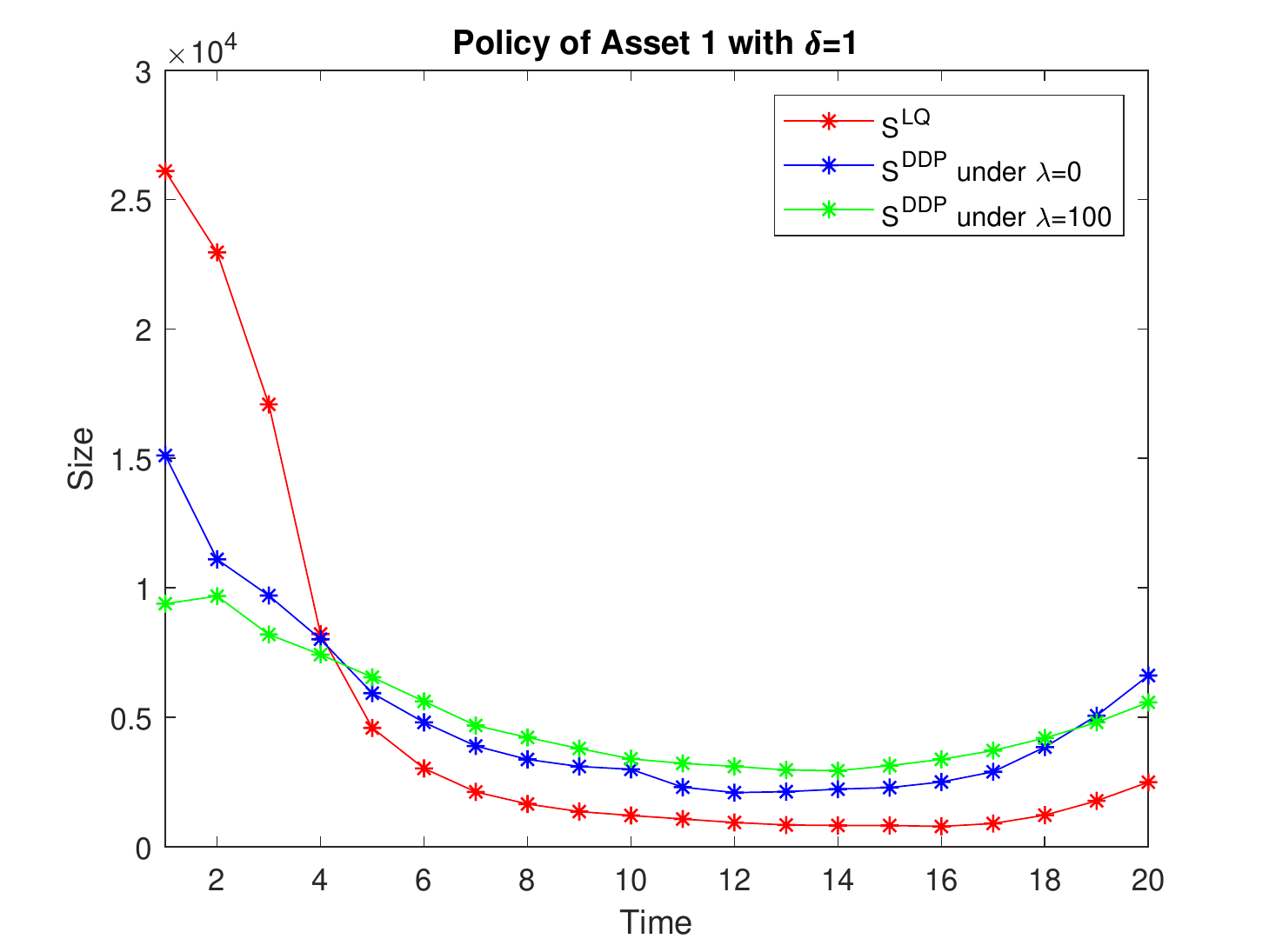}
	}
	\subfigure{
		\includegraphics[width=5.1cm]{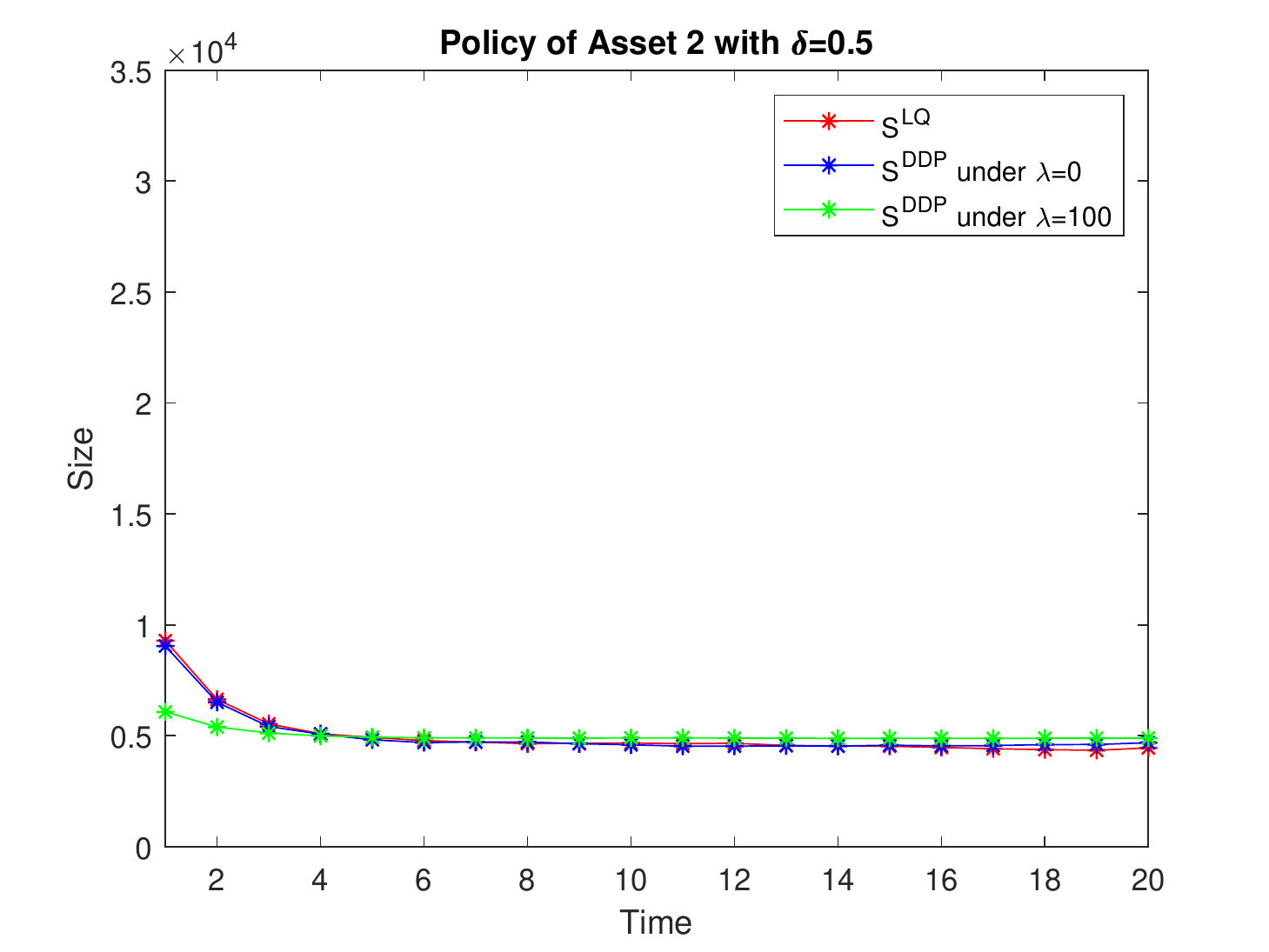}
	}
	\subfigure{
		\includegraphics[width=5.1cm]{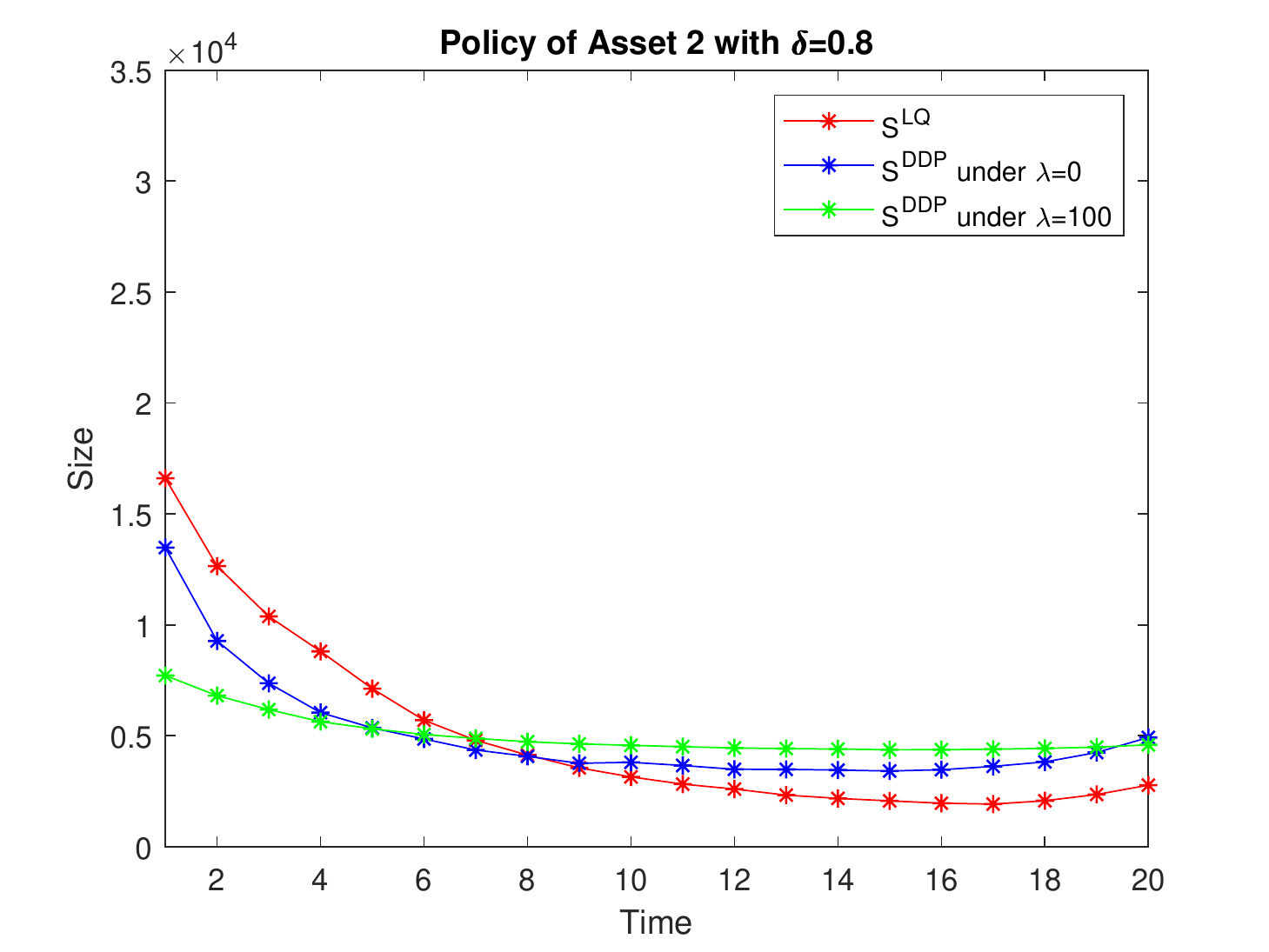}
	}
	\subfigure{
		\includegraphics[width=5.1cm]{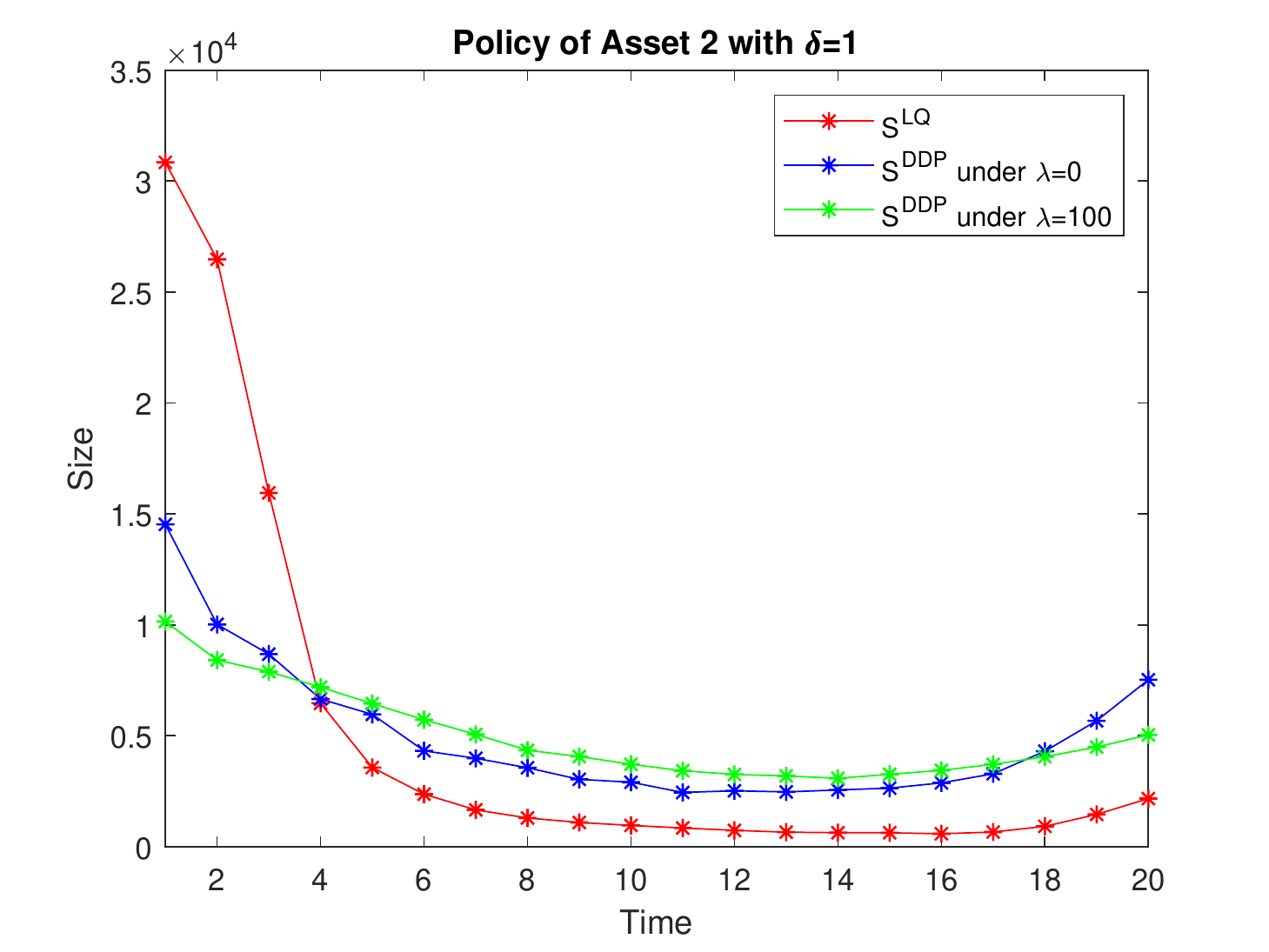}
	}
	\subfigure{
		\includegraphics[width=5.1cm]{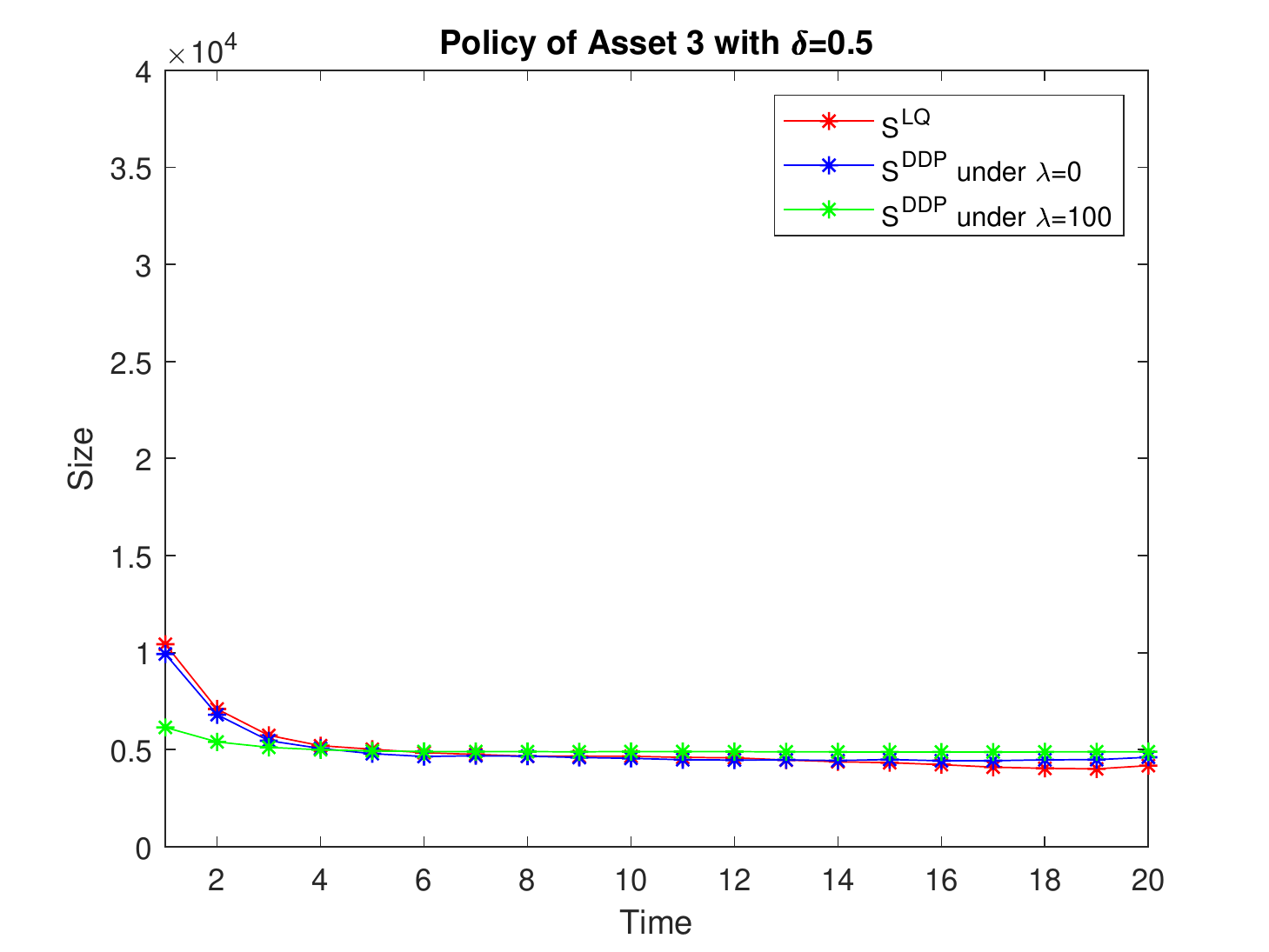}
	}
	\subfigure{
		\includegraphics[width=5.1cm]{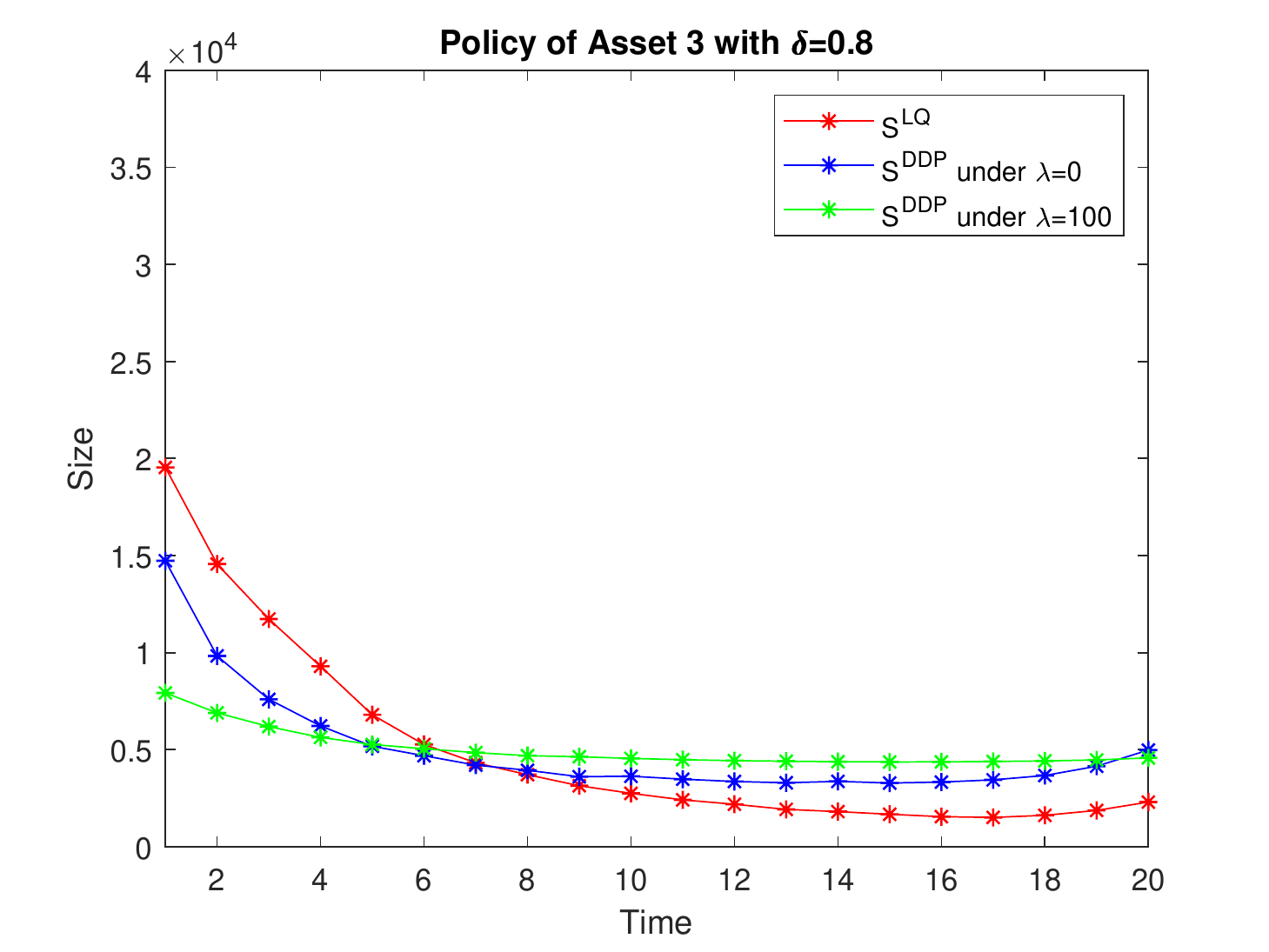}
	}
	\subfigure{
		\includegraphics[width=5.1cm]{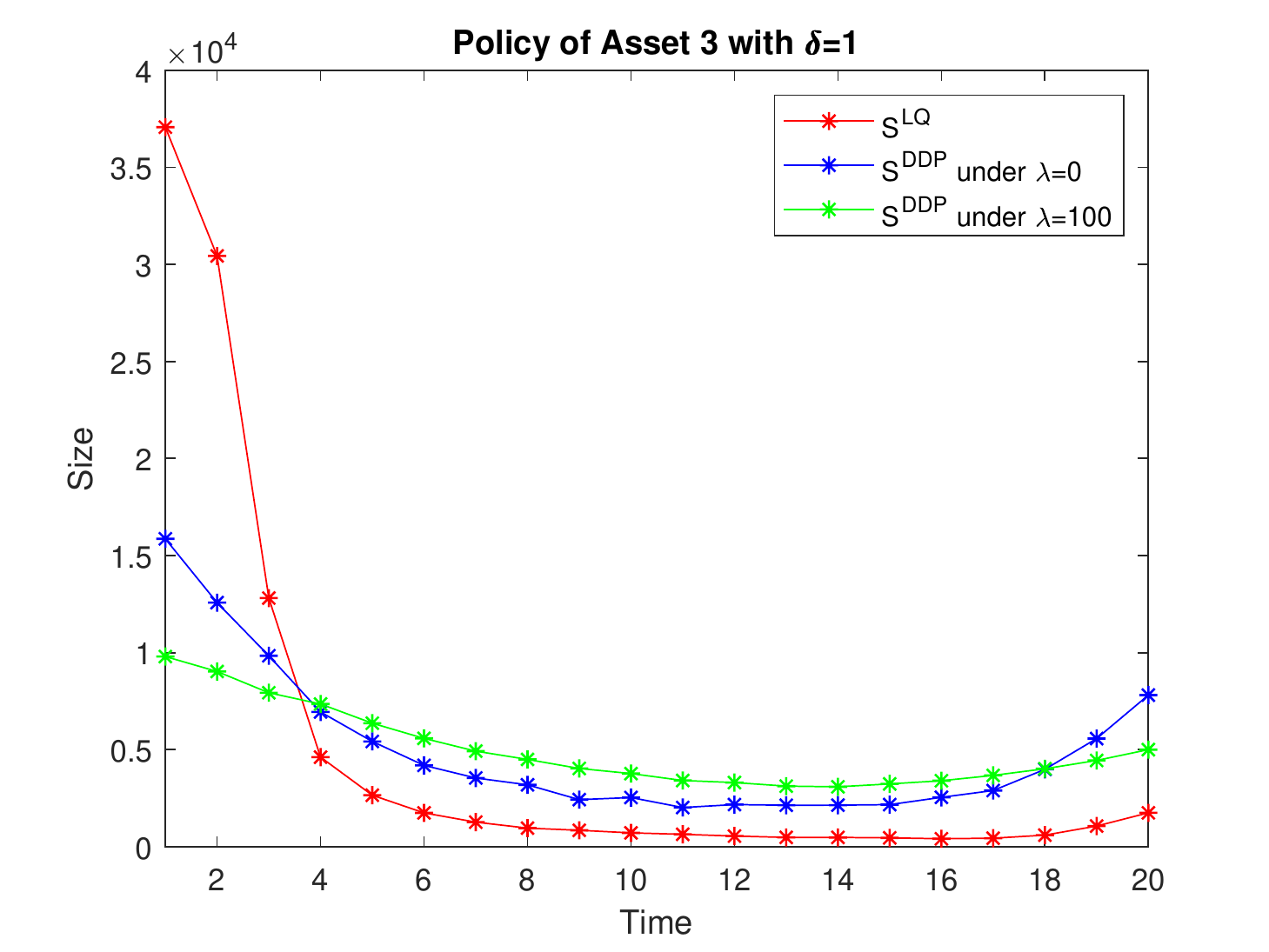}
	}
	\caption{A simulation comparison between two policies $\mathbf{S}^{DDP}$ and $\mathbf{S}^{LQ}$. We simulate 10,000 sample paths for random noises
	$(\boldsymbol{\epsilon}_{t}, \boldsymbol{\eta}_{t}, t=1, \cdots, 20)$ to drive the model. Let $\mathbf{X}_1=[3,3]^{tr}$. The top panel plots the average values
	 of $\mathbf{X}$ in different time periods over all these sample paths. As we increase $\delta$, the decay in signal $\mathbf{X}$ slows, indicating a stronger
	 autocorrelation in the information process. The remaining rows display the average quantities of assets that the trader needs to buy under the policies $\mathbf{S}^{DDP}$
	 and $\mathbf{S}^{LQ}$ across these 10,000 paths of $\mathbf{X}$ during each period.}
	\label{exm1:t3}
\end{figure}

Examining $\mathbf{S}^{DDP}$ will shed more insights into this improved policy. In Figure \ref{exm1:t3}, we compare it with $\mathbf{S}^{LQ}$, which is the policy derived from the simplified
auxiliary problem, by simulating them on the same sample paths of the signal process $\mathbf{X}_t$. Let us assume that the trader receives a large signal at $t=1$, i.e., the initial value
$\mathbf{X}_1$ is large. The top panel displays the evolution of the two-dimensional signal $\mathbf{X}_t=(X^{1}_t, X^{2}_t)$ over time under different autocorrelation coefficient $\delta$. The other three rows illustrate how these two policies respond to
the changes in $\mathbf{X}_t$ in terms of the respective purchase amounts of the three assets. We can see that, in response to the ``good" initial signal, the suboptimal strategy $\mathbf{S}^{LQ}$
(the red curves in Figure \ref{exm1:t3}) immediately increases its purchase. This behavior is economically sensible. Under our choice of $\mathbf{C}$ in the dynamic of (\ref{exm1:dyn3}), a high
current value of $\mathbf{X}$ implies that the prices of the assets are likely to move up in the future. To avoid the high transaction cost that the trader might pay consequently, she would like to buy more at the current price immediately.

However, this policy is suboptimal in the presence of such market frictions as the price impact and the no-sale constraint. The blue curves in the figure illustrate how the optimal policy
$\mathbf{S}^{DDP}$ should behave. Interestingly, it executes transactions much more slowly in response to the same signal $\mathbf{X}$ compared with $\mathbf{S}^{LQ}$.
Moreover, the autocorrelation of the signal process accentuates the difference in the trading speeds between the two policies. As we increase the value of $\delta$ from the left column to
the right in Figure \ref{exm1:t3}, $\mathbf{S}^{LQ}$ and $\mathbf{S}^{DDP}$ become distinct.  We also examine the effect of the temporary impact by changing $\lambda$ in the
experiments. The green curves are corresponding to the case in which $\lambda=100$. In comparison with the red curves ($\lambda=0$), the trader further smooths her transactions
in order to avoid the excessive costs associated with the temporary impact of trading (cf. the first term in (\ref{exm1:objective_new})).

\subsection{Inventory Management with Lost Sales and Lead Time}
\label{exm2}

In this section, we consider a single-item inventory management problem with stochastic demands, a constant lead time and lost sales. Assume that a manager has a finite planning horizon of $T$
periods. In period $t$, $t=1,2,\cdots, T$, a random demand amounting to $d_{t}$ will arise. All the demands across different periods are supposed to be independent and have the identical
distribution. The manager needs to use the current inventory to meet the demand in each period and meanwhile determines an amount of $a_{t}$ to order. Denote $L$ to be the order lead time;
that is, the order placed in period $t$ will arrive in period $t+L$. Hence, the manager's decision making should be based on a state vector of $L$ components
$\mathbf{x}_t=(x_{0,t},x_{1,t},\cdots, x_{L-1,t})$, where $x_{0,t}$ is the amount of the current inventory in period $t$ and $x_{l,t}$ is the order arriving in the subsequent periods $t+l$ for
$l=1, \cdots, L-1$. If the current inventory is not sufficient, we assume that the unfulfilled demands will be immediately lost. After receiving $x_{1,t}$ at the beginning of the next period, the inventory
level transits to $(x_{0,t}-d_{t})^++x_{1,t}$ and the manager starts a new decision-making loop. From the above discussion, we can easily see that the state vector for period $t+1$ should
be given by
\begin{eqnarray}
\label{exm2:dy}
\mathbf{x}_{t+1}=((x_{0,t}-d_t)^++x_{1,t},x_{2,t},\cdots,x_{L-1,t},a_t).
\end{eqnarray}
Note that this dynamic is not linear.

The manager faces three types of costs: procurement cost associated with orders, inventory holding cost, and lost-sale penalty. For notational simplicity, we ignore the first type of cost
in our model by letting the unit cost of procurement be $0$. As argued in \cite{jm}, this assumption will not hurt the generality of the setup. Let $h$ and $p$ denote the marginal cost of holding
inventory and the penalty of lost sales, respectively. Then the manager attempts to minimize the discounted total cost over $T+L$ periods, namely,
\begin{eqnarray}
\label{exm2:tg}
	\min_{ \substack{a_{t} \in \mathbb{Z}_{+},\\ 1\le t\le T+L}}\mathbb{E}\Big[\sum_{t=1}^{T+L}\gamma^t\left(h(x_{0,t}-d_t)^++p(d_t-x_{0,t})^+\right)\Big],
\end{eqnarray}
where
\[
q(\mathbf{x}_t,d_t):=h(x_{0,t}-d_t)^++p(d_t-x_{0,t})^+
\]
is the sum of the inventory cost and the lost-sale penalty in period $t$,  $\gamma \in (0, 1)$ is the discount factor used by the manager, and $\mathbb{Z}_{+}$ stands for the set of all nonnegative
integers.

The lost-sales model was first formulated in \cite{ks} and further explored in \cite{mor1,mor2}. It is well known that the model is intractable, especially for a large lead time $L$.
\cite{zpa, zpb} presents insightful structural analysis on this standard problem and, based on that, tests several plausible heuristics. He finds that the following myopic policy yields
analytical value functions and performs reasonably well. Rather than considering the entire time horizon, the myopic policy chooses the order quantity $a_{t}$ in period $t$
to minimize the cost from period $t$ to period $t+L$. That is, letting
\begin{eqnarray}
\label{exm2:myopic1}
	a^{\textrm{my}}_t=\arg\min_{a_{t} \in \mathbb{Z}_{+}}\mathbb{E}\left[\sum_{s=t}^{t+L}\gamma^{s-t}q(\mathbf{x}_s,d_s)\right].
\end{eqnarray}
Note that the order $a_{t}$ arrives in period $t+L$ and has nothing to do with the inventory prior to that period. Thus we can easily show that the optimization in (\ref{exm2:myopic1}) is equivalent to
\begin{eqnarray}
\label{exm2:myopic2}
	a^{\textrm{my}}_t=\arg\min_{a_{t} \in \mathbb{Z}_{+}}\mathbb{E}[\gamma^{L}q(\mathbf{x}_{t+L},d_{t+L})].
\end{eqnarray}
This policy apparently neglects the evolution of the inventory system after period $t+L$.

Relatedly, \cite{cdj}  develop a new numerical approach to approximate the optimal value function of this example using a selected number of points in a bounded rectangular domain.
Their method hinges on the $L^{\sharp}$-convex property of the value function. \cite{bgy} analyze the asymptotic optimality of a given heuristic in an infinite-horizon lost-sales inventory
model with positive lead time. \cite{bs} apply the information relaxation based dual method to assess the above myopic policy.

\medskip

\noindent The following numerical experiments test the performance of the DDP algorithm by using it to assess and improve several heuristic policies. We assume that the stochastic demand $d_t$
follows a geometric distribution with mean $m$. As pointed out by \cite{zpa}, this distribution is more likely to produce extreme demand scenarios. Two possible lead times, $L=4$ and $L=10$, are
considered. As in the previous optimal execution problem, we need to choose a proper state selector $G$ to sample the representative states $\mathbf{x}_t$ in each period $t$. Let $\theta = h/ (p+h)$ and define
\[
s_l=\min\left\{s: \mathbb{P}\left(\sum_{m=l}^{L}d_m>s\right)\le\theta \right\}
\]
for $l=0, \cdots, L-1$. Both \cite{mor1} and \cite{zpa,zpb} show that, starting with initial state $\mathbf{x}_1=0$, the inventory process under the optimal policy will never leave the region
\begin{eqnarray}
\label{exm2:state_redu}
\mathcal{X}_t=\left\{\mathbf{x}_t\ge 0: \sum_{m=l}^{L-1}x_{m,t}\le s_l,\ l=0,\cdots,L-1 \right\}.
\end{eqnarray}
In light of these results, we take $G$ to be the discrete uniform distribution over the compact set $\mathcal{X}$.

In the interest of space, we defer the explicit expressions of all the basis functions used in this section to Appendix \ref{app:num}. To evaluate the penalty function, we need to
calculate the expectations of these basis functions. This step may be computationally expensive when $L$ is large. As mentioned in Section \ref{sec:mc}, we suggest using low-discrepancy
sequences from the QMC literature to develop effective approximations. A detailed explanation of this approximation can also be found in \ref{app:num}.

Along each sample path of demand $\mathbf{d}_t=(d_t,d_{t+1},\cdots,d_{T+L})^{tr}$, the DDP algorithm solves the following deterministic inner optimization problem
\begin{eqnarray}
\label{exm2:optimization}
	J(\mathbf{x}_{t}, \mathbf{d}_{t}):=\inf_{ a\in \mathbb{Z}^{T-t+1}_{+}}\sum_{s=0}^{T+L-t}\Big\{\gamma^sq(\mathbf{x}_{t+s},d_{t+s},a_{t+s})+z_t(a,\mathbf{d}_t)\Big\}
\end{eqnarray}
at each time step $t$ for the dual value determination. It can be reduced to an integer DC program. \cite{mmm} employ a special form of continuous relaxation
(known as ``lin-vex extension" in their paper) to find an exact solution to DC optimization programs with integer constraints. However, to save the computational effort, we
take an alternative approach here by simply relaxing the integer constraint
$a\in \mathbb{Z}^{T-t+1}_{+}$ to $a \ge 0$ when solving (\ref{exm2:optimization}). The relaxation enables us to apply the sequential-convex-programming method in \ref{app:dc} to obtain
a lower bound for $J(\mathbf{x}_{t}, \mathbf{d}_{t})$. The numerical experiments show that the convergence of the DDP algorithm is not affected by this continuous relaxation.

\medskip

\noindent Table \ref{exm2:t1} displays the performance of our DDP algorithm in improving some heuristic policies. In addition to the myopic policy given in (\ref{exm2:myopic22}), we also consider
several alternative approximate policies as follows:
\begin{itemize}
\item Lookahead. The above myopic policy ignores the long run impact of the current order. To remedy this, we may introduce a $\tilde{V}$ to approximately capture the future impact of the
present order. A lookahead policy stems from solving
\begin{eqnarray}
\label{exm2:myopic22}
	a^{\textrm{LA}}_t=\arg\min_{a_{t} \in \mathbb{Z}_{+}}\mathbb{E}\left[\sum_{s=t}^{t+L}\gamma^{s-t}q(\mathbf{x}_s,d_s)+\gamma^{L+1}\tilde{V}(\mathbf{x}_{t+L+1},d_{t+L+1})\right].
\end{eqnarray}
In the experiment, we try the total cost function at time $t+L+1$, $q(\mathbf{x}_{t+L+1},d_{t+L+1})$, as the approximation $\tilde{V}$.
\item Linear programming approximation. We omit the details here because the idea is similar to the LP approximation in the previous example.
\end{itemize}

\begin{table}[htbp]
	\centering
	\renewcommand\arraystretch{1}
	\begin{tabular}{cccccccc}
		\toprule
		& &  & L=4  & & &  \\
		\midrule
		Approximation & Iteration & Dual Values & (SE)  & Primal Values & (SE)   & Gap  \\
		\midrule
		& 0     &       &              & 563.72  & (0.42)  & 4.36\% \\
		Myopic  & 1     & 539.16 & (0.38)   & -  & -  & - \\
		& 2     & 539.86 & (0.09)   & 542.00  & (0.43)  & 0.39\% \\
		& 3     & 539.88 & (0.08)   &       &       &  \\
		\midrule
		& 0     &       &              & 560.13  & (0.41)  & 3.63\% \\
		Lookahead & 1     & 539.78 & (0.36)   & -  & -  & - \\
		& 2     & 539.88 & (0.08)   & 542.13  & (0.43)  & 0.41\% \\
		& 3     & 539.89 & (0.08)   &       &       &  \\

		\midrule
		& 0     &       &              & 566.31  & (0.50)  & 5.00\% \\
		Linear & 1     & 537.86 & (0.40)   & -  & -  & - \\
		& 2     & 539.80 & (0.08)   & 542.08  & (0.40)  & 0.41\% \\
		& 3     & 539.84 & (0.08)   &       &       &  \\
		\bottomrule
		\toprule
		& &  & L=10  & & &  \\
		\midrule
		Approximation & Iteration & Dual Values & (SE)  & Primal Values & (SE)   & Gap  \\
		\midrule
		& 0     &       &              & 829.63  & (0.28)  & 7.36\% \\
		Myopic   & 1     & 768.58 & (0.36)   & -  & -  & - \\
		 & 2     & 770.93 & (0.08)   & -  & -  & - \\
		& 3     & 771.80 & (0.08)   & 779.36  & (0.29)  & 0.96\% \\
		& 4     & 771.89 & (0.07)   &       &       &  \\
		\midrule
		& 0     &       &              & 827.14  & (0.30)  & 7.04\% \\
		Lookahead  & 1     & 768.91 & (0.36)   & -  & -  & - \\
		& 2     & 771.10 & (0.08)   & -  & -  & - \\
		& 3     & 771.83 & (0.08)   & 779.54  & (0.29)  & 0.99\% \\
		& 4     & 771.84 & (0.08)   &       &       &  \\
		\midrule
		& 0     &       &              & 844.28  & (0.30)  & 9.50\% \\
		Linear  & 1     & 764.10 & (0.35)   & -  & -  & - \\
		& 2     & 770.59 & (0.09)   & -  & -  & - \\
		& 3     & 771.61 & (0.08)   & -  & -  & - \\
		& 4     & 771.89 & (0.08)   & 779.80  & (0.30)  & 1.02\% \\
		& 5     & 771.83 & (0.08)   &       &       &  \\
		\bottomrule
	\end{tabular}%
  \captionsetup{font=scriptsize}
	\caption{Convergence results of the DDP algorithm under different initial policies in the example of inventory management. The default parameters used in this experiment are
	$\{m=4,\ h=1,\ p=9,\ T=30\}$. Two lead times are implemented, i.e., $L=4$ and $L=10$. The three types of heuristic policies used as the inputs are Myopic, Lookahead, and Linear.
	We sample $500$ states for $L=4$ and $1000$ states for $L=10$ in each time step from the distribution $G$ to compute the dual values. In LP, we simulate 100 states to set up
	the corresponding constraints. All computation experiments are conducted on a PC equipped with an Intel Xeon 32-core 2.93 GHz CPU and 12.0 GB of RAM. The computation
	environment is Windows 7 and MATLAB R2017a and parallel pool. The average computational time is 323.6s per iteration for $L=4$ and 2365.1s for $L=10$.}
	\label{exm2:t1}%
\end{table}%

It is worth mentioning that this problem is solvable through the associated Bellman equation when the lead time $L=4$. Using (\ref{exm2:state_redu}), the total number of the states that we need to
visit in each period in this case is 60,129.  By brute-force searching for the best order quantities in all these states, we find that the true optimal value of the problem when $L=4$ should be $541.82$
under the parameter values we set up for the experiment. The structure of Table \ref{exm2:t1} remains similar to that of Table \ref{exm1:t1}. We can see that the DDP algorithm manages to
significantly reduce down the duality gaps of all these heuristics.

\begin{figure}[htbp]
	\centering
	\subfigure{
		\includegraphics[width=5cm]{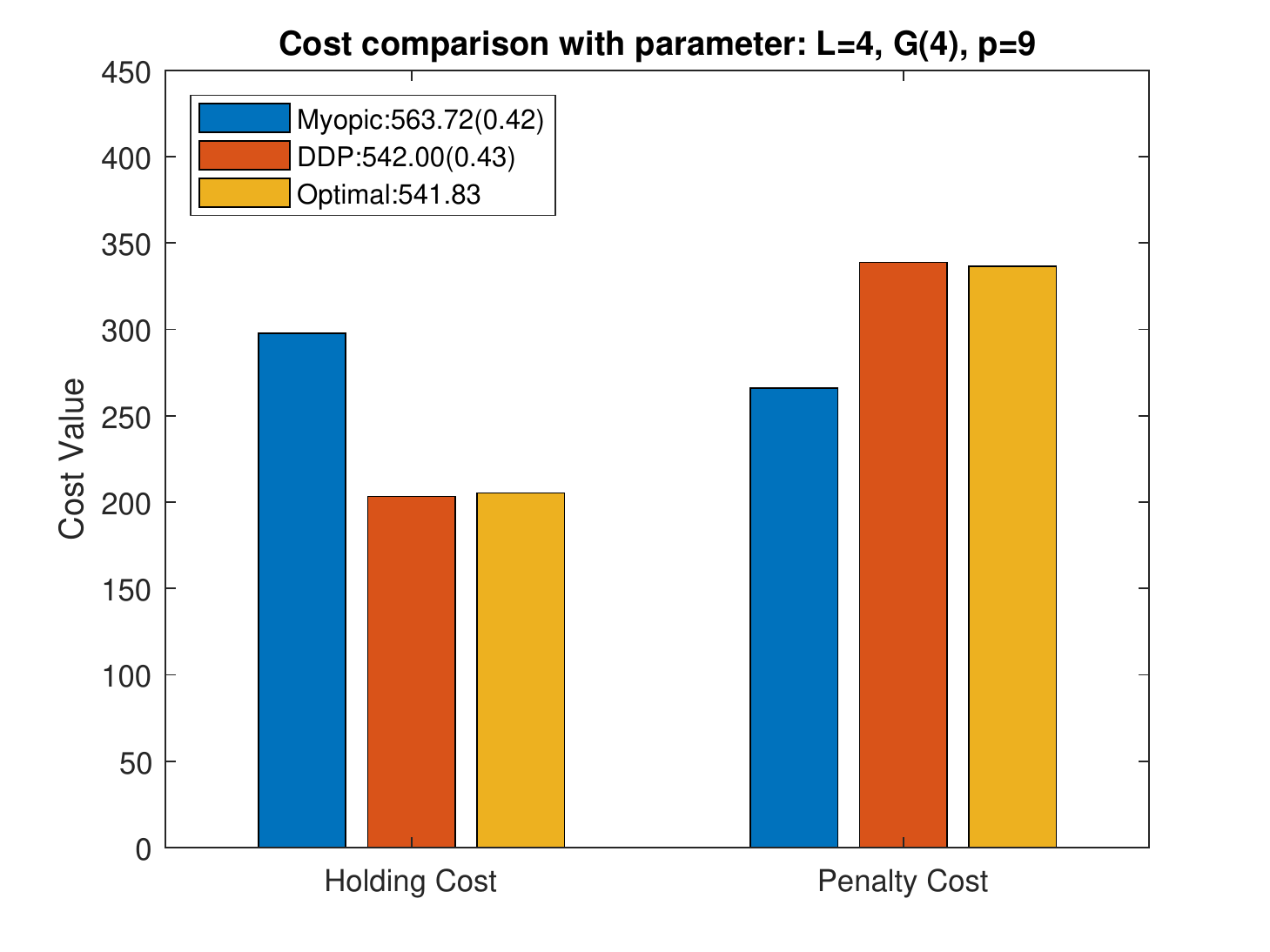}
	}
	\subfigure{
		\includegraphics[width=5cm]{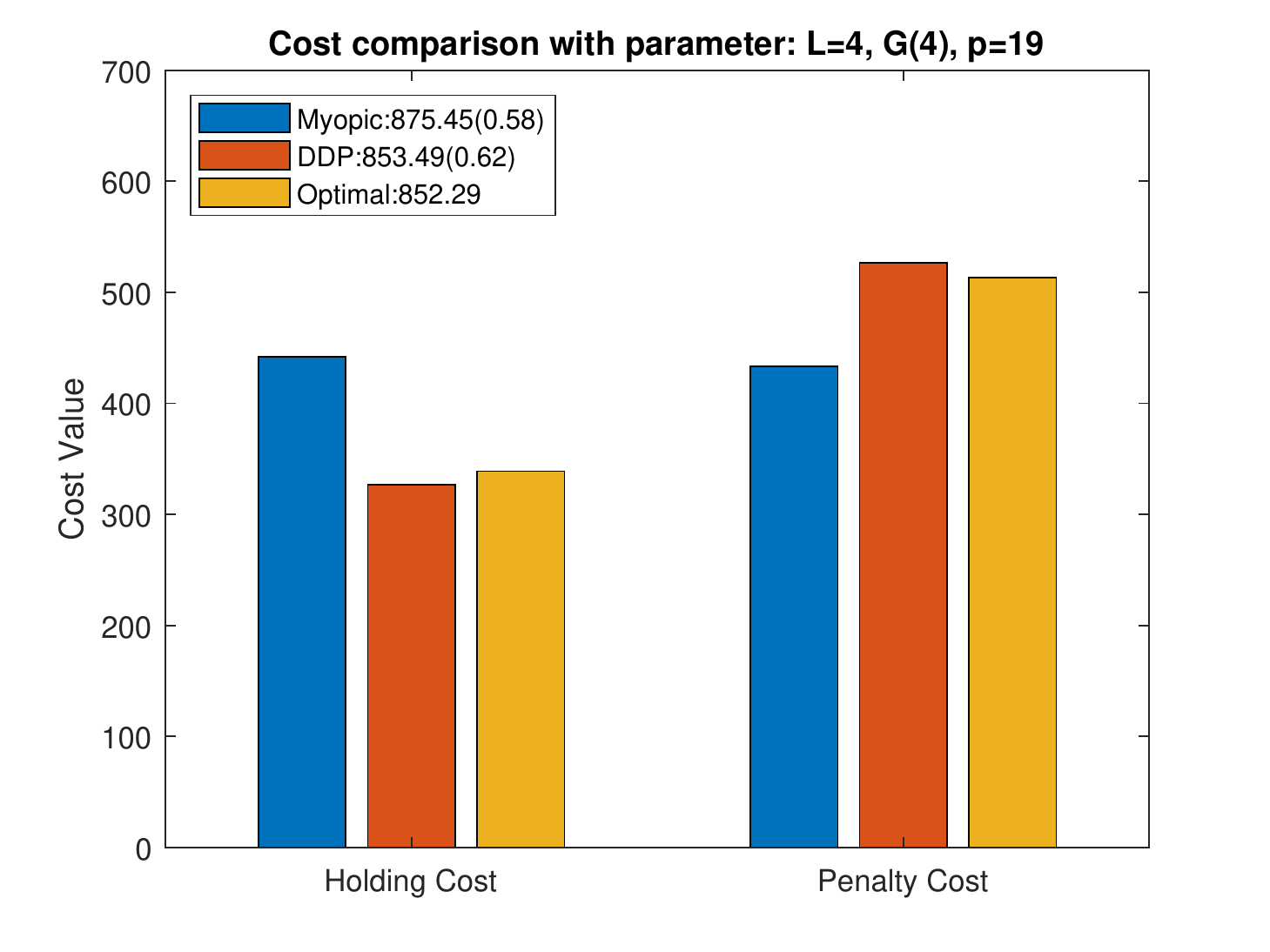}
	}
	\\
	\subfigure{
		\includegraphics[width=5cm]{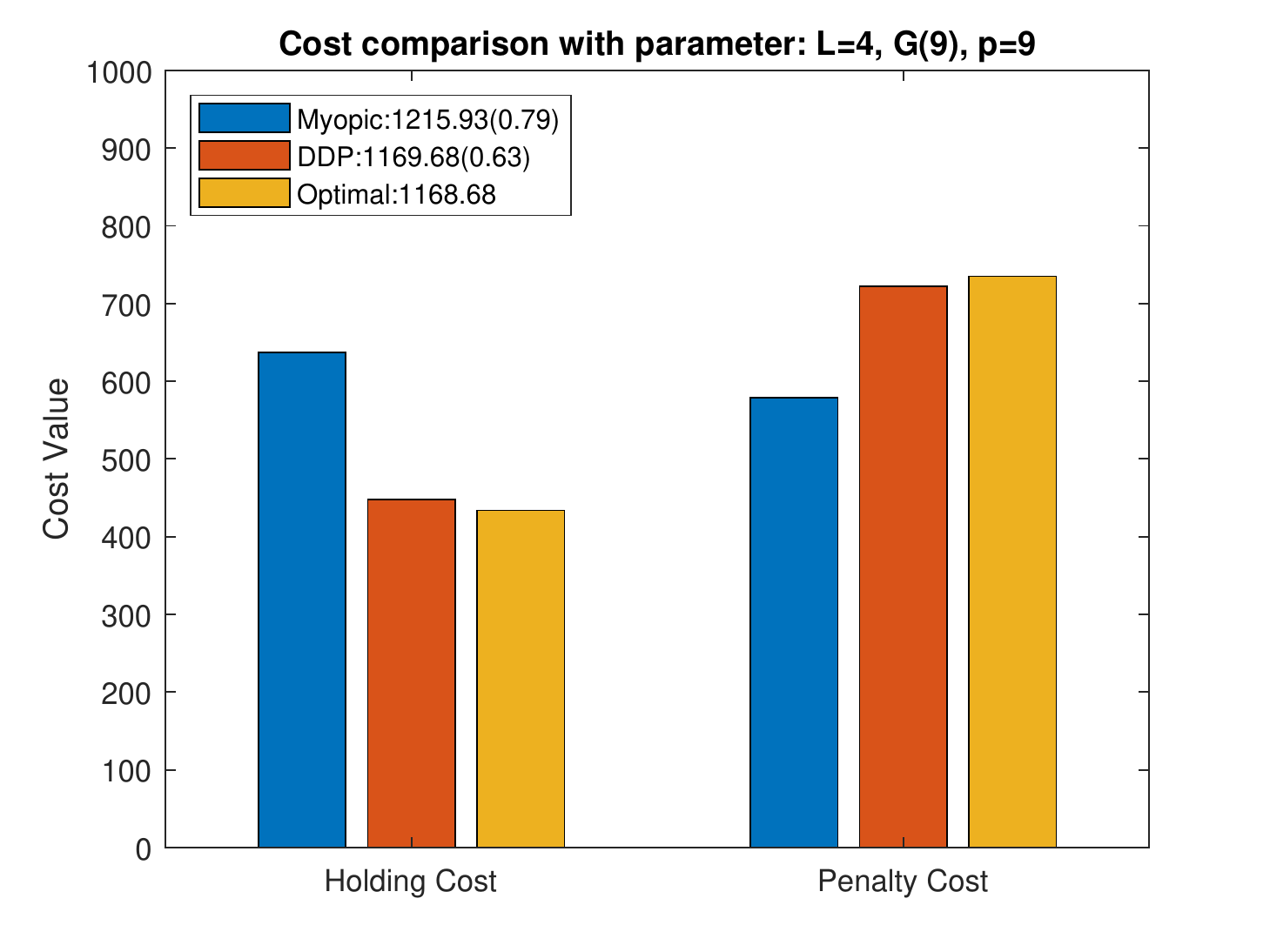}
	}
	\subfigure{
		\includegraphics[width=5cm]{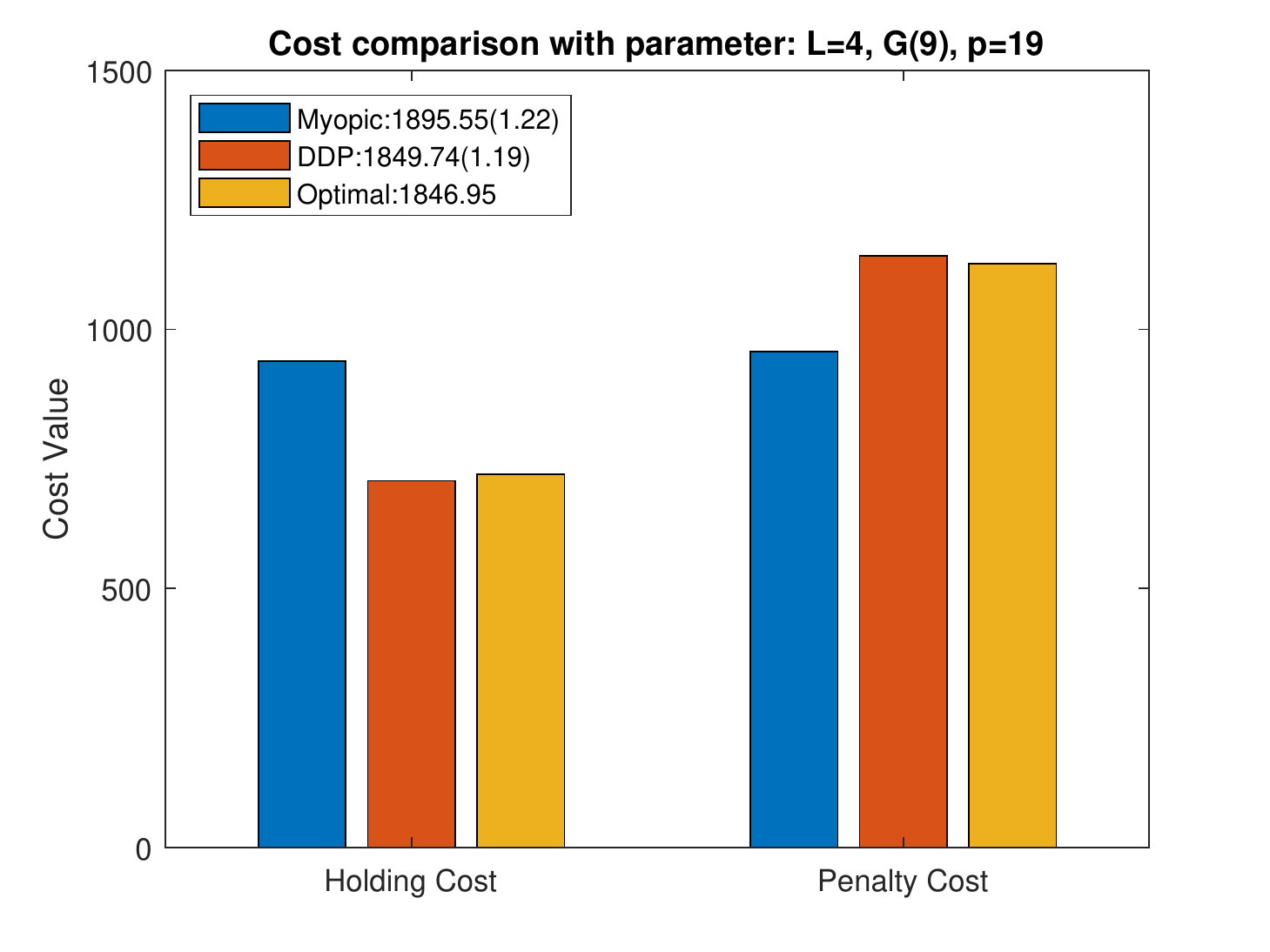}
	}

	\caption{Cost comparison in $L=4$. We compare the holding and penalty costs under three policies: the myopic policy $a^{my}$, the policy improved from the DDP algorithm
	$a^{DDP}$, and the optimal one $a^{op}$. We simulate 10,000 sample paths of random demands and evaluate both $a^{my}$ and $a^{DDP}$ based on the same set of sample paths.
	We use brute-force searching to solve the Bellman equation to obtain the optimal value for $a^{op}$. $G(4)/G(9)$ stands for the geometric random demand with mean 4/9.
	In the legend of each subfigure, we report the total cost of each policy and the corresponding standard error in the brackets. Note that $a^{DDP}$, the resulted policy from our
	DDP algorithm, behaves exactly the same as the optimal one.}
	\label{exm2:t4}
\end{figure}

\begin{figure}[htbp]
	\centering
	\subfigure{
		\includegraphics[width=5cm]{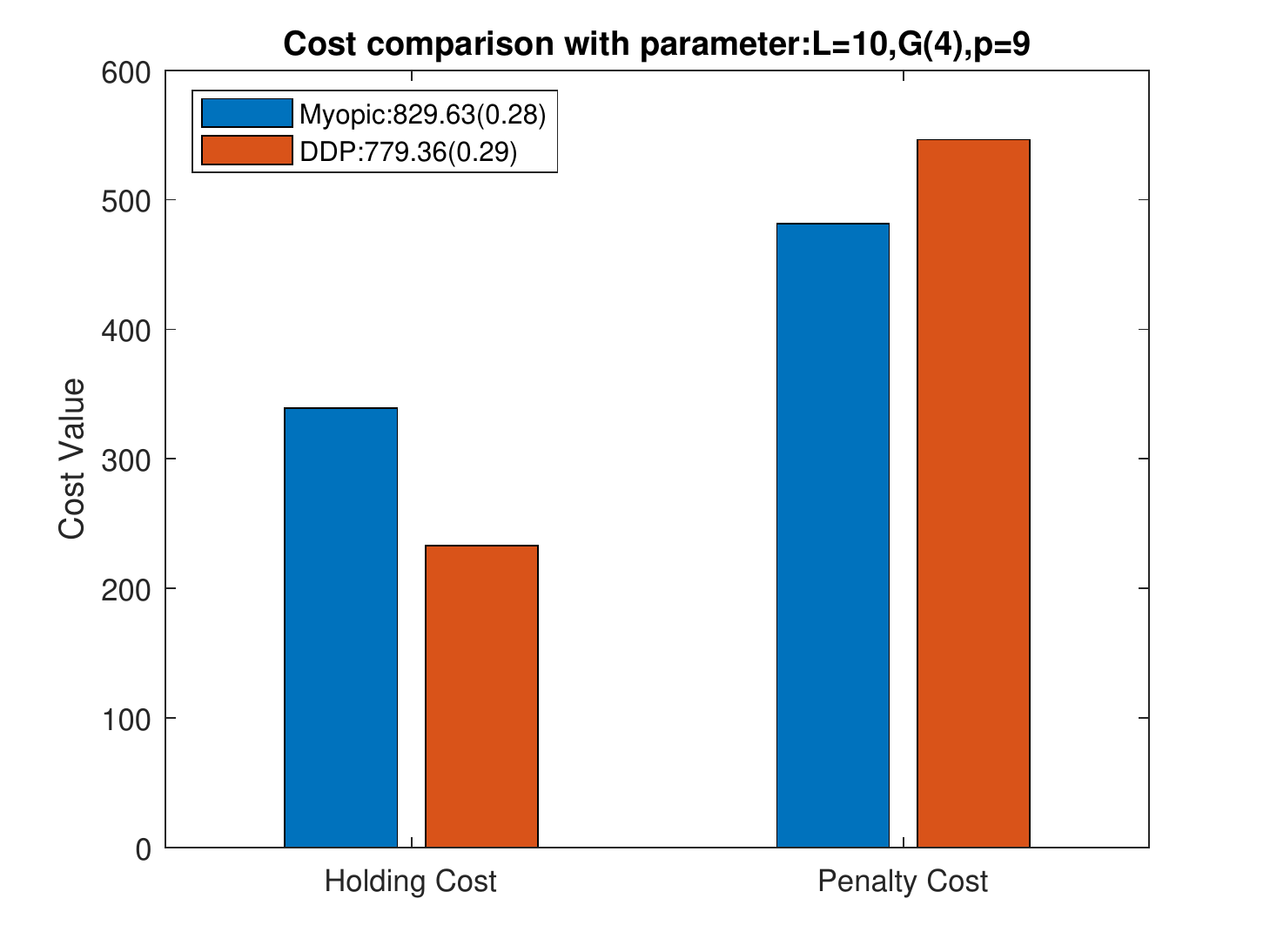}
	}
	\subfigure{
		\includegraphics[width=5cm]{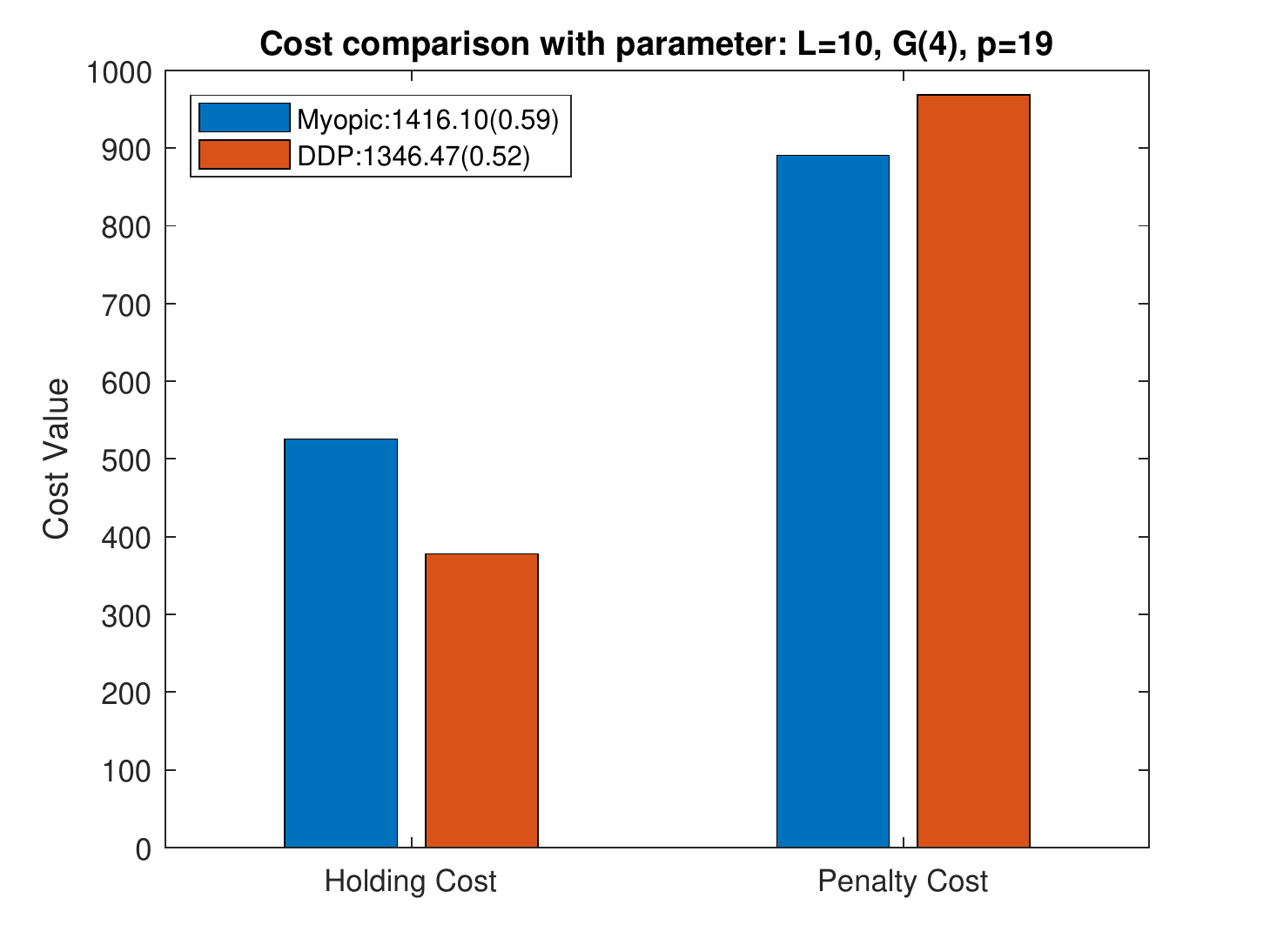}
	}
	\\
	\subfigure{
		\includegraphics[width=5cm]{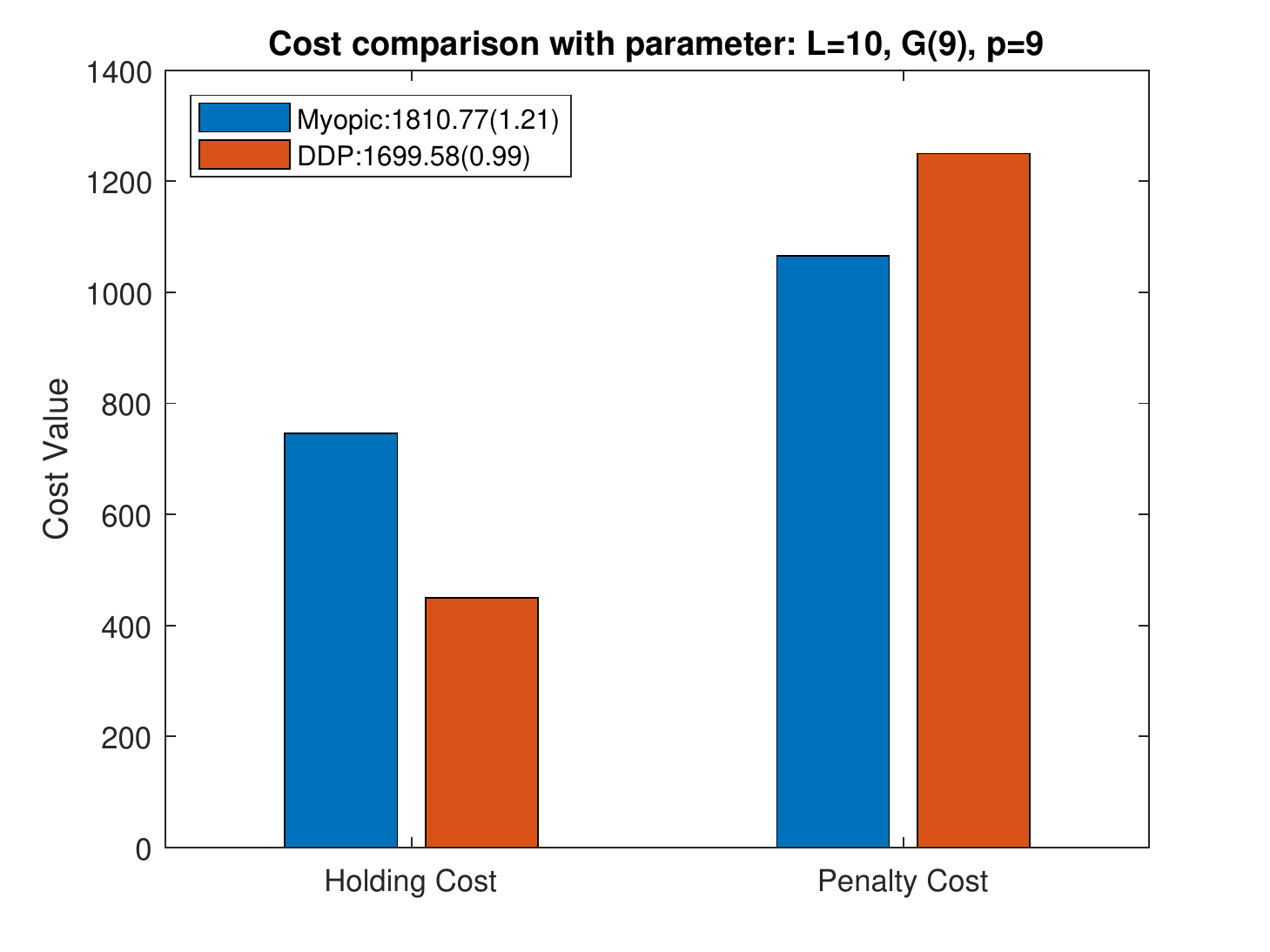}
	}
	\subfigure{
		\includegraphics[width=5cm]{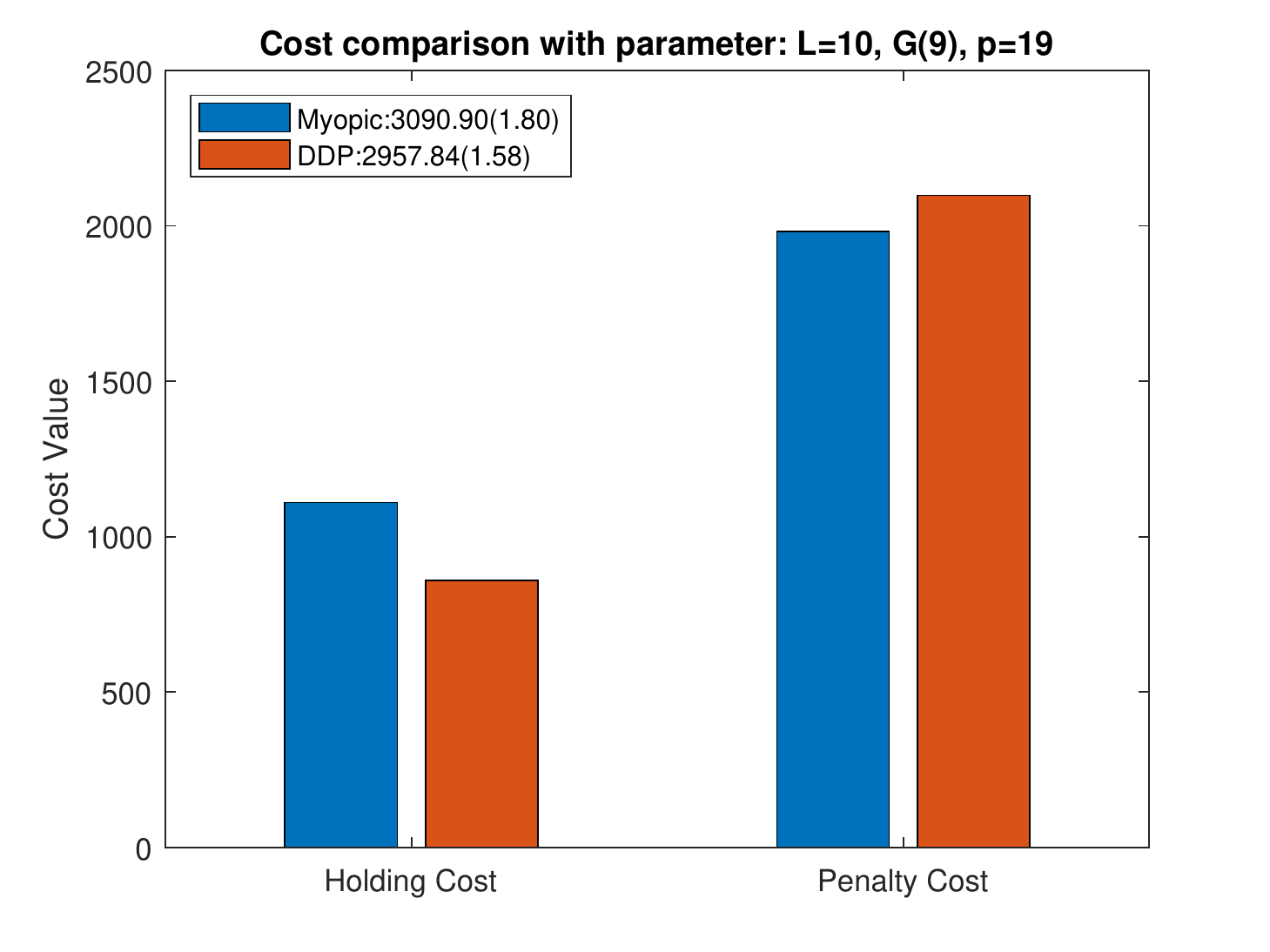}
	}

	\caption{Cost comparison in $L=10$. Note that we do not report the costs associated with the optimal policy because it is impossible to apply the Bellman equation to solve for the optimal solution
	due to the high dimensionality.}
	\label{exm2:t5}
\end{figure}

\begin{figure}[htbp]
		\centering
		\includegraphics[width=2.7in]{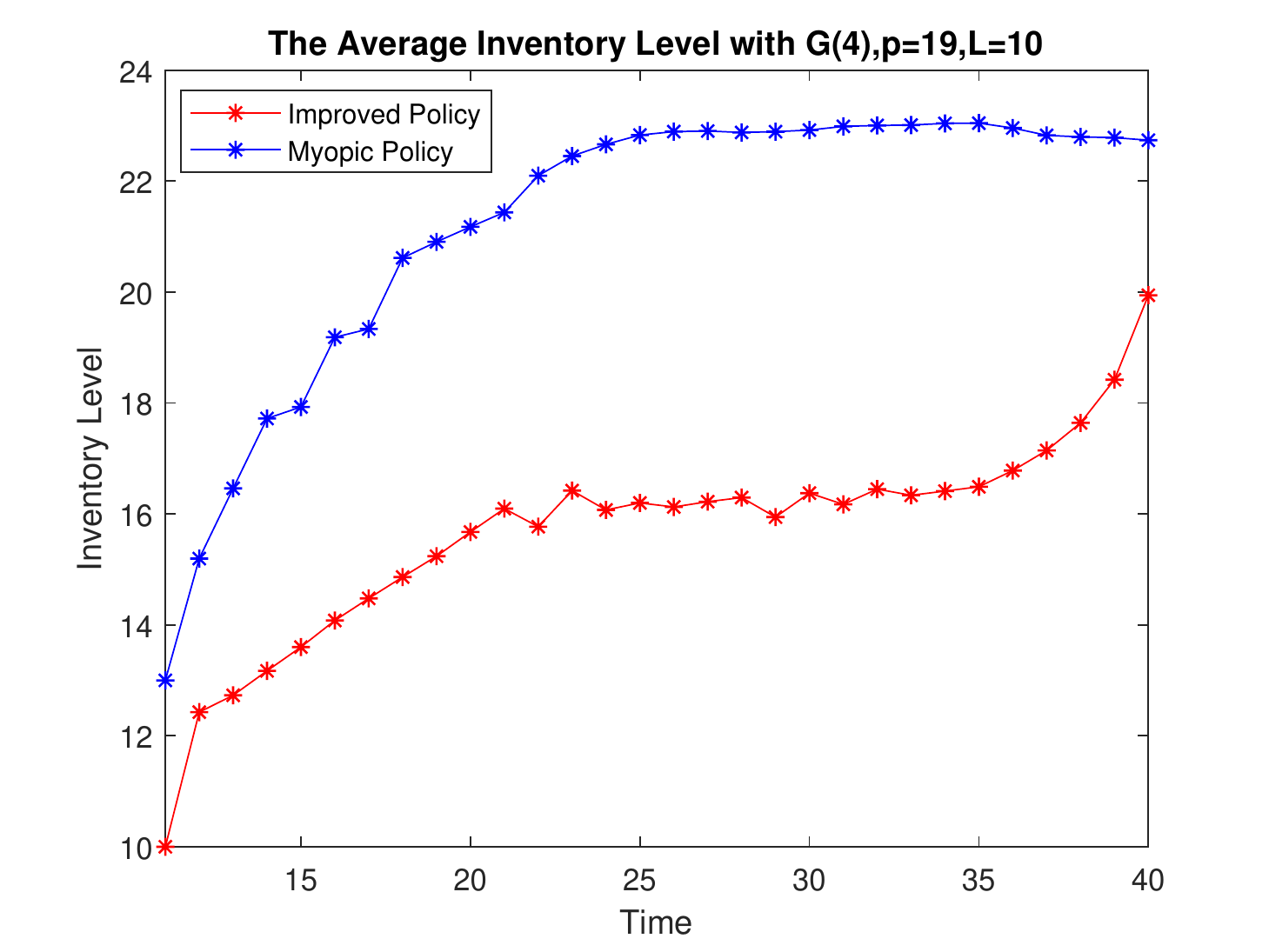}
		\caption{Average inventory level for a system controlled by myopic policy $a^{my}$ and improve policy $a^{DDP}$. The parameter settings are: $L=10$, $p=19$, $T=30$ and the
		geometric demand distribution with mean 4. In this figure we sample 10,000 random  demands $\{d_t, 1\le t\le T+L\}$. Under the same sample path, we run the myopic and improved
		policies. The two curves in the figure display the average inventory level at each time step across all the sampled demands.}
		\label{exm2:p1}
\end{figure}

Figures \ref{exm2:t4} and \ref{exm2:t5} help us gain more insight about where the improvement of the policy that the DDP finally converges to comes from. We simulate the inventory system under both
the myopic policy and the policy obtained from the DDP method. The two types of costs, the inventory holding cost
\begin{eqnarray*}
\label{exm2:cost1}
	\mathbb{E}\Big[\sum_{t=1}^{T+L}\gamma^t h(x_{0,t}-d_t)^+\Big]
\end{eqnarray*}
and the lost sales penalty cost
\begin{eqnarray*}
\label{exm2:costs}
	\mathbb{E}\Big[\sum_{t=1}^{T+L}\gamma^t p(d_t-x_{0,t})^+\Big],
\end{eqnarray*}
are calculated and compared in the figure. The myopic policy focuses on the short-term performance and neglects the long-run impact of orders on the inventory level. Therefore it incurs a smaller
lost-sale penalty than the optimal one. However, this comes at the expense of the inventory holding cost. In contrast, the improved policy that resulted from the DDP algorithm strikes a better balance
between these two costs. The inventory cost under it is smaller than what the myopic policy causes, which leads to a better overall cost performance. Figure \ref{exm2:p1} further compares the
average inventory level for a system controlled by both policies and subject to the same demand shocks over the time horizon. It clearly demonstrates that the system tends to build up more
inventory, thus incurring more holding costs, if the manager uses the myopic policy.

\section{Conclusions}

In this paper we present a duality-driven iterative approach (DDP) for solving a general SDP problem. The duality gap yielded by the method can be used to assess the performance of a given policy.
More importantly, repeatedly applying the dual operation on the basis of the technique of information relaxation will lead to policy improvement and convergence to the optimality. To implement the
DDP, we also develop a regression Monte Carlo method. In conjunction with such techniques as DC programming and parallel computing, our method demonstrates numerical
effectiveness and accuracy in dealing with multidimensional complex SDP problems.

% Acknowledgments here

\bigskip
\noindent \textbf{Acknowledgments}: This research is supported by the Research Grant Council Hong Kong through the scheme of General Research Fund (Grant No. 14237616 and 14207918), and the National Natural Science Foundation of China (Grant No. 71991474, 71721001, and U1811462).

\newpage
\appendix
\centerline{\huge{\textbf{Appendix}}}
\section{Proofs of Main Results in Section \ref{sec:value_iteration}}
\label{app:proof}

Before proving Proposition \ref{pro:super}, we establish a Bellman equation-like characterization of the value of the inner optimization problem
in the dual formulation. Let
\begin{eqnarray*}
\mathfrak{J}_{t}(\xi | t, x_t)= \inf_{ a \in A|t}
	\left[\sum_{s=t}^{T-1} r_s(x_s,a_s,\xi_s)+ r_T(x_T)+z_t(a, \xi) \right]
\end{eqnarray*}
with $z_t(a,\xi)$ being defined as in (\ref{penalty_def}). We have

\begin{lem}
\label{bellman_equation_dual}
For $0\le t\le T-1$,
\begin{eqnarray*}
\mathfrak{J}_{t}(\xi | t, x_t)= \inf_{a_t \in A_t}\Big\{\mathbb{E}[r_t(x_t,a_t,\xi_t)+W_{t+1}(x_{t+1})|x_{t}=x]-W_{t+1}(x_{t+1})+\mathfrak{J}_{t+1}(\xi | t+1, x_{t+1}) \Big\}.
\end{eqnarray*}
\end{lem}
\textit{Proof of Lemma \ref{bellman_equation_dual}.} Note that $\mathfrak{J}_{t}(\xi | t, x_t)$ admits the following representation
\begin{eqnarray*}
	\mathfrak{J}_{t}(\xi | t, x_t)&=&\inf_{a_t \in A_t} \Big\{ r_t(x_t,a_t,\xi_t)+\mathbb{E}[r_s(x_s,a_s,\xi_s)+W_{s+1}(f_s(x_s,a_s,\xi_s))]-(r_s(x_s,a_s,\xi_s) \\
  &&+W_{s+1}(f_{s}(x_s,a_s,\xi_s)))
	+\inf_{ a \in A|(t+1)} \Big[\sum_{s=t+1}^{T-1} r_s(x_s,a_s,\xi_s) +r_T(x_T)+  z_{t+1}(a,\xi)\Big] \Big\}.
\end{eqnarray*}
The conclusion trivially follows. $\hfill\square$
\medskip

\noindent\textit{Proof of Proposition \ref{pro:super}.}
By the definition of the dual operator $\mathcal{D}$, we have
\begin{eqnarray}
\label{ec:pro_1}
\mathcal{D}W_t(x)=\mathbb{E}\left[\inf_{ a \in A|t}
	\left[\sum_{s=t}^{T-1} r_s(x_s,a_s,\xi_s)+ r_T(x_T)+z_t(a, \xi) \right]\Big|x_t=x\right]=\mathbb{E}\Big[\mathfrak{J}_{t}(\xi | t, x_t)\Big|x_t=x\Big].
\end{eqnarray}
According to Lemma \ref{bellman_equation_dual}, for any action $a_t \in A_t$,
\[
\mathfrak{J}_{t}(\xi | t, x_t) \le \mathbb{E}[r_t(x_t,a_t,\xi_t)+W_{t+1}(x_{t+1})|x_{t}]-W_{t+1}(x_{t+1})+\mathfrak{J}_{t+1}(\xi | t+1, x_{t+1}).
\]
Therefore,
\begin{eqnarray}
\label{ec:pro_2}
\resizebox{.95\hsize}{!}{$
 \mathbb{E}\Big[\mathfrak{J}_{t}(\xi | t, x_t)\Big|x_t=x\Big]\le \mathbb{E}\Big[\mathbb{E}[r_t(x_t,a_t,\xi_t)+W_{t+1}(x_{t+1})|x_{t}]-W_{t+1}(x_{t+1})+\mathfrak{J}_{t+1}(\xi | t+1, x_{t+1})\Big|x_t=x\Big].
 $}
\end{eqnarray}

Note that, by the iterated law of conditional expectation, we have
\[
\mathbb{E}\big[\mathbb{E}[W_{t+1}(x_{t+1})|x_t]-W_{t+1}(x_{t+1})|x_{t}=x\big]=0.
\]
Moreover,
\[
\mathbb{E}[\mathfrak{J}_{t+1}(\xi | t+1, x_{t+1})|x_t]=\mathbb{E}[\mathbb{E}[\mathfrak{J}_{t+1}(\xi | t+1, x_{t+1})]|x_t=x]=\mathbb{E}[\mathcal{D}W_{t+1}(x_{t+1})|x_t=x].
\]
Both equalities lead to that the right hand of (\ref{ec:pro_2}) is equal to
\begin{eqnarray*}
\mathbb{E}[r_t(x_t,a_t,\xi_t)+\mathcal{D}W_{t+1}(x_{t+1})|x_t=x]
\end{eqnarray*}
In conjunction with (\ref{ec:pro_1}), we have
\begin{eqnarray*}
\mathcal{D}W_t(x) \le \inf_{a_t \in A_t}  \mathbb{E}[r_t(x,a_t,\xi_t) +\mathcal{D}W_{t+1}(f_t(x,a_t,\xi_t))]. \hfill\square
\end{eqnarray*}

\medskip

\noindent\textit{Proof of Theorem \ref{thm:convergence}.}
(i) To show this, we need the following two claims:
\begin{itemize}
\item  For any sequence $W=(W_0,W_1,...,W_T)$, $\mathcal{D}^nW$ is a subsolution for any $n\ge1$.
\item For a subsolution sequence $W=(W_0,W_1,...,W_T)$, $W_{t}: \mathbb{R}^n \rightarrow \mathbb{R}$,  we have
\begin{eqnarray*}
W_t(x)\le(\mathcal{D}W)_{t}(x) \text{ for } 0 \le t \le T-1
\end{eqnarray*}
and  $W_T(x)=(\mathcal{D}W)_{T}(x)=r_{T}(x)$.

\end{itemize}

According to Proposition \ref{pro:super}, $\mathcal{D}^nW$ is a subsolution for any $n\ge1$. Now we turn to the second claim. It suffices to prove that, for any sequence of subsolution $W$, we have $\mathfrak{J}_{t}(\xi | t, x_t)\ge W_t(x)$ for all $t \ge 0$ with $\mathfrak{J}_{t}(\xi | t, x_t)$ defined in Lemma \ref{bellman_equation_dual}. Indeed, invoking (\ref{ec:pro_1}),
\[
\mathcal{D}W_{t}(x) = \mathbb{E}\Big[\mathfrak{J}_{t}(\xi | t, x_t)\Big|x_t=x\Big] \ge W_{t}(x).
\]

We prove the claim of $\mathfrak{J}\ge W$ by performing induction on $t$. For $t=T$, it is clearly true since $\mathfrak{J}_T= r_T = W_{T}$ by the definition.
Suppose for $s\ge t+1$, the claim $\mathfrak{J}_s\ge W_s$ holds. Then at time $t$, according to Lemma \ref{bellman_equation_dual},
\begin{eqnarray*}
\mathfrak{J}_{t}(\xi | t, x_t)=&& \inf_{a_t \in A_t}\Big\{\mathbb{E}[r_t(x_t,a_t,\xi_t)+W_{t+1}(x_{t+1})]-W_{t+1}(x_{t+1})+\mathfrak{J}_{t}(\xi | t, x_{t+1}) \Big\}\\
\ge &&\inf_{a_t \in A_t}\mathbb{E}[r_t(x_t,a_t,\xi_t)+W_{t+1}(x_{t+1})].
\end{eqnarray*}
As $W$ is a subsolution sequence, we have
\begin{eqnarray*}
\mathfrak{J}_{t}(\xi | t, x_t)\ge \inf_{a_t \in A_t}\mathbb{E}[r_t(x_t,a_t,\xi_t)+W_{t+1}(x_{t+1})]\ge W_t(x_t).
\end{eqnarray*}
That completes the induction loop.
\smallskip

\noindent (ii) Consider a subsolution sequence $W$ such that $\mathcal{D}W=W$. We use induction again to prove this part.
At the last period $T$, we know that $(\mathcal{D}W)_{T}(x)=r_{T}(x)$ for all $x$ according to the definition of the operator $\mathcal{D}$.
Hence,
\[
W_T(x)=(\mathcal{D}W)_T(x)=r_T(x)=V_T(x).
\]
In words, $W$ and $V$ coincide at $T$. Now we assume this claim also holds for $s \ge t+1$ for some $t$. By Proposition \ref{pro:super},
the sequence $\mathcal{D}W$ constitutes a subsolution. Therefore,
\begin{eqnarray*}
(\mathcal{D}W)_t(x)&\ge&\inf_{a_t \in A_t}\mathbb{E}[r_t(x_{t},a_t,\xi_t)+(\mathcal{D}W)_{t+1}(f_t(x_t,a_t,\xi_t))|x_{t}=x]\\
&=&\inf_{a_t \in A_t}\mathbb{E}[r_t(x,a_t,\xi_t)+W_{t+1}(f_t(x,a_t,\xi_t))].
\end{eqnarray*}
In addition, the right hand side of the above inequality equals to, by the induction hypothesis,
\begin{eqnarray*}
\inf_{a_t \in A_t}\mathbb{E}[r_t(x,a_t,\xi_t)+V_{t+1}(f_t(x,a_t,\xi_t))]=V_t(x).
\end{eqnarray*}
where the last equality is due to the fact that $V$, as the true value function, should satisfy the Bellman equation. These lead to $(\mathcal{D}W)_t(x) \ge V_t(x)$.
On the other hand, $(\mathcal{D}W)_t(x) \le V_t(x)$ because of the weak duality property. Thus, we should have $(\mathcal{D}W)_t(x) \le V_t(x)$, which completes the induction loop.

\smallskip

\noindent (iii) First, we claim that, if for some $n$ and $t$ the equality $(\mathcal{D}^nW)_{t+1}(x)=V_{t+1}(x)$ holds for all $x$, then
\[
(\mathcal{D}^{k}W)_{t}(x)=V_{t}(x),
\]
for any $k\ge n+1$. This claim can be easily proved by induction. We omit the detail here for the interest of space. Once this claim is established, noting that $(\mathcal{D}^{1}W)_{T}(x)=V_{T}(x)$ is true,
it is easy to see that $(\mathcal{D}^{2}W)_{T-1}(x)=V_{T-1}(x)$ must be true for all state $x$. Using the above claim again, we can reach the following conclusion:
\[
(\mathcal{D}^{3}W)_{T-2}(x)=V_{T-2}(x)\quad \textrm{and}\quad (\mathcal{D}^{3}W)_{T-1}(x)=V_{T-1}(x).
\]
Repeatedly using the above argument leads to, for a general $k$,
 \[
(\mathcal{D}^{k}W)_{t}(x)=V_{t}(x)\ \textrm{for}\ t \ge T+1-k.
\]
In particular, when $k=T+1$, we have $(\mathcal{D}^{T+1}W)_{t}(x)=V_{t}(x)$ for $t \ge 0$. The theorem statement is proved. $\hfill\Box$

\section{The DDP Method in LQC}
\label{LQ_iteration}

We need the following technical lemma in calculating the duality of LQC.
\begin{lem}
	\label{LQ}
	We consider the quadratic programming
	\[
	J_{t}=\sum_{s=t}^{T-1}\left(x^{tr}_s Q_{s} x_s +a^{tr}_{t} R_{t}a_{t}+2\alpha_t^{tr}x_t+2\beta_t^{tr}a_t\right) + x^{tr}_T Q_{T} x_T,\]
	with the equality constraints
	\[x_{t+1}=D_t x_t+B_t a_t+\xi_t, \]
	where $\xi_t\in \mathbb{R}^n,\ \alpha_t\in \mathbb{R}^m,\ \beta_t\in \mathbb{R}^m$ are given vectors, $Q_t\in \mathbb{R}^{n\times n}$ and $R_t\in \mathbb{R}^{m\times m}$
	are positive semi-definite symmetric and positive definite symmetric matrix, respectively. Then, the optimal solution and minimum cost are given by
	\begin{eqnarray*}
	a_t&=&-L_tx_t-\theta_{t}^{-1}m_t,\\
	J_t&=&x_t^{tr}K_tx_t+2n_t^{tr}x_t+\sum_{s=t}^{T-1}\left(\xi_s^{tr}K_{s+1}\xi_s+2n^{tr}_{s+1}\xi_s-m_s^{tr}\theta_{s}^{-1}m_s \right),
	\end{eqnarray*}
	with
	\begin{align*}
	m_t &= B^{tr}_t (K_{t+1}\xi_t+ n_{t+1})+\beta_t, \\
	n_t &= (D_t -B_t L_t )^{tr}[n_{t+1} +K_{t+1} \xi_t]+\alpha_t-L_t^{tr}\beta_t, \quad  0 \le t \le T-1, \quad n_T=0.
	\end{align*}
	$\theta_t$, $L_t$ and $K_t$ are defined as
	\begin{eqnarray*}
		&&K_{t}=D^{tr}_{t}\left(K_{t+1}-K_{t+1}B_{t}\theta_{t}^{-1}B^{tr}_{t}K_{t+1}\right)D_{t}+Q_{t},\quad t=0, \cdots, T-1.\quad K_{T}=Q_{T}.\\
		&&L_{t}=\theta_{t}^{-1}B^{tr}_{t}K_{t+1}D_{t},\quad \theta_t=R_{t}+B^{tr}_{t}K_{t+1}B_{t}.
	\end{eqnarray*}
\end{lem}
\noindent \textit{Proof of Lemma \ref{LQ}.} This statement can be established as a straightforward application of the well known analytical expression of the solution to a quadratic program with equality
constraints; see \cite{nw}, Chapter 16. $\hfill\Box$

\smallskip

\noindent Now we proceed to demonstrate how to apply the DDP algorithm to the LQC problem in detail.
\begin{itemize}
\item[$-$] \textbf{Problem description}: Solve
\begin{eqnarray*}
&&\min_{\alpha\in \mathcal{A}_{\mathbb{F}}}\mathbb{E}\Big[\sum_{t=0}^{T-1}\left(x^{tr}_t Q_{t} x_t +\alpha^{tr}_{t} R_{t}\alpha_{t}\right) + x^{tr}_T Q_{T} x_T\Big],\\
s.t.  \quad && x_{t+1}=D_{t}x_t+B_{t}\alpha_t+\xi_t,\quad t=0, \cdots, T-1.
\end{eqnarray*}
\item[$-$] \textbf{Solution}: It is well known that the above control problem admits closed form solutions. For $t=0, \cdots, T-1$, the optimal policy should be $\alpha^*_{t}(x)=-L_{t}x$,
where the matrix $L_{t} \in \mathbb{R}^{m \times n}$ is given by
\[
L_{t}=(R_{t}+B^{tr}_{t}K_{t+1}B_{t})^{-1}B^{tr}_{t}K_{t+1}D_{t}.
\]
Here all matrices $K_{t}\in \mathbb{R}^{n \times n}$ are positive semidefinite symmetric, and we can use the following recursive relationship to determine them:
%\begin{eqnarray*}
%&&K_{T}=Q_{T};\\
%&&K_{t}=D^{tr}_{t}\left(K_{t+1}-K_{t+1}B_{t}(B^{tr}_{t}K_{t+1}B_{t}+R_{t})^{-1}B^{tr}_{t}K_{t+1}\right)D_{t}+Q_{t},\quad t=0, \cdots, T-1.
%\end{eqnarray*}
\begin{eqnarray*}
&&K_{T}=Q_{T};\\
&&K_{t}=D^{tr}_{t}\left(K_{t+1}-K_{t+1}B_{t}(R_{t}+B^{tr}_{t}K_{t+1}B_{t})^{-1}B^{tr}_{t}K_{t+1}\right)D_{t}+Q_{t},\quad t=0, \cdots, T-1.
\end{eqnarray*}
Under such a linear policy, the optimal cost function equals
\begin{eqnarray*}
V_{t}(x)=x^{tr}K_{t}x+\sum_{s=t}^{T-1}\mathbb{E}\left[\xi^{tr}_{s}K_{s+1}\xi_{s}\right].
\end{eqnarray*}
\end{itemize}

\noindent \textit{Proof of Proposition \ref{pro:lq}.} Consider a policy of the linear form
\begin{eqnarray}
\label{app_LQ_a0}
\alpha_{t}(x)=P_t x+E_t.
\end{eqnarray}
The subsequent calculation shows that we can achieve the optimal policy and value function of the LQC problem in two iterations by the DDP algorithm.
\begin{enumerate}
\item \textbf{First iteration}. Under linear policy (\ref{app_LQ_a0}), it is easy to verify that the cost-to-go function is quadratic with respect to states:
\begin{eqnarray*}
W^0_{t}(x)=x^{tr}H_{t}x+2F_t^{tr}x+C_t,
\end{eqnarray*}
with
\begin{eqnarray*}
H_t&=&Q_t+P_t^{tr}R_tP_t+(D_t+B_tP_t)^{tr}H_{t+1}(D_t+B_tP_t),\quad H_T=Q_T,\\
F_t&=&P_t^{tr}R_t^{tr}E_t+(D_t+B_tP_t)^{tr}H_{t+1}B_tE_t, \quad F_T=0,\\
C_t&=&C_{t+1}+E_t^{tr}R_tE_t+E_t^{tr}B_t^{tr}H_{t+1}B_tE_t+\mathbb{E}[\xi_t H_{t+1}\xi_t]+2F_{t+1}^{tr}B_tE_t,\quad C_T=0.
\end{eqnarray*}

Given $W^0_t$, we can construct the penalty function $z^{1}_{t}(a, \xi)$ by
\begin{eqnarray*}
z^{1}_{t}(a, \xi)&=&\sum_{s=t}^{T-1}\left\{\mathbb{E}[W^{0}_{s+1}(D_{s}x_s+B_{s}a_s+\xi_s)]-W^{0}_{s+1}(D_{s}x_s+B_{s}a_s+\xi_s)\right\}\nonumber\\
&=&\sum_{s=t}^{T-1}\left\{-2\xi_s^{tr}(F_s+H_{s+1}(D_{s}x_s+B_{s}a_s))-\xi_s^{tr}H_{s+1}\xi_s+\mathbb{E}[\xi_s^{tr}H_{s+1}\xi_s]\right\}.
\end{eqnarray*}
Then, the dual value in the first iteration satisfies
\begin{eqnarray}
\label{dual1_LQ}
\underline{V}^{1}_{t}(x)&=&\mathbb{E}\Big[ \inf_{a \in A|t}\Big\{\sum_{s=t}^{T-1}\Big(x_s^{tr}Q_s x_s+a_s^{tr}R_s a_s -2\xi_s^{tr}(F_s+H_{s+1}(D_{s}x_s+B_{s}a_s))\nonumber\\
&&-\xi_s^{tr}H_{s+1}\xi_s+\mathbb{E}[\xi_s^{tr}H_{s+1}\xi_s]\Big)+ x_T^{tr}Q_T x_T\Big\}\Big| x_{t}=x\Big].
\end{eqnarray}
Using Lemma \ref{LQ}, we can explicitly solve the inner optimization problem in (\ref{dual1_LQ}). It is a quadratic program. That leads to
\begin{eqnarray*}
\label{dual2_LQ}
\underline{V}^{1}_{t}(x) = V_t(x)-\mathbb{E}\left[\sum_{s=t}^{T-1}m_s^{tr}(R_{s}+B^{tr}_{s}K_{s+1}B_{s})^{-1}m_s\right], \quad t=0, \cdots, T-1,
\end{eqnarray*}
with
\begin{align*}
m_t &= B^{tr}_t (K_{t+1}-H_{t+1})\xi_t+ B^{tr}_t n_{t+1}, \\
n_t &= (D_t -B_t L_t )^{tr}[n_{t+1} +(K_{t+1}-H_{t+1}) \xi_t], \quad  0 \le t \le T-1, \quad n_T=0.
\end{align*}
\item \textbf{Second iteration}.
Note that $\underline{V}^{1}_{t}(x)$ is represented as the optimal value function $V_t(x)$ minus some constant. Hence, it is easy to see that $z^{2}_{t}(a, \xi)$ is the optimal
penalty function if we use $\underline{V}^{1}_{t}(x)$ to construct it. From this observation, we can calculate out that the dual value after the second iteration satisfy
\[\underline{V}^{2}_{t}(x) = V_t(x),\quad \alpha^{2}_{t}(x) = \alpha^*_t(x).\square \]
\end{enumerate}

\section{DC Optimization}
\label{app:dc}

In this appendix we briefly review some primary facts about DC functions and the related optimization problem. A function $f$ is called a DC function if there exist convex functions,
$g$ and $h$: $\mathbb{R}^{n}\rightarrow \mathbb{R}$ such that $f$ can be decomposed to the difference between $g$ and $h$:
\[
f(x)=g(x)-h(x),\quad \forall x \in \mathbb{R}^{n}.
\]
The set of DC functions has a very rich structure. For instance, Lemma \ref{lem1} points out that the class of DC functions is closed under some algebraic operations such as addition, multiplication,
and max/min.
 \begin{lem}[Theorem 4.1 in \cite{hpt}]
\label{lem1}
If $f_{1}$ and $f_{2}$ are two DC functions, then the following functions are also DC:\\
(a) $\lambda_{1}f_{1}(x)+\lambda_{2}f_{2}(x)$ for any constants $\lambda_{1}$ and $\lambda_{2}$,\\
(b) $\max\{f_{1}(x), f_{2}(x)\}$ and $\min\{f_{1}(x), f_{2}(x)\}$,\\
(c) $f_{1}(x)f_{2}(x)$.
\end{lem}

The standard form of a DC programming problem is given by
\begin{eqnarray}
\label{dc function}
\min \quad && f_0(x)-g_0(x)\\
s.t. \quad && f_i(x)-g_i(x)\le 0, \quad i=1,\cdots m,\nonumber\\
&& x\in \mathcal{X},\nonumber
\end{eqnarray}
where $\mathcal{X}\in \mathbb{R}^n$ is a nonempty closed convex set, and $f_i$'s, $g_i$'s are all convex in $\mathcal{X}$. Recently a sequential-convex-programming based DC algorithm and
its variations emerge as an effective approach to solving the problem. The idea of this approach is to create a sequence of values $\{x^{k}\}$ by solving convex programs sequentially so that $\{x^{k}\}$ converges to a local minimum
of (\ref{dc function}). Given a convex function $g$, a real vector $v$ is called its \emph{subgradient} at $x$ if $v$ satisfies
\[
g(y)\ge g(x)+v^{T}(y-x)\quad \text{for all } y,
\]
where $v^{T}$ is the transpose of vector $v$. Let $\partial g(x)$ be the set of all the subgradients of function $g$ at $x$. Using this notation, we can present the overarching structure of the
method in the following table:
\begin{framed}
	{\small \begin{center}
			\textbf{Table II: A Sequential Convex Programming Method}
		\end{center}
		\begin{itemize}[labelwidth=-0.25in,itemindent=-0.5in]
			\item
			\textbf{Step 0.} Choose $x^0\in \mathcal{X}$ arbitrarily. Set $k=0$.
			\item
			\textbf{Step 1.} Compute $s_{g_i}^k\in \partial g_i(x^k)$ for $i=0,1,\cdots,m.$
			\item
			\textbf{Step 2.}  Solve
			 \begin{eqnarray*}
			     x^{k+1}\in \arg \min_{y\in \mathcal{C}(x^k,\{s^k_{g_i}\}_{i=1}^{m})} \big\{f_0(y)-[g_0(x^k)+(s^k_{g_0})^T(y-x^k)] \big\}
			  \end{eqnarray*}
			\quad \quad with the feasible set $\mathcal{C}(x^k,\{s^k_{g_i}\}_{i=1}^{m})$ being given by
			\begin{eqnarray*}
			    \mathcal{C}(x^k,\{s^k_{g_i}\}_{i=1}^{m})=\big\{y\in \mathcal{X}: f_i(y)-[g_i(x^k)+(s^{k}_{g_i})^T(y-x^k)]\le 0, \ i=1,\cdots,m \big\}.
			\end{eqnarray*}
			\item
			\textbf{Step 3.}
			Set $k\leftarrow k+1$ and go to Step 1.
	\end{itemize}}
\end{framed}

Note that in Step 2, we linearize all the convex functions $g_{i}$, $i=1, \cdots, m$, through their subgradients, thereby relaxing the original problem into a tractable convex program.
A number of literature shows that the resulted sequence $\{x^k\}$ converge to a KKT point of (\ref{dc function}) under some regularity conditions; see, e.g., see \cite{yr}, \cite{ls2009}, \cite{l}, and \cite{bv}.

\section{Convergence of the Monte Carlo DDP Algorithm}
\label{app:convergence}

\smallskip

\subsection{Review of the algorithm}
\label{ec_cov_no}
Let us go through the major steps of the regression-based DDP algorithm proposed in Section \ref{sec:mc}. It is summarized in the following table.
\begin{framed}
{\small \begin{center}
\textbf{Table III: Implementation Details of Regression Based Monte Carlo DDP}
\end{center}
\begin{itemize}[labelwidth=-0.25in,itemindent=-0.5in]
\item
\textbf{Step 0.} Initialization:
\begin{itemize}[leftmargin=0in]
\item \textbf{Step 0a.} Choose a sequence of distribution functions $(G_{1}, \cdots, G_{T})$ and a set of basis functions
\indent \indent $\{\psi_{1}, \cdots, \psi_{M}\}$.
\item \textbf{Step 0b.} Simulate states for each period $t$ from these distributions: $x^{(l)}_{t} \sim G_{t}$ for $1 \le l \le L$ and
\indent \indent $1 \le t \le T$.
\item \textbf{Step 0c.} Construct the initial approximation
\[
\underline{\widehat{\mathfrak{V}}}^{0}_{t}(x):=\sum_{m=1}^{M}\widehat{\beta}^{0}_{t, m}\psi_{m}(x),\quad 1 \le t \le T-1.
\]
\indent \indent One way to do it is to evaluate the value of a policy $\alpha$ of being at state $x^{(l)}_{t}$ for all $l$ and $t$ and use
\indent \indent the basis functions to extrapolate these values to the entire state space. See Section 9.1 of \\
\indent \indent \cite{p} for the discussion on sampling and approximating the value of a policy.
\end{itemize}
\item
\textbf{Step 1.} Use the regression method to implement the dual iteration:
\begin{itemize}[leftmargin=0in]
\item \textbf{Step 1a.} Starting with the approximation from the last iteration:
\[
\underline{\widehat{\mathfrak{V}}}^{n-1}_{t}(x):=\sum_{m=1}^{M}\widehat{\beta}^{n-1}_{t, m}\psi_{m}(x),\quad 1 \le t \le T-1,
\]
\quad \quad define a penalty function sequence such that $\mathfrak{z}^{n}_{T}(a, \xi)=0$ and
\begin{eqnarray*}
\resizebox{.9\hsize}{!}{$
\mathfrak{z}^{n}_{t}(a, \xi)=\sum_{s=t}^{T-1}\left\{\mathbb{E}[r_{s}(x_{s}, a_{s}, \xi_{s})+\underline{\widehat{\mathfrak{V}}}^{n-1}_{s+1}(f_{s}(x_{s}, a_{s}, \xi_{s}))]
-(r_{s}(x_{s}, a_{s}, \xi_{s})+\underline{\widehat{\mathfrak{V}}}^{n-1}_{s+1}(f_{s}(x_{s}, a_{s}, \xi_{s})))\right\}
$}
\end{eqnarray*}
\indent \indent for any $0 \le t \le T-1$, with $a=(a_{0}, \cdots, a_{T-1}) \in A$ and $\xi=(\xi_{0}, \cdots, \xi_{T-1})$.
\item \textbf{Step 1b.} At each point $x^{(l)}_{t}$, simulate one path of $\xi^{(l)}|t=(\xi^{(l), t}_{t}, \xi^{(l), t}_{t+1}, \cdots, \xi^{(l), t}_{T-1})$ independently
\indent \indent and solve the optimization program (\ref{j}-\ref{constraint}) for $\mathfrak{J}^{(l)}_{t, n}$.
\item \textbf{Step 1c.} Use the least-square method to fit the data $(x^{(1)}_{t}, \mathfrak{J}^{(1)}_{t, n}), \cdots, (x^{(L)}_{t}, \mathfrak{J}^{(L)}_{t, n})$ to obtain a new
\indent \indent expansion on the dual:
\[
\underline{\widehat{\mathfrak{V}}}^{n}_{t}(x):=\sum_{m=1}^{M}\widehat{\beta}^{n}_{t, m}\psi_{m}(x),\quad 1 \le t \le T-1,
\]
\indent \indent where $\widehat{\beta}^{n}_{t}=(\widehat{B}^{t}_{\psi\psi})^{-1}\widehat{B}^{t, n}_{\mathfrak{J}\psi}$ whenever $\widehat{B}^{t,n}_{\psi\psi}$ is invertible.
Here the $(i,j)$-element of matrix  $\widehat{B}^{t,n}_{\psi\psi}$ and the \indent \indent $k$-th element of vector $\widehat{B}^{t, n}_{\mathfrak{J}\psi}$ are defined in
(\ref{bij}) and (\ref{bj}), respectively. See the discussion below for \indent \indent the case in which the numerical inversion $\widehat{B}^{t,n}_{\psi\psi}$ is not stable.
\item \textbf{Step 1d.} At $x_{0}$, simulate $L$ independent paths of $\xi^{(l)}|0=(\xi^{(l), 0}_{1}, \xi^{(l), 0}_{2}, \cdots, \xi^{(l), 0}_{T})$, $1 \le l \le L$, \indent \indent and solve the optimization program (\ref{j}-\ref{constraint}) with $t=0$ for $\mathfrak{J}^{(l)}_{0,n}$. Let
\[
\underline{\widehat{\mathfrak{V}}}^{n}_{0}(x_{0})=\frac{1}{L}\sum_{l=1}^{L}\mathfrak{J}^{(l)}_{0,n}.
\]
\end{itemize}
\item
\textbf{Step 2.}
Let $n=n+1$ and go to Step 1.
\end{itemize}}
\end{framed}

In the implementation of  Step 1c, we find that $\widehat{B}^{t,n}_{\psi\psi}$ could be nearly singular for some sampled $(x^{(1)}_{t,n}, \cdots, x^{(L)}_{t,n})$. That will
result in numerical instability on $\widehat{\beta}^{n}_{t}$, and in turn, the final dual output $\underline{\widehat{\mathfrak{V}}}^{T+1}_{0}$. To prevent
$\underline{\widehat{\mathfrak{V}}}^{T+1}_{0}$ from being extremely large or small due to the singularity of $\widehat{B}^{t,n}_{\psi\psi}$, we truncate the
output at a pre-specified sufficiently large $K$ in the numerical experiments, i.e.,
\[
 \underline{\widehat{\mathfrak{V}}}^{T+1}_{0}(x)=\max\left\{-K,\min\left\{K,\frac{1}{L}\sum_{l=1}^{L}\mathfrak{J}^{(l)}_{0, T+1}\right\}\right\}.
\]
Lemma \ref{lem_A} provides an upper bound on the probability that the matrix $\widehat{B}^{t,n}_{\psi\psi}$ is close to singularity. As both $L$ and $M$ tend to infinity, the
probability of near-singular $\widehat{B}^{t,n}_{\psi\psi}$ will vanish.

The output of our regression based algorithm $\{\underline{\widehat{\mathfrak{V}}}_{t}(x), 0 \le t \le T\}$ can also be used to simulate for an upper-bound estimate for the true
value of the original problem. The key steps are summarized in Table IV. Note that all the policies are suboptimal. It is obvious to see that
$\widehat{\overline{\mathfrak{V}}}_{0}$ will converge to one upper bound for the true value as $K \rightarrow +\infty$. Furthermore, we may construct a confidence interval
based on $\underline{\widehat{\mathfrak{V}}}_{0}$ and $\widehat{\overline{\mathfrak{V}}}_{0}$. Let $\underline{\sigma}$ and $\overline{\sigma}$ be the sample
standard deviations of $\{\mathfrak{J}^{(l)}_{0}, l=1, \cdots, L\}$ in Step 1d of Table III and $\{\sum_{t}r_{t}(x^k_{t}, \mathfrak{a}^k_{t}), k=1, \cdots, K\}$ in Step 3 of Table IV.
Then, we can form the following interval:
\begin{eqnarray}
\label{confidence}
\left(\underline{\widehat{\mathfrak{V}}}_{0}-z_{\delta/2}\frac{\underline{\sigma}}{\sqrt{L}},\ \widehat{\overline{\mathfrak{V}}}_{0}+z_{\delta/2}\frac{\overline{\sigma}}{\sqrt{K}}\right),
\end{eqnarray}
with $z_{\delta}$ being the $1-\delta$ quantile of the standard normal distribution. By Theorem \ref{thm:error} and Remarks \ref{rem1}, \ref{rem2}, this interval (\ref{confidence}) will provide a valid asymptotic confidence interval for $V_{0}$.
\begin{framed}
	{\small \begin{center}
			\textbf{Table IV: Direct Policy Valuation}
		\end{center}
		\begin{itemize}[labelwidth=-0.25in,itemindent=-0.5in]
			\item
			\textbf{Step 0.} Initialization: start from the initial state $x_0$ and choose a large number $K$.
			\item
			\textbf{Step 1.} Do for $k=1, \cdots, K$
			\begin{itemize}[leftmargin=0in]
			\item \textbf{Step 1a.} Set $t=0$ and let $x^k_{t}=x_{0}$.
			\item \textbf{Step 1b.} At $x^k_{t}$, solve the best action $\mathfrak{a}^k_{t}$, given the value function at the next step is approximated
			\indent \indent \indent by $\underline{\widehat{\mathfrak{V}}}_{t+1}$. That is,
			\begin{eqnarray*}
	        \mathfrak{a}^k_{t}=\argmin_{a_{t} \in A_{t}}\mathbb{E}\left[r_t(x^k_{t}, a_{t}, \xi_{t})+\underline{\widehat{\mathfrak{V}}}_{t+1}(f_t(x^k_{t},a_t,\xi_t))\right].
            \end{eqnarray*}
			\item \textbf{Step 1c.} Simulate $\xi^k_{t}$ and generate the state for the next step through $x^k_{t+1}=f_{t}(x^k_{t}, \mathfrak{a}^k_{t}, \xi^k_{t})$.
			\item \textbf{Step 1d.} Set $t\leftarrow t+1$ and go to Step 1b until $t=T$.
			\end{itemize}
			\item \textbf{Step 3.} Compute
			\[
			\widehat{\overline{\mathfrak{V}}}_{0}(x_{0}):=\frac{1}{K}\sum_{k=1}^{K}\sum_{t}r_{t}(x^k_{t}, \mathfrak{a}^k_{t}).
			\]
	\end{itemize}}
\end{framed}

\subsection{One example of exploration pitfall}
\label{app:ee}

As noted in Section \ref{sec:mc}, the state sampler $G$ is crucial to ensure the convergence of the DDP algorithm. This subsection presents one example to illustrate a possible
exploration pitfall if we use a policy-dependent sampler to draw the states on which we estimate the dual values.

Consider the following 2-period SDP problem:
\begin{eqnarray*}
	\min_u\ \mathbb{E}\Big[\sum_{t=0}^{2}-(x_t-10)^{+}|x_{0}=x\Big].
\end{eqnarray*}
Here, the control $u_t$ can only be taken from the set $\{0,1,2\}$ and the dynamic satisfies
\[
x_{t+1}=20+10u_t(u_t-2)-u_t\xi_t=\left\{
\begin{array}{lcl}
20,       &      & {u_t=0}\\
10- \xi_t,    &      & {u_t=1}\\
20-2\xi_t,    &      & {u_t=2}
\end{array}. \right.
\]
The random noise $\xi_t$ follows the uniform distribution $U(0,10)$. It is easy to see that the optimal value functions of the problem at $t=0, 1, 2$
are given by
 \[
 V_2(x)=-(x-10)^+,\quad V_1(x)=-10-(x-10)^+,\ \textrm{and}\ V_0(x)=-20-(x-10)^+,
 \]
 respectively. And the corresponding optimal policy is $u_{t}(x)=0$ for all $t$ and $x$.

 Suppose that the set of basis functions we take is
 \[
 \Psi(x)=[\psi_{1}(x), \psi_{2}(x), \psi_{3}(x)]:=[1,x,(10-x)^+].
 \]
And we are given by an initial policy $u_{t}=1$ for all $t=0, 1, 2$; that is, the policy always selects the action of 1 no matter what state and period the planner is at. Instead of using an independent state sampler as suggested
in Step0a of Table III, let us consider the situation that we rely on such $u$ to drive the system to obtain the states that we may estimate the dual values later. Denote them by $(x^{(1)}_{t}, \cdots, x^{(L)}_{t})$, $t=0, 1, 2$. Note
that all of them are in $(0, 10)$.

Evaluating the value of this policy on these states, we know that all the values are $\widehat{\mathfrak{V}}^{0}_{t}(x^{(l)}_{t})=0$. If we use the regression technique to extrapolate these
values to the entire state space, we need to solve
\begin{eqnarray}
\label{counter_example}
\inf_{\beta} \frac{1}{L}\sum_{l=1}^{L}\left(\Psi(x^{l}_{t})\beta-\underline{\widehat{\mathfrak{V}}}^{0}_{t}(x^{l}_{t})\right)^2
\end{eqnarray}
for the regression coefficients $\hat{\beta}^{0}_{t}$. It is easy to see that this is an underdetermined problem in the sense that infinitely many $\beta$ are the minimizer of the term on the right hand side of (\ref{counter_example}).

Take one solution $\hat{\beta}^{0}_{t, i}=0$ with $i=1, 2, 3$ for all $t$; that is, we extrapolate $\widehat{\mathfrak{V}}^{0}_{t}(x)=0$ to the entire space as the approximate value used in Step 1a of Table III. Substitute it into the expression
of the penalty. Following Step 1b in Table III, we solve the inner optimization problem (\ref{j}-\ref{constraint}) at $x^{(l)}_{t}\in (0,10)$ and obtain
\begin{eqnarray}
\label{counter_example2}
\mathfrak{J}^{(l)}_{0, 1}=-20,\quad \mathfrak{J}^{(l)}_{1, 1}=-10,\quad \mathfrak{J}^{(l)}_{2, 1}=0.
\end{eqnarray}
Note that all these values have nothing to do with the random noise $\xi$. Use $\mathfrak{J}^{(l)}_{1, 1}$ as an example to explain how the above is calculated. As $\underline{\widehat{\mathfrak{V}}}^{0}_{t}(x)=0$ for all $t$,
the penalty function $\mathfrak{z}^{1}_{t}(u, \xi)$ should also be zero. Then
\begin{eqnarray*}
\mathfrak{J}^{(l)}_{1, 1}(x^{(l)}_{1}, \xi^{(l),1}_{1})&=&\inf_{u_1}\{-(x^{(l)}_{1}-10)^+ -(x_2^{(l)}-10)^+\}\\
&=&\inf_{u_1}\{-(x^{(l)}_{1}-10)^+ -(20+10u_1(u_1-2)-u_1\xi^{(l),1}_{1}-10)^+\}.
\end{eqnarray*}
Apparently, $u_1=0$ is the optimal solution to this inner optimization problem. We thus have
\[
\mathfrak{J}^{(l)}_{1, 1}(x^{(l)}_{1}, \xi^{(l),1}_{1})=(x^{(l)}_{1}-10)^+ -10= -10
\]
because $x^{(l)}_{1} \in (0, 10)$.

Under (\ref{counter_example2}), after we fit these values using the basis functions according to the  Step 1c in Table III, we know that
\[
\hat{\beta}_0^1=[-20,0,0], \quad \hat{\beta}_1^1=[-10,0,0],\quad \hat{\beta}_2^1=[0,0,0].
\]
That is,
\[
\widehat{\mathfrak{V}}^{1}_{0}(x)=-20, \quad \widehat{\mathfrak{V}}^{1}_{1}(x)=-10,\quad \widehat{\mathfrak{V}}^{1}_2(x)=0
\]
for all $x$. Repeat the calculation for more rounds of dual operation and we find that the dual value will not change, i.e., $\widehat{\mathfrak{V}}^{n}_{t}(x)=\widehat{\mathfrak{V}}^{1}_{t}(x)$ for all $x$ and $t=0, 1, 2$. No convergence to the optimal
value function will occur.

The above example shows that using control policies to generate the representative states may lead our DDP algorithm to be stuck in a suboptimal solution. The cause is that all the sampled states we select at the beginning are in $(0, 10)$ and no one falls in
$(10, 20)$, the other part of the state space. Lacking the related information in $(10, 20)$, the extrapolation from the regression cannot produce correct estimation for the value in that interval.

\subsection{Proof of Theorem \ref{thm:error}}
\label{app:con}

Now we turn to prove Theorem \ref{thm:error}. Below we will use $C$ to represent a generic constant, which is independent of $M$ and $L$. Note that it may change
step by step. In the theorem statement, we also use the following concept of Lebesgue constant. Consider a sequence of basis functions $\{\psi_m(x), m \ge 1\} $.
Given a function $f$ such that $\|f\|_{\infty}\ne 0$ and $\|f\|_{\infty}< \infty$, we use the standard least square method to find a proper expansion of $\{\psi_m(x), m \ge 1\} $
to approximate $f$; that is, let
\[
\widehat{\beta}_f=\arg\min_{\alpha}\mathbb{E}^G[\|f(x)-\Psi_M^{tr}(x)\alpha\|^2]
\]
and then $f \approx \Psi_M^{tr}\widehat{\beta}_f$.
\begin{dfn}[Lebesgue constant]
\label{dfn:lebesgue}
Define
\begin{eqnarray}
\label{lm}
l_M=\sup\left\{\frac{\|\Psi^{tr}_M(x)\widehat{\beta}_f\|_{\infty}}{\|f\|_{\infty}}:\|f\|_{\infty}\ne 0, \|f\|_{\infty}< \infty\right\}.
\end{eqnarray}
\end{dfn}

\medskip

\subsubsection{Technical Lemmas}
We need to establish several lemmas first.
\begin{lem}
\label{lem_inf}
For any function $f(x)$ and $g(x)$,
\[
\inf_x(f(x)+g(x))\ge \inf_xf(x)+\inf_xg(x)
\]
and
\[
\inf_xf(x)- \inf_xg(x)\ge\inf_x(f(x)-g(x)).
\]
\end{lem}

\begin{lem}
\label{lem_ineq}
For any $x$, $y\in \mathbb{R}$, let constant
\[
K\ge |y|.
\]
Then we have
\[
\Big|\max\Big\{-K,\min\big\{K,x\big\}\Big\}-y\Big|\le \Big|x-y\Big|.
\]

\end{lem}
\noindent \textit{Proof of Lemma \ref{lem_ineq}.} It can be easily verified . $\hfill\Box$\\

From Assumption \ref{as5} and \ref{as1}, we can establish the non-multicollinearity of basis functions as shown in the following lemma.
\begin{lem}
	\label{g}
	Under Assumption \ref{as5} and \ref{as1}, the smallest eigenvalue of matrix $B^t_{\psi \psi}$ is bounded away from zero uniformly in $M$.
\end{lem}
\noindent \textit{Proof of Lemma \ref{g}.} Note that $B^t_{\psi \psi}=\mathbb{E}^G[\Psi(X_{t})\Psi^{tr}(X_{t})]$ is a nonnegative definite matrix. Thus, its smallest
eigenvalue satisfies
\begin{eqnarray}
\label{lem:ec2}
\lambda_{\min}(B^t_{\psi \psi})=\min_{||w||_2=1}w^{tr}\mathbb{E}^G[\Psi(X_{t})\Psi^{tr}(X_{t})]w.
\end{eqnarray}
Moreover, by Assumption \ref{as1}, there exists an $\epsilon>0$ such that $dG/dF(x)>\epsilon$ for $x \in \mathcal{X}$. We have
\begin{eqnarray}
\label{lem:ec1}
\mathbb{E}^G[\Psi(X_{t})\Psi^{tr}(X_{t})]=\int_{\mathbb{R}^n} \Psi(x)\Psi^{tr}(x) \frac{dG}{dF} (x) dF(x) \ge \epsilon \int_{\mathbb{R}^n} \Psi(x)\Psi^{tr}(x) dF(x).
\end{eqnarray}
The orthogonality of the basis functions in Assumption \ref{as5} implies that the right hand side of the above inequality is given by $\epsilon \cdot I$, where
$I$ is an identity matrix. For any vector $w \in \mathbb{R}^n$, the inequality (\ref{lem:ec1}) implies that
\begin{eqnarray*}
w^{tr}\mathbb{E}^G[\Psi(X_{t})\Psi^{tr}(X_{t})] w \ge \epsilon w^{tr}w.
\end{eqnarray*}
In conjunction with (\ref{lem:ec2}), we have
\begin{eqnarray*}
\lambda_{\min}(B^t_{\psi \psi}) \ge \epsilon. \hfill\square
\end{eqnarray*}

Consider one sequence of i.i.d. random vectors $X_1,\cdots, X_L \in \mathbb{R}^d$ and  another sequence of i.i.d. random variables
$Y_1,\cdots, Y_L \in \mathbb{R}$. Suppose that all of $(X_i)_{1\le i\le L}$ are square integrable and their second moments are bounded
above by a constant. Furthermore, $(Y_i)_{1\le i \le L}$ is assumed to be essentially bounded, i.e., there exists $\|Y\|_{\infty}$ such that
\[
\max_{1\le i\le L}|Y_i| \le \|Y\|_{\infty}.
\]
Then, we have
\begin{lem}
\label{lem_exp}
There exists a constant $C$, independent of $L$ and $d$, such that
\begin{eqnarray*}
\mathbb{E}\Big[\Big\|\frac{1}{L}\sum_{l=1}^{L}X_lY_l-\mathbb{E}[X_lY_l]\Big\|_2\Big]\le \frac{C\sqrt{d}}{\sqrt{L}}\Big\|Y\Big\|_{\infty}.
\end{eqnarray*}
\end{lem}
\noindent \textit{Proof of Lemma \ref{lem_exp}.}
Let $X_l^k$ denote the $k$-th element of vector $X_l$.
According to Jensen's inequality,
\begin{eqnarray*}
\mathbb{E}\Big[\Big\|\frac{1}{L}\sum_{l=1}^{L}X_lY_l-\mathbb{E}[X_lY_l]\Big\|_2\Big]
&=&\mathbb{E}\Big[\Big\{\sum_{k=1}^d\Big(\frac{1}{L}\sum_{l=1}^{L}X_l^kY_l-\mathbb{E}[X_l^kY_l]\Big)^2\Big\}^{1/2}\Big]\\
&\le &\Big\{\mathbb{E}\Big[\sum_{k=1}^d\Big(\frac{1}{L}\sum_{l=1}^{L}X_l^kY_l-\mathbb{E}[X_l^kY_l]\Big)^2\Big]\Big\}^{1/2}.
\end{eqnarray*}
Observe that
\begin{eqnarray*}
\mathbb{E}\Big[\Big(\frac{1}{L}\sum_{l=1}^{L}X_l^kY_l-\mathbb{E}[X_l^kY_l]\Big)^{2}\Big]
=\frac{1}{L^2}\sum_{l=1}^{L}\mathbb{E}\Big[\Big(X_l^kY_l-\mathbb{E}[X_l^kY_l]\Big)^{2}\Big].
\end{eqnarray*}
Each summand on the right hand side of above equality satisfies
\begin{eqnarray*}
\mathbb{E}\Big[\Big(X_l^kY_l-\mathbb{E}[X_l^kY_l]\Big)^{2}\Big]
\le\mathbb{E}\Big[\Big(X_l^kY_l\Big)^2\Big]
\le C\Big\|Y\Big\|_{\infty}^2,
\end{eqnarray*}
if we take
\[
C=\max_{1\le k\le d}E[(X_{l}^k)^2].
\]
Accordingly, we have
\begin{eqnarray*}
\mathbb{E}\Big[\Big\|\frac{1}{L}\sum_{l=1}^{L}X_lY_l-\mathbb{E}[X_lY_l]\Big\|_2\Big]\le \frac{C\sqrt{d}}{\sqrt{L}}\Big\|Y\Big\|_{\infty}.\hfill\Box
\end{eqnarray*}
\\
The next lemma gives bound on the probability that the sample matrix $\widehat{B}^{t,n}_{\psi \psi}$ deviates from its mean $B^t_{\psi \psi}$. More precisely,  given a $\delta>0$, for any time $t$ and iteration $n$, let
\[
A_t^n(\delta)=\{ \|I-(B^t_{\psi \psi})^{-1}\widehat{B}^{t,n}_{\psi \psi}\|_2\ge \delta\},
\]
where $I$ is the identity matrix. We have
\begin{lem}
	\label{lem_A}
	There exists a constant $C$, independent of $M$ and $L$, such that for any $\delta$,
	\[
	\mathbb{P}(A_t^n(\delta))\le 2M \exp\{-\frac{L\delta^2}{CM^2}\}.
	\]
\end{lem}
\noindent \textit{Proof of Lemma \ref{lem_A}.} It is Lemma 2.1 in \cite{cc}. This inequality is also known as the matrix Bernstein inequality in the literature;
see also \cite{tr}. $\hfill\Box$\\
\\
In the following lemma we develop an upper bound estimate on the distance between sample value $\mathfrak{J}_{t,n}(\xi | t, x)$ and the optimal value $V_t(x)$. To be more precisely,
\begin{lem}
	\label{lem_J}
	Given the initial state $x_t=x$ and the truncated randomness sequence $\xi | t=(\xi_{t}, \cdots, \xi_{T-1})$, the corresponding optimization problem $\mathfrak{J}_{t,n}(\xi|t,x)$, satisfies
	\begin{eqnarray*}
		\Big\|\mathfrak{J}_{t,n}(\xi | t, x)-V_t(x)\Big\|_{\infty} &\le& 2 \sum_{s=t+1}^T \Big\|\underline{\widehat{\mathfrak{V}}}_{s}^{n-1}(x)-V_s(x)\Big\|_{\infty}.
	\end{eqnarray*}
\end{lem}
\textit{Proof of Lemma \ref{lem_J}.}
Recall that $\mathfrak{J}_{t,n}(\xi | t, x)$ is defined by
\begin{eqnarray*}
\mathfrak{J}_{t,n}(\xi | t, x)= \inf_{a \in A|t}\Big(\sum_{s=t}^{T-1}r_{s}(x_s, a_{s}, \xi_{s}) + r_{T}(x_{T})+\mathfrak{z}^{n}_{t}(a, \xi)\Big)
\end{eqnarray*}
with
\begin{eqnarray*}
{\mathfrak{z}}^{n}_{t}(a, \xi)=\sum_{s=t}^{T-1}\Big\{\mathbb{E}\big[r_{s}(x_{s}, a_{s}, \xi_{s})+\underline{\widehat{\mathfrak{V}}}_{s+1}^{n-1}(x_{s+1})\big]
-\big(r_{s}(x_{s}, a_{s}, \xi_{s})
+\underline{\widehat{\mathfrak{V}}}_{s+1}^{n-1}(x_{s+1})\big)\Big\}.
\end{eqnarray*}
Following similar arguments as the proof of Lemma \ref{bellman_equation_dual}, we can show that $\mathfrak{J}_{t,n}$ admits the following recursive representation:
\begin{eqnarray*}
  \resizebox{.99\hsize}{!}{$
\mathfrak{J}_{t,n}(\xi | t, x)= \inf_{a \in A_t}\Big(\mathbb{E}[r_{t}(x_t, a, \xi_{t})+  \underline{\widehat{\mathfrak{V}}}_{t+1}^{n-1}(f_t(x,a,\xi_t))]- \underline{\widehat{\mathfrak{V}}}_{t+1}^{n-1}(f_t(x,a,\xi_t))+\mathfrak{J}_{t+1,n}(\xi | t+1, f(x,a,\xi_t))\Big).
$}
\end{eqnarray*}
By Lemma \ref{lem_inf}, we know that
\begin{eqnarray}
\label{ec_Jtn}
&&\mathfrak{J}_{t,n}(\xi | t, x)\nonumber\\
&\ge&
  \resizebox{.9\hsize}{!}{$
\inf_{a \in A_t}\mathbb{E}\Big[r_{t}(x_t, a, \xi_{t})+  \underline{\widehat{\mathfrak{V}}}_{t+1}^{n-1}(f_t(x,a,\xi_t))\Big]+ \inf_{a \in A|t}\Big\{\mathfrak{J}_{t+1,n}(\xi | t+1, f(x,a,\xi_t))-\underline{\widehat{\mathfrak{V}}}_{t+1}^{n-1}(f_t(x,a,\xi_t))\Big\},
$}
\nonumber \\
&=:&J_{1} + J_{2}.
\end{eqnarray}

Consider the part of $J_1$ on the right hand side of (\ref{ec_Jtn}). Note that $V_t(x)$ satisfies the Bellman equation,
\begin{eqnarray}
\label{pf_bellman}
V_t(x)=\inf_{a\in A_t}\mathbb{E}\Big[r_t(x,a,\xi_t)+V_{t+1}(f_t(x,a,\xi_t))\Big].
\end{eqnarray}
Then, the difference between $J_1$  and $V_t(x)$ should be
\begin{eqnarray*}
&&J_1-V_t(x)\\
&=&\inf_{a \in A_t}\mathbb{E}\Big[r_{t}(x_t, a, \xi_{t})+  \underline{\widehat{\mathfrak{V}}}_{t+1}^{n-1}(f_t(x,a,\xi_t))\Big]-\inf_{a\in A_t}\mathbb{E}\Big[r_t(x,a,\xi_t)+V_{t+1}(f_t(x,a,\xi_t))\Big]\\
&\ge& \inf_{a\in A_t}\mathbb{E}\Big[\underline{\widehat{\mathfrak{V}}}_{t+1}^{n-1}(f_t(x,a,\xi_t))-V_{t+1}(f_t(x,a,\xi_t))\Big],
\end{eqnarray*}
where the inequality in the last line is because of Lemma \ref{lem_inf}. Furthermore,  since $f_{t}(x, a, \xi_t) \in \mathcal{X}$ for any state $x$, action $a$ and random noise $\xi_t$, we have
\begin{eqnarray*}
\underline{\widehat{\mathfrak{V}}}_{t+1}^{n-1}(f_t(x,a,\xi_t))-V_{t+1}(f_t(x,a,\xi_t))
\ge -\sup_{x \in \mathcal {X}}\Big|\underline{\widehat{\mathfrak{V}}}_{t+1}^{n-1}(x)-V_{t+1}(x)\Big|
= -\Big\|\underline{\widehat{\mathfrak{V}}}_{t+1}^{n-1}(x)-V_{t+1}(x)\Big\|_{\infty}.
\end{eqnarray*}
Taking expectation with respect to $\xi_t$ and taking infimum over all possible actions $a\in A_t$ on both side of above inequality will lead to
\[
\inf_{a\in A_t}\mathbb{E}\Big[\underline{\widehat{\mathfrak{V}}}_{t+1}^{n-1}(f_t(x,a,\xi_t))-V_{t+1}(f_t(x,a,\xi_t))\Big]\ge -\Big\|\underline{\widehat{\mathfrak{V}}}_{t+1}^{n-1}(x)-V_{t+1}(x)\Big\|_{\infty}.
\]
That implies,
\begin{eqnarray}
\label{J_1}
J_1-V_t(x)\ge -\Big\|\underline{\widehat{\mathfrak{V}}}_{t+1}^{n-1}(x)-V_{t+1}(x)\Big\|_{\infty}.
\end{eqnarray}

Next we turn to $J_{2}$, the second part on the right hand of (\ref{ec_Jtn}). Substitute the definition of $\mathfrak{J}_{t+1,n}(\xi | t+1, f(x,a,\xi_t))$ into $J_{2}$. After some term
rearrangements, we obtain
\begin{eqnarray}
\label{ec:j2_1}
J_{2}&=&\inf_{a \in A|t}\Big\{\mathfrak{J}_{t+1,n}(\xi | t+1, f(x,a,\xi_t))-\underline{\widehat{\mathfrak{V}}}_{t+1}^{n-1}(f_t(x,a,\xi_t))\Big\}\nonumber\\
&=&\inf_{a \in A|t}\Big\{\sum_{s=t+1}^{T-1}\Big(\mathbb{E}\big[r_{s}(x_s, a_s, \xi_{s})+  \underline{\widehat{\mathfrak{V}}}_{s+1}^{n-1}(x_{s+1})\big]-\underline{\widehat{\mathfrak{V}}}_{s}^{n-1}(x_{s})\Big)+\Big(r_T(x_T)-\underline{\widehat{\mathfrak{V}}}_{T}^{n-1}(x_{T})\Big)\Big\}.\nonumber\\
\end{eqnarray}
Applying Lemma \ref{lem_inf} to the first term on the right hand side of the above equality,
\begin{eqnarray}
\label{ec:j2_2}
&&\inf_{a \in A|t}\Big\{\sum_{s=t+1}^{T-1}\Big(\mathbb{E}\big[r_{s}(x_s, a_s, \xi_{s})+  \underline{\widehat{\mathfrak{V}}}_{s+1}^{n-1}(x_{s+1})\big]-\underline{\widehat{\mathfrak{V}}}_{s}^{n-1}(x_{s})\Big)\Big\}\nonumber\\
&\ge&\sum_{s=t+1}^{T-1}\inf_{a \in A|t}\Big\{\inf_{a\in A_s}\mathbb{E}\Big[r_s(x_s,a,\xi_s)+\underline{\widehat{\mathfrak{V}}}_{s+1}^{n-1}(f_s(x_s,a,\xi_s))\Big]
-\underline{\widehat{\mathfrak{V}}}_{s}^{n-1}(x_{s})\Big\}\nonumber\\
&\ge& \sum_{s=t+1}^{T-1}\inf_{x\in \mathcal{X}}\Big\{\inf_{a\in A_s}\mathbb{E}\Big[r_s(x,a,\xi_s)+\underline{\widehat{\mathfrak{V}}}_{s+1}^{n-1}(f_s(x,a,\xi_s))\Big]-\underline{\widehat{\mathfrak{V}}}_{s}^{n-1}(x)\Big\}.
\end{eqnarray}
Here the last inequality is obvious because every summand of the sum in the second line, as a function of state variable $x$, is bounded below by its minimum over the space $\mathcal{X}$.
Similarly, we have
\begin{eqnarray}
\label{ec:j2_3}
r_T(x_T)-\underline{\widehat{\mathfrak{V}}}_{T}^{n-1}(x_{T})\ge\inf_{x\in \mathcal{X}}\Big\{r_T(x)-\underline{\widehat{\mathfrak{V}}}_{T}^{n-1}(x)\Big\}.
\end{eqnarray}
From (\ref{ec:j2_1}-\ref{ec:j2_3}),
\begin{eqnarray}
\label{ec:j2}
  \resizebox{.95\hsize}{!}{$
J_{2}\ge \sum_{s=t+1}^{T-1}\inf_{x\in \mathcal{X}}\Big\{\inf_{a\in A_s}\mathbb{E}\Big[r_s(x,a,\xi_s)+\underline{\widehat{\mathfrak{V}}}_{s+1}^{n-1}(f_s(x,a,\xi_s))\Big]
-\underline{\widehat{\mathfrak{V}}}_{s}^{n-1}(x)\Big\}+\inf_{x\in \mathcal{X}}\Big\{r_T(x)-\underline{\widehat{\mathfrak{V}}}_{T}^{n-1}(x)\Big\}.
$}
\end{eqnarray}

We add and subtract the optimal value function $V$ simultaneously in every summand of the sum on the right hand side of (\ref{ec:j2}). This operation will not change its value. That is,
\begin{eqnarray}
\label{ec:j2_4}
&&\inf_{x\in \mathcal{X}}\Big\{\inf_{a\in A_s}\mathbb{E}\Big[r_s(x,a,\xi_s)+\underline{\widehat{\mathfrak{V}}}_{s+1}^{n-1}(f_s(x,a,\xi_s))\Big]-\underline{\widehat{\mathfrak{V}}}_{s}^{n-1}(x)\Big\}\nonumber\\
&=&\inf_{x\in \mathcal{X}}\Big\{\inf_{a\in A_s}\mathbb{E}\Big[r_s(x,a,\xi_s)+V_{s+1}(f_s(x,a,\xi_s))+\Big(\underline{\widehat{\mathfrak{V}}}_{s+1}^{n-1}(f_s(x,a,\xi_s))-V_{s+1}(f_s(x,a,\xi_s))\Big)\Big]\nonumber\\
&&-V_{s}(x)+\Big(V_{s}(x)-\underline{\widehat{\mathfrak{V}}}_{s}^{n-1}(x)\Big)\Big\}\nonumber\\
&\ge&
\resizebox{.99\hsize}{!}{$
\inf_{x\in \mathcal{X}}\Big\{\inf_{a\in A_s}\mathbb{E}\Big[r_s(x,a,\xi_s)+V_{s+1}(f_s(x,a,\xi_s))\Big]-V_{s}(x)\Big\}+\inf_{x\in \mathcal{X}}\Big\{\inf_{a\in A_s}\Big(\underline{\widehat{\mathfrak{V}}}_{s+1}^{n-1}(f_s(x,a,\xi_s))-V_{s+1}(f_s(x,a,\xi_s))\Big)\Big\}
$}
\nonumber\\
&&+\inf_{x\in \mathcal{X}}\Big(V_{s}(x)-\underline{\widehat{\mathfrak{V}}}_{s}^{n-1}(x)\Big),
\end{eqnarray}
where we use Lemma \ref{lem_inf} again to obtain the last inequality. Thanks to the Bellman equation, we know that the first term on the right hand side of the inequality (\ref{ec:j2_4}) is $0$.
In addition, following similar arguments leading to (\ref{J_1}), we can establish
\[
\inf_{x\in \mathcal{X}}\Big\{\inf_{a\in A_s}\Big(\underline{\widehat{\mathfrak{V}}}_{s+1}^{n-1}(f_s(x,a,\xi_s))-V_{s+1}(f_s(x,a,\xi_s))\Big)\Big\} \ge -\Big\|\underline{\widehat{\mathfrak{V}}}_{s+1}^{n-1}(x)-V_{s+1}(x)\Big\|_{\infty}
\]
and
\[
\inf_{x\in \mathcal{X}}\Big(V_{s}(x)-\underline{\widehat{\mathfrak{V}}}_{s}^{n-1}(x)\Big) \ge -\Big\|\underline{\widehat{\mathfrak{V}}}_{s}^{n-1}(x)-V_{s}(x)\Big\|_{\infty}.
\]
As a consequence, we have
\begin{eqnarray*}
\resizebox{.99\hsize}{!}{$
\inf_{x\in \mathcal{X}}\Big\{\inf_{a\in A_s}\mathbb{E}\Big[r_s(x,a,\xi_s)+\underline{\widehat{\mathfrak{V}}}_{s+1}^{n-1}(f_s(x,a,\xi_s))\Big]-\underline{\widehat{\mathfrak{V}}}_{s}^{n-1}(x)\Big\}
\ge-\Big\|\underline{\widehat{\mathfrak{V}}}_{s}^{n-1}(x)-V_{s}(x)\Big\|_{\infty}-\Big\|\underline{\widehat{\mathfrak{V}}}_{s+1}^{n-1}(x)-V_{s+1}(x)\Big\|_{\infty}.
$}
\end{eqnarray*}

Summing the above inequality over $s=t+1$ to $T-1$, (\ref{ec:j2}) implies that
\begin{eqnarray}
\label{J_2}
J_2\ge -2\sum_{s=t+2}^{T}\|\underline{\widehat{\mathfrak{V}}}_{s}^{n-1}(x)-V_{s}(x)\|_{\infty}-\|\underline{\widehat{\mathfrak{V}}}_{t+1}^{n-1}(x)-V_{t+1}(x)\|_{\infty}.
\end{eqnarray}
Hence,
\begin{eqnarray*}
	\mathfrak{J}_{t,n}(\xi | t, x)-V_t(x) = J_{1}-V_{t}(x)+J_{2} \ge -2 \sum_{s=t+1}^T \|\underline{\widehat{\mathfrak{V}}}_{s}^{n-1}(x)-V_s(x)\|_{\infty}.
\end{eqnarray*}

To finish the proof, we need to derive the upper bound for $\mathfrak{J}_{t,n}(\xi | t, x)-V_t(x)$. By Assumption \ref{as3}, let $a^*(x)$ be the optimal solution to the Bellman equation;
that is,
\[
a^*_t=\arg\inf_{a\in A_t}\mathbb{E}\Big[r_t(x_t,a,\xi_t)+V_{t+1}(x_{t+1})\Big|x_t=x\Big].\]
Such a policy $a_t^*(x)$ must be a suboptimal solution to the optimization in the definition of $\mathfrak{J}_{t,n}(\xi | t, x)$. Therefore,
\begin{eqnarray}
\label{ec_J*}
	\mathfrak{J}_{t,n}(\xi | t, x)\le  \sum_{s=t}^{T-1}\Big\{\mathbb{E}\big[r_s(x_s^*,a_s^*(x_s^*),\xi_s)+\underline{\widehat{\mathfrak{V}}}_{s+1}^{n-1}(x_{s+1}^*)\big]-\underline{\widehat{\mathfrak{V}}}_{s+1}^{n-1}(x^*_{s+1})\Big\}+r_T(x^*_T),
\end{eqnarray}
with $x^*_{t+1}=f_t(x^*_t,a^*_t(x^*_t),\xi_t)$. On the other hand, we may rewrite $V_t(x^*_t)$ using the following telescoping sum:
\begin{eqnarray*}
V_t(x^*_t)&=& \sum_{s=t}^{T-1}\Big\{V_{s}(x^*_{s})-V_{s+1}(x^*_{s+1})\Big\}+r_T(x^*_T).
\end{eqnarray*}
Note that $V_{T}(\cdot) \equiv r_{T}(\cdot)$. By the Bellman equation, for all $s=t, \cdots, T_1$ and $x^*_{s}$,
\begin{eqnarray*}
	V_{s}(x^*_{s})=\mathbb{E}\big[r_s(x^*_s,a_s^*(x_s^*),\xi_s)+V_{s+1}(x_{s+1}^*)\big].
\end{eqnarray*}
Therefore,
\begin{eqnarray}
\label{ec_V*}
	V_t(x^*_t)&=& \sum_{s=t}^{T-1}\Big\{\mathbb{E}\big[r_s(x_s^*,a_s^*(x_s^*),\xi_s)+V_{s+1}(x_{s+1}^*)\big]-V_{s+1}(x^*_{s+1})\Big\}+r_T(x^*_T).
	\end{eqnarray}
Subtract the above two relation (\ref{ec_J*}) and (\ref{ec_V*}),
\begin{eqnarray*}
&&\mathfrak{J}_{t,n}(\xi | t, x_t^*)-V_t(x_t^*)\\
&\le& \sum_{s=t}^{T-1}\Big\{\mathbb{E}\big[\underline{\widehat{\mathfrak{V}}}_{s+1}^{n-1}(x_{s+1}^*)-V_{s+1}(x_{s+1}^*)\big]+V_{s+1}(x_{s+1}^*)-\underline{\widehat{\mathfrak{V}}}_{s+1}^{n-1}(x^*_{s+1})\Big\}.
\end{eqnarray*}
Following the similar procedures leading to (\ref{J_1}), for $t\le s\le T-1$, we can show that
 \[
 \mathbb{E}\big[\underline{\widehat{\mathfrak{V}}}_{s+1}^{n-1}(x_{s+1}^*)-V_{s+1}(x_{s+1}^*)\big]\le \big\| \underline{\widehat{\mathfrak{V}}}_{s+1}^{n-1}(x)-V_{s+1}(x) \big\|_{\infty}
 \]
and
\[
V_{s+1}(x_{s+1}^*)-\underline{\widehat{\mathfrak{V}}}_{s+1}^{n-1}(x^*_{s+1})\le \big\|V_{s+1}(x)-\underline{\widehat{\mathfrak{V}}}_{s+1}^{n-1}(x)\big\|_{\infty}.
\]
Consequently, we have
\begin{eqnarray*}
	\mathfrak{J}_{t,n}(\xi | t, x)-V_t(x) \le 2 \sum_{s=t+1}^T \Big\|\underline{\widehat{\mathfrak{V}}}_{s}^{n-1}(x)-V_s(x)\Big\|_{\infty}.
\end{eqnarray*}
In summary, the combination of these two bounds implies
\begin{eqnarray*}
	\Big\|\mathfrak{J}_{t,n}(\xi | t, x)-V_t(x)\Big\|_{\infty} \le 2 \sum_{s=t+1}^T \Big\|\underline{\widehat{\mathfrak{V}}}_{s}^{n-1}(x)-V_s(x)\Big\|_{\infty}.\hfill\square\\
\end{eqnarray*}

From the Lemma \ref{lem_J}, it turns out that
\begin{cor}
\label{cor_V}
Let
\[
\underline{\mathfrak{V}}_t^n(x)=\mathbb{E}[\mathfrak{J}_{t,n}(\xi | t, x)].
\]
Then we have
\begin{eqnarray}
\label{ec_jv1}
	\|\underline{\mathfrak{V}}_t^n(x)-V_t(x)\|_{\infty} &\le& 2 \sum_{s=t+1}^T \|\underline{\widehat{\mathfrak{V}}}_{s}^{n-1}(x)-V_s(x)\|_{\infty}.
\end{eqnarray}
\end{cor}
\textit{Proof of Corollary \ref{cor_V}.}
This can be easily verified by Jensen's inequality. $\hfill\square$\\
\\
In the next lemma, we attempt to bound the sampling error
\[
\|\Psi_M^{tr}(x)\widehat{\beta}_t^n - \Psi_M^{tr}(x)\beta_t^n\|_{\infty}
\]
when $\widehat{B}^{t,k}_{\psi \psi}$ gives a ``good" approximation to $B^t_{\psi \psi}$. To be precise, define event $A(\delta,n)$ to be
	\[
	% A(\delta,n)=\bigcup_{1\le k\le n,\atop T-k+1\le t\le T} A_t^k(\delta)=\bigcup_{1\le k\le n,\atop T-k+1\le t\le T}\{ \|I-(B^t_{\psi \psi})^{-1}\widehat{B}^{t,k}_{\psi \psi}\|_2\ge \delta\} .
	A(\delta,n)=\bigcup_{\substack{1\le k\le n,\\T-k+1\le t\le T}} A_t^k(\delta)=\bigcup_{\substack{1\le k\le n,\\T-k+1\le t\le T}}\{ \|I-(B^t_{\psi \psi})^{-1}\widehat{B}^{t,k}_{\psi \psi}\|_2\ge \delta\} .
	\]
\begin{lem}
	\label{lem_error}
	Let \[
	\delta=\frac{1}{2M^{1/2}L^{1/4}}
	\]
	in the above definition of $A(\delta, n)$.
	 There exists a constant $C$, independent of $M $ and $L$, such that for $1\le n\le T$ and $T-n+1\le t\le T$,
	\[\mathbb{E}\Big[1_{A(\delta,n)^c}\cdot\Big\|\Psi_M^{tr}(x)(\widehat{\beta}_t^n-\beta_t)\Big\|_{\infty}\Big]\le C\Big(\frac{M^{3/2}}{L^{1/4}}\Big)\mathbb{E}\Big[\Big\|\mathfrak{J}_{t,n}(\xi|t,x)\Big\|_{\infty}\cdot1_{A(\delta,n-1)^c}\Big].\]
	\end{lem}
\textit{Proof of Lemma \ref{lem_error}.}
%In the following we focus on the complement of set $A$, \[A^c=\bigcap_{1\le n\le T+1,\atop T-n+1\le t\le T}\{\|B^t_{\psi \psi}-\widehat{B}^{t,n}_{\psi \psi}\|_2< \delta\}.\]
Recall that if $A$ is a real symmetric matrix, then all the eigenvalues of this matrix is real. In this proof, we use $\lambda_{\max}(A)$ and $\lambda_{\min}(A)$ to denote the maximum and minimum eigenvalue for a general symmetric $A$. By the definitions of $\widehat{\beta}^n_t$ and $\beta^n_t$, we have
\[
\widehat{\beta}^n_t=(\widehat{B}_{\psi \psi}^{t,n})^{-1}\cdot \frac{1}{L}\sum_{l=1}^{L}\Psi_M(x_{t,n}^{(l)})\mathfrak{J}_{t,n}^{(l)}\quad \textrm{and}\quad
\beta^n_t=(B^t_{\psi \psi})^{-1}\mathbb{E}[\Psi_M(x_{t,n}^{(l)})\mathfrak{J}_{t,n}^{(l)}].
\]
Hence,
\[
\Psi_M^{tr}(x)(\widehat{\beta}^n_t-\beta^n_t)= \Psi_M^{tr}(x)(\widehat{B}_{\psi \psi}^{t,n})^{-1}\frac{1}{L}\sum_{l=1}^{L}\Psi_M(x_{t,n}^{(l)})\mathfrak{J}_{t,n}^{(l)}-\Psi_{M}^{tr}(x)(B^t_{\psi \psi})^{-1}\mathbb{E}[\Psi_M(x_{t,n}^{(l)})\mathfrak{J}_{t,n}^{(l)}].
\]
We simultaneously add and subtract
\[
\Psi_{M}^{tr}(x)(B^t_{\psi \psi})^{-1}\frac{1}{L}\sum_{l=1}^{L}\Psi_M(x_{t,n}^{(l)})\mathfrak{J}_{t,n}^{(l)}
\]
on the right hand side of the above equality. That results in
\begin{eqnarray*}
&&\Psi_M^{tr}(x)(\widehat{\beta}^n_t-\beta^n_t)\\
&=& \resizebox{.9\hsize}{!}{$
\Psi_M^{tr}(x)\Big[(\widehat{B}_{\psi \psi}^{t,n})^{-1}-(B^t_{\psi \psi})^{-1}\Big]\frac{1}{L}\sum_{l=1}^{L}\Psi_M(x_{t,n}^{(l)})\mathfrak{J}_{t,n}^{(l)}+\Psi_M^{tr}(x)(B_{\psi \psi}^{t})^{-1}\Big\{\frac{1}{L}\sum_{l=1}^{L}\Psi_M(x_{t,n}^{(l)})\mathfrak{J}_{t,n}^{(l)}-\mathbb{E}[\Psi_M(x_{t,n}^{(l)})\mathfrak{J}_{t,n}^{(l)}]\Big\}.
$}
\end{eqnarray*}
Then the triangle inequality implies
\begin{eqnarray*}
&&\|\Psi_M^{tr}(x)(\widehat{\beta}^n_t-\beta^n_t)\|_{\infty}\le\epsilon_{t,n}^{(1)}+\epsilon_{t,n}^{(2)},
\end{eqnarray*}
where $\epsilon_{t,n}^{(1)}$ and $\epsilon_{t,n}^{(2)} $ are defined as
\begin{eqnarray*}
\epsilon_{t,n}^{(1)}&=&\sup_{x\in \mathcal{X}}\Big|\Psi_M^{tr}(x)\Big[(\widehat{B}_{\psi \psi}^{t,n})^{-1}-(B^t_{\psi \psi})^{-1}\Big]\frac{1}{L}\sum_{l=1}^{L}\Psi_M(x_{t,n}^{(l)})\mathfrak{J}_{t,n}^{(l)}\Big|
\end{eqnarray*}
and
\begin{eqnarray*}
\epsilon_{t,n}^{(2)}&=&\sup_{x\in \mathcal{X}}\Big|\Psi_{M}^{tr}(x)(B^t_{\psi \psi})^{-1}\Big\{\frac{1}{L}\sum_{l=1}^{L}\Psi_M(x_{t,n}^{(l)})\mathfrak{J}_{t,n}^{(l)}-\mathbb{E}[\Psi_M(x_{t,n}^{(l)})\mathfrak{J}_{t,n}^{(l)}]\Big\}\Big|.
\end{eqnarray*}

By the Cauchy-Schwartz inequality, it is easy to see that $\epsilon_{t,n}^{(1)}$ is bounded by
\begin{eqnarray}
\label{ec_ep1}
	\epsilon_{t,n}^{(1)}&\le&\frac{1}{L}\sum_{l=1}^L\sup_{x\in \mathcal{X}}\Big|\Psi_M^{tr}(x)\Big[(\widehat{B}_{\psi \psi}^{t,n})^{-1}-(B^t_{\psi \psi})^{-1}\Big]\Psi_M(x_{t,n}^{(l)})\mathfrak{J}_{t,n}^{(l)}\Big|\nonumber\\
	&\le& \frac{1}{L}\sum_{l=1}^L\sup_{x\in \mathcal{X}}\Big\|\Psi_M^{tr}(x)\Big\|_2\cdot\Big\|(\widehat{B}_{\psi \psi}^{t,n})^{-1}-(B^t_{\psi \psi})^{-1}\Big\|_2\cdot\Big\|\Psi_M(x_{t,n}^{(l)})\mathfrak{J}_{t,n}^{(l)}\Big\|_2.
\end{eqnarray}
Under Assumption \ref{as2}, there exists a constant $C$ such that
\[
\sup_{x\in \mathcal{X}}\|\Psi_M^{tr}(x)\|_2\le CM.
\]
We next develop an upper bound for the last term on the right hand side of (\ref{ec_ep1}). Note that
\[
\Big\|\Psi_M(x_{t,n}^{(l)})\mathfrak{J}_{t,n}^{(l)}\Big\|_2
=  \Big[\sum_{m=1}^{M}\big(\psi_m(x_{t,n}^{(l)})\mathfrak{J}_{t,n}^{(l)}\big)^2\Big]^{\frac{1}{2}} \le \Big[\sum_{m=1}^{M}\big(\psi_m(x_{t,n}^{(l)}
\big)^2\Big]^{\frac{1}{2}} \cdot\Big\|\mathfrak{J}_{t,n}(\xi|t,x)\Big\|_{\infty}.
\]
By Assumption \ref{as2},
\[
\Big[\sum_{m=1}^{M}\big(\psi_m(x_{t,n}^{(l)})\big)^2\Big]^{\frac{1}{2}}\le CM.
\]
Hence,
\begin{eqnarray}
\label{proof_l}
\Big\|\Psi_M(x_{t,n}^{(l)})\mathfrak{J}_{t,n}^{(l)}\Big\|_2
\le CM\Big\|\mathfrak{J}_{t,n}(\xi|t,x)\Big\|_{\infty}.
\end{eqnarray}

To bound
\[
\|(\widehat{B}_{\psi \psi}^{t,n})^{-1}-(B^t_{\psi \psi})^{-1}\|_2,
\]
by Cauchy-Schwarz inequality,
\begin{eqnarray*}
\Big\|\Big(\widehat{B}_{\psi \psi}^{t,n}\Big)^{-1}-\Big(B^t_{\psi \psi}\Big)^{-1}\Big\|_2
&=&\Big\|\Big(I-(B_{\psi\psi}^t)^{-1}\widehat{B}_{\psi \psi}^{t,n}\Big)\Big((B_{\psi\psi}^t)^{-1}\widehat{B}_{\psi \psi}^{t,n}\Big)^{-1}
\Big(B_{\psi\psi}^t\Big)^{-1}\Big\|_2\\
&\le& \Big \|I-(B_{\psi\psi}^t)^{-1}\widehat{B}_{\psi \psi}^{t,n}\Big\|_2\cdot\Big\|\Big((B_{\psi\psi}^t)^{-1}\widehat{B}_{\psi \psi}^{t,n}\Big)^{-1}\Big\|_2\cdot\Big\|\Big(B_{\psi\psi}^t\Big)^{-1}\Big\|_2.
\end{eqnarray*}
From the definition of $A(\delta, n)$, we know that
\begin{eqnarray}
\label{ec_lam5}
\Big \|I-(B_{\psi\psi}^t)^{-1}\widehat{B}_{\psi \psi}^{t,n}\Big\|_2\le (2M^{1/2}L^{1/4})^{-1}
\end{eqnarray}
on the set of $A(\delta, n)^c$. Using Example 5.6.6 in \cite{hj},
\[
 \Big\|\Big(B_{\psi\psi}^t\Big)^{-1}\Big\|_2=\lambda_{\max}\Big(\big(B_{\psi\psi}^t\big)^{-1}\Big)=\frac{1}{\lambda_{\min}(B^t_{\psi\psi})}.
\]

As for
\[
\Big\|\Big((B_{\psi\psi}^t)^{-1}\widehat{B}_{\psi \psi}^{t,n}\Big)^{-1}\Big\|_2,
\]
it is well known that
\begin{eqnarray}
\label{ec_lam1}
\lambda_{\min}((B_{\psi\psi}^t)^{-1}\widehat{B}_{\psi \psi}^{t,n})=\min_{\|w\|=1}w^{tr}(B_{\psi\psi}^t)^{-1}\widehat{B}_{\psi \psi}^{t,n}w;
\end{eqnarray}
see Theorem 4.2.2 in \cite{hj}.
For any vector $w$ with $\|w\|=1$,
\[
w^{tr}(B_{\psi\psi}^t)^{-1}\widehat{B}_{\psi \psi}^{t,n}w=w^{tr}I w+w^{tr}((B_{\psi\psi}^t)^{-1}\widehat{B}_{\psi \psi}^{t,n}-I)w=1-w^{tr}(I-(B_{\psi\psi}^t)^{-1}\widehat{B}_{\psi \psi}^{t,n})w,
\]
where $I$ is an identity matrix. Hence,
\begin{eqnarray}
\label{ec_lam2}
\min_{\|w\|=1}w^{tr}(B_{\psi\psi}^t)^{-1}\widehat{B}_{\psi \psi}^{t,n}w = 1-\max_{\|w\|=1}w^{tr}(I-(B_{\psi\psi}^t)^{-1}\widehat{B}_{\psi \psi}^{t,n})w.
\end{eqnarray}
On the other hand, it is easy to show that
\begin{eqnarray}
\label{ec_lam4}
\max_{\|w\|=1}w^{tr}(I-(B_{\psi\psi}^t)^{-1}\widehat{B}_{\psi \psi}^{t,n})w=\|I-(B_{\psi\psi}^t)^{-1}\widehat{B}_{\psi \psi}^{t,n})\|_2.
\end{eqnarray}
Combining (\ref{ec_lam5}-\ref{ec_lam4}) yields
\begin{eqnarray}
\label{ec_lambda_min}
1_{A(\delta,n)^c}\cdot \lambda_{\min}((B_{\psi\psi}^t)^{-1}\widehat{B}_{\psi \psi}^{t,n})\ge 1_{A(\delta,n)^c}\cdot(1-(2M^{1/2}L^{1/4})^{-1}) \ge
\frac{1}{2}\cdot1_{A(\delta,n)^c}.
\end{eqnarray}
where the last inequality is due to the fact that $M, L \ge 1$. Thus,
\begin{eqnarray}
\label{proof_bd}
1_{A(\delta,n)^c}\cdot \Big\|\Big(\widehat{B}_{\psi \psi}^{t,n}\Big)^{-1}-\Big(B^t_{\psi \psi}\Big)^{-1}\Big\|_2
\le 1_{A(\delta,n)^c}\cdot \frac{1}{\lambda_{\min}(B_{\psi\psi}^t)M^{1/2}L^{1/4}}.
\end{eqnarray}

By (\ref{ec_ep1}), (\ref{proof_l}), and (\ref{proof_bd}), we have
\begin{eqnarray*}
	\mathbb{E}\Big[\epsilon_{t,n}^{(1)}\cdot1_{A(\delta,n)^c}\Big]
	\le \frac{CM^{3/2}}{\lambda_{\min}(B^t_{\psi \psi})L^{1/4}}\mathbb{E}\Big[\Big\|\mathfrak{J}_{t,n}(\xi|t,x)\Big\|_{\infty}\cdot1_{A(\delta,n)^c}\Big].
\end{eqnarray*}
From Lemma \ref{g}, $\lambda_{\min}(B^t_{\psi \psi})$ is  bounded below by some constant. Therefore, the right hand side of the above inequality is
further bounded by
\begin{eqnarray*}
	C\frac{M^{3/2}}{L^{1/4}}\mathbb{E}\Big[\Big\|\mathfrak{J}_{t,n}(\xi|t,x)\Big\|_{\infty}\cdot1_{A(\delta,n-1)^c}\Big],
\end{eqnarray*}
if we change the constant properly.

Using Cauchy-Schwartz inequality again, $\epsilon_{t,n}^{(2)}$ satisfies
\begin{eqnarray*}
\epsilon_{t,n}^{(2)}&\le&\sup_{x\in \mathcal{X}}\Big\|\Psi_{M}^{tr}(x)\Big\|_{2}\cdot\Big\|(B^t_{\psi \psi})^{-1}\Big\|_{2}\cdot\Big\|\frac{1}{L}\sum_{l=1}^{L}\Psi_M(x_{t,n}^{(l)})\mathfrak{J}_{t,n}^{(l)}-\mathbb{E}[\Psi_M(x_{t,n}^{(l)})\mathfrak{J}_{t,n}^{(l)}]\Big\|_2.
\end{eqnarray*}
Let $D_{t,n} $ denote the last term of right hand side in above inequality.
Note that
\[
\sup_{x\in \mathcal{X}}\|\Psi_{M}^{tr}(x)\|_{2}\le CM
\]
by Assumption \ref{as2} and
\[
\|(B^t_{\psi \psi})^{-1}\|_{2}=\lambda_{\min}^{-1}(B^t_{\psi \psi}).
\]
We have
\begin{eqnarray*}
\epsilon_{t,n}^{(2)}\le \frac{CM}{\lambda_{\min}(B^t_{\psi \psi})}\cdot D_{t,n}.
\end{eqnarray*}

The definition of $A(\delta,n)$ implies that
\[
A(\delta,n)^c \subseteq A(\delta,n-1)^c.
\]
Therefore,
\begin{eqnarray}
\label{ec_epe}
\mathbb{E}\Big[\epsilon_{t,n}^{(2)}\cdot1_{A(\delta,n)^c}\Big]\le\mathbb{E}\Big[\epsilon_{t,n}^{(2)}\cdot1_{A(\delta,n-1)^c}\Big]
\le \frac{CM}{\lambda_{\min}(B^t_{\psi \psi})}\cdot \mathbb{E}\Big[D_{t,n}\cdot1_{A(\delta,n-1)^c}\Big].
\end{eqnarray}
Let $\mathcal{G}_n$ be a $\sigma$-algebra defined as follows,
\[
\mathcal{G}_n=\sigma\left(\left\{(x^{(1)}_{t,k}, \cdots, x^{(L)}_{t,k}),(\xi^{(l)}|t)\right\},1\le k\le n\right).
\]
By the iterated law of conditional expectation, the expectation term on the right hand side of (\ref{ec_epe}) equals
\begin{eqnarray}
\label{ec_D1}
\mathbb{E}\Big[D_{t,n}\cdot1_{A(\delta,n-1)^c}\Big]= \mathbb{E}\left[ \mathbb{E} \left[D_{t,n}\cdot1_{A(\delta,n-1)^c}\Big|\mathcal{G}_{n-1}\right]\right].
\end{eqnarray}
Since the event $A(\delta,n-1)^c$ is measurable with respect to $\mathcal{G}_{n-1}$, we have
\begin{eqnarray}
\mathbb{E}\Big[D_{t,n}\cdot1_{A(\delta,n-1)^c}\Big]= \mathbb{E}\left[ \mathbb{E} \left[D_{t,n}\Big|\mathcal{G}_{n-1}\right]\cdot1_{A(\delta,n-1)^c}\right].
\end{eqnarray}
Following the proof of Lemma \ref{lem_exp},  we can show
\begin{eqnarray}
\label{ec_D2}
\mathbb{E} \Big[D_{t,n}\Big|\mathcal{G}_{n-1}\Big]\le \frac{C\sqrt{M}}{\sqrt{L}}\Big\|\mathfrak{J}_{t,n}(\xi|t,x)\Big\|_{\infty}.
\end{eqnarray}
In light of (\ref{ec_epe}-\ref{ec_D2}),
\begin{eqnarray*}
\mathbb{E}\Big[\epsilon_{t,n}^{(2)}\cdot1_{A(\delta,n)^c}\Big]\le C\frac{M^{3/2}}{\lambda_{\min}(B_{\psi\psi}^t)\sqrt{L}}\mathbb{E}\Big[\Big\|\mathfrak{J}_{t,n}(\xi|t,x)\Big\|_{\infty}\cdot 1_{A(\delta,n-1)^c}\Big].
\end{eqnarray*}

Using again the fact that $\lambda_{\min}(B_{\psi\psi}^t)$ is bounded below, the right hand side of above can be bounded by
\[
C\frac{M^{3/2}}{\sqrt{L}}\mathbb{E}\Big[\Big\|\mathfrak{J}_{t,n}(\xi|t,x)\Big\|_{\infty}\cdot 1_{A(\delta,n-1)^c}\Big],
\]
by changing constant $C$ properly. Finally we put the upper bounds of $\epsilon_{t,n}^{(1)}$ and $\epsilon_{t,n}^{(2)}$ together to conclude
\[
\mathbb{E}\Big[1_{A(\delta,n)^c}\cdot\Big\|\Psi_M^{tr}(x)(\widehat{\beta}_t^n-\beta_t)\Big\|_{\infty}\Big]\le C\Big(\frac{M^{3/2}}{L^{1/4}}\Big)\mathbb{E}\Big[\Big\|\mathfrak{J}_{t,n}(\xi|t,x)\Big\|_{\infty}\cdot1_{A(\delta,n-1)^c}\Big].\hfill\square
\]

\subsubsection{Proof of Theorem \ref{thm:error}} Let
\[
\delta=\frac{1}{2M^{1/2}L^{1/4}}
\]
as in Lemma \ref{lem_error} and define $A(\delta, T)$ accordingly. We decompose
\[
\mathbb{E}\Big[\Big|\underline{\widehat{\mathfrak{V}}}^{T+1}_0(x)-V_0(x)\Big|\Big]=\mathbb{E}\Big[1_{A(\delta,T)}\cdot\Big|\underline{\widehat{\mathfrak{V}}}^{T+1}_0(x)-V_0(x)\Big|\Big]+\mathbb{E}\Big[1_{A(\delta,T)^c}\cdot\Big|\underline{\widehat{\mathfrak{V}}}^{T+1}_0(x)-V_0(x)\Big|\Big].
\]
\emph{Step 1}. We plan to develop a bound for
\[
\mathbb{E}\Big[1_{A(\delta,T)}\cdot\Big|\underline{\widehat{\mathfrak{V}}}^{T+1}_0(x)-V_0(x)\Big|\Big].
\]
To this end, according to Lemma \ref{lem_A},
\[
\mathbb{P}(A(\delta,T))\le
% \sum_{1\le n\le T,\atop T-n+1\le t\le T}
\sum_{\substack{1\le n\le T,\\ T-n+1\le t\le T}}
\mathbb{P}(\|I-(B^t_{\psi \psi})^{-1}\widehat{B}^{t,n}_{\psi \psi}\|_2\ge \delta)\le T(T+1)M
\exp\left(-\frac{L^{1/2}}{CM^3}\right).
\]
Since $V_0(x)$ is bounded by Assumption \ref{as3} and $\underline{\widehat{\mathfrak{V}}}^{T+1}_{0}(x)$ is also truncated by pre-specified constant $K$ as stated in Section \ref{ec_cov_no}, there should exist a constant $C$ such that
\begin{eqnarray}
\label{ec_1}
\mathbb{E}[1_{A(\delta,T)}\cdot|\underline{\widehat{\mathfrak{V}}}^{T+1}_{0}(x)-V_0(x)|]\le C\mathbb{P}(A(\delta,T))\le CT(T+1)M \exp\left(-\frac{L^{1/2}}{CM^3}\right).
\end{eqnarray}

\medskip

\noindent \emph{Step 2}. We intend to establish the relationship between $\|\underline{\widehat{\mathfrak{V}}}_{t}^{n}(x)-V_{t}(x)\|_{\infty}$ and $\|\underline{\widehat{\mathfrak{V}}}_{t}^{n-1}(x)-V_{t}(x)\|_{\infty}$ for $1\le n\le T$ and $T-n+1\le t\le T$. Our claim is that
\begin{eqnarray}
\label{proof_induction}
&&\mathbb{E}\big[1_{A(\delta,n)^c}\|\underline{\widehat{\mathfrak{V}}}_{t}^{n}(x)-V_{t}(x)\|_{\infty}\big]\nonumber\\
&\le&(1+l_M)\Delta+C\frac{M^{3/2}}{L^{1/4}}+(2l_M+C\frac{M^{3/2}}{L^{1/4}})\sum_{s=t+1}^T\mathbb{E}[1_{A(\delta,n-1)^c}\|\underline{\widehat{\mathfrak{V}}}_{s}^{n-1}(x)-V_s(x)\|_{\infty}].
\end{eqnarray}

To show this, by adding and subtracting the term $\Psi_M^{tr}(x)\beta_t^n$ at the same time within $\|\underline{\widehat{\mathfrak{V}}}_{t}^{n}(x)-V_{t}(x)\|_{\infty}$, we have
\begin{eqnarray*}
\Big\|\underline{\widehat{\mathfrak{V}}}^n_{t}(x)-V_{t}(x)\Big\|_{\infty}&=&\Big\|\Psi_M^{tr}(x)\widehat{\beta}_t^n-V_{t}(x)\Big\|_{\infty}\\
&=&\Big\|\Big(\Psi_M^{tr}(x)\widehat{\beta}_t^n -\Psi_M^{tr}(x)\beta_t^n\Big)+\Big(\Psi_M^{tr}(x)\beta_t^n  -V_t(x)\Big)\Big\|_{\infty}.
\end{eqnarray*}
Then, according to the triangle inequality,
\begin{eqnarray}
\label{error}
\Big\|\underline{\widehat{\mathfrak{V}}}^n_{t}(x)-V_{t}(x)\Big\|_{\infty}
\le\Big\|\Psi_M^{tr}(x)\widehat{\beta}_t^n - \Psi_M^{tr}(x)\beta_t^n\Big\|_{\infty}+\Big\| \Psi_M^{tr}(x)\beta_t^n - V_t(x)\Big\|_{\infty}.
\end{eqnarray}
Note that Lemma \ref{lem_error} provides the upper bound on $\mathbb{E}[1_{A(\delta,n)^c}\|\Psi_M^{tr}(x)\widehat{\beta}_t^n - \Psi_M^{tr}(x)\beta_t^n\|_{\infty}]$. Henceforth we only need to consider how to bound the second part in right hand side of (\ref{error}).

Let
\[
\beta_t=\arg\min_{\alpha}\mathbb{E}^G[(V_t(x)-\Psi^{tr}_M(x)\alpha)^2].
\]
We add and subtract the term $\Psi_M^{tr}(x)\beta_t$ simultaneously in the next equation and use the triangle inequality again,
\begin{eqnarray}
\label{pf_tr1}
\Big\| \Psi_M^{tr}(x)\beta_t^n - V_t(x)\Big\|_{\infty}
&=&\Big \| \Big(\Psi_M^{tr}(x)\beta_t^n - \Psi_M^{tr}(x)\beta_t\Big)+\Big( \Psi_M^{tr}(x)\beta_t -V_t(x)\Big)\Big\|_{\infty}\nonumber\\
&\le & \Big\| \Psi_M^{tr}(x)\beta_t^n - \Psi_M^{tr}(x)\beta_t\Big\|_{\infty}+\Big\|\Psi_M^{tr}(x)\beta_t-V_t(x)\Big\|_{\infty}.
\end{eqnarray}
Under the basis function set $\Psi_M(x)$, $\Psi_M^{tr}(x)(\beta_t^n-\beta_t)$ is the least square estimation of function $\underline{\mathfrak{V}}_t^n(x)-V_t(x)$. Recall the definition of $l_M$ in (\ref{lm}). Then
\begin{eqnarray*}
 \| \Psi_M^{tr}(x)\beta_t^n - \Psi_M^{tr}(x)\beta_t\|_{\infty}\le l_M\|\underline{\mathfrak{V}}_t^n(x)-V_t(x)\|_{\infty}.
\end{eqnarray*}
To bound $\|\Psi_M^{tr}(x)\beta_t-V_t(x)\|_{\infty}$, Lemma 2.4 in \cite{cc} shows that
\[
\|\Psi_M^{tr}(x)\beta_t-V_t(x)\|_{\infty}\le (l_M+1)\Delta,
\]
with $\Delta$ representing the approximation error as defined in the Theorem statement.

Therefore $\| \Psi_M^{tr}(x)\beta_t^n - V_t(x)\|_{\infty}$ satisfies
\begin{eqnarray}
\label{pf_tr}
\| \Psi_M^{tr}(x)\beta_t^n - V_t(x)\|_{\infty}
\le l_M\|\underline{\mathfrak{V}}_t^n(x)-V_t(x)\|_{\infty}+(l_M+1)\Delta.
\end{eqnarray}

From Lemma \ref{lem_error}, relationship (\ref{error}) and (\ref{pf_tr}), we have
\begin{eqnarray}
\label{ec_mid}
	&&\mathbb{E}\big[1_{A(\delta,n)^c}\| \underline{\widehat{\mathfrak{V}}}_{t}^{n}(x) - V_t(x)\|_{\infty}\big]\\
	\le&& C\frac{M^{3/2}}{L^{1/4}}\mathbb{E}\big[\|\mathfrak{J}_{t,n}(\xi|t,x)\|_{\infty}1_{A(\delta,n-1)^c}\big]+l_M\mathbb{E}\big[1_{A(\delta,n)^c}\|\underline{\mathfrak{V}}_t^n(x)-V_t(x)\|_{\infty}\big]+(l_M+1)\Delta.\nonumber
\end{eqnarray}
We need to bound each term in the last line of above inequality.
%For the error $\|\underline{\mathfrak{V}}_t^n(x)-V_t(x)\|_{\infty}$ in (\ref{ec_mid}), by the definition of $\underline{\mathfrak{V}}_t^n(x)$ and Jensen's inequality, we know
%\begin{eqnarray*}
%	\Big\|\underline{\mathfrak{V}}_t^n(x)-V_t(x)\Big\|_{\infty}
%	=\Big\|\mathbb{E}\Big[\mathfrak{J}_{t,n}(\xi | t, x)-V_t(x)\Big]\Big\|_{\infty}
%	\le\mathbb{E}\Big[\Big\|\mathfrak{J}_{t,n}(\xi | t, x)-V_t(x)\Big\|_{\infty}\Big].
%\end{eqnarray*}
According to the Corollary \ref{cor_V}, we have
\begin{eqnarray}
\label{ec_jv1_2}
	\mathbb{E}\Big[\Big\|\underline{\mathfrak{V}}_t^n(x)-V_t(x)\Big\|_{\infty}\cdot1_{A(\delta,n)^c}\Big] &\le& 2 \sum_{s=t+1}^T\mathbb{E}\Big[\Big \|\underline{\widehat{\mathfrak{V}}}_{s}^{n-1}(x)-V_s(x)\Big\|_{\infty}\cdot1_{A(\delta,n-1)^c}\Big].
\end{eqnarray}
For $\|\mathfrak{J}_{t,n}(\xi|t,x)\|_{\infty}$, it satisfies
\[
\Big\|\mathfrak{J}_{t,n}(\xi|t,x)\Big\|_{\infty}\le \Big\|\mathfrak{J}_{t,n}(\xi|t,x)-V_t(x)\Big\|_{\infty}+\Big\|V_t(x)\Big\|_{\infty}.
\]
As the optimal value function $V_t(x)$ is bounded on compact set $\mathcal{X}$ in Assumption \ref{as3}, there exists a constant $C$ such that
\begin{eqnarray}
\label{ec_jv2}
\mathbb{E}\Big[\Big\|\mathfrak{J}_{t,n}(\xi|t,x)\Big\|_{\infty}\cdot1_{A(\delta,n-1)^c}\Big]&\le& C+\mathbb{E}\Big[\Big\|\mathfrak{J}_{t,n}(\xi|t,x)-V_t(x)\Big\|_{\infty}\cdot1_{A(\delta,n-1)^c}\Big]\nonumber\\
&\le&
C+2 \sum_{s=t+1}^T\mathbb{E}\Big[ \Big\|\underline{\widehat{\mathfrak{V}}}_{s}^{n-1}(x)-V_s(x)\Big\|_{\infty}\cdot1_{A(\delta,n-1)^c}\Big].
\end{eqnarray}

We combine (\ref{ec_mid}-\ref{ec_jv2}),
\begin{eqnarray*}
&&\mathbb{E}\big[1_{A(\delta,n)^c}\|\underline{\widehat{\mathfrak{V}}}_{t}^{n}(x)-V_{t}(x)\|_{\infty}\big]\\
&\le&(1+l_M)\Delta+C\frac{M^{3/2}}{L^{1/4}}+(2l_M+C\frac{M^{3/2}}{L^{1/4}})\sum_{s=t+1}^T\mathbb{E}[1_{A(\delta,n-1)^c}\|\underline{\widehat{\mathfrak{V}}}_{s}^{n-1}(x)-V_s(x)\|_{\infty}].\nonumber
\end{eqnarray*}

\medskip

\noindent \emph{Step 3}. From (\ref{proof_induction}) in \emph{Step 2}, we use induction on $n$ to show that for $1\le t\le T$ and  $ n\ge T-t+1$,
\begin{eqnarray}
\label{ec_ind}
\mathbb{E}\Big[1_{A(\delta,n)^c}\cdot\Big\|\underline{\widehat{\mathfrak{V}}}_t^n(x)-V_t(x)\Big\|_{\infty}\Big]
	\le\Big(1+2l_M+C\frac{M^{3/2}}{L^{1/4}}\Big)^{T-t}\Big[(1+l_M)\Delta+C\frac{M^{3/2}}{L^{1/4}}\Big].
\end{eqnarray}

We omit the calculation detail in the interest of space.
\medskip

\noindent \emph{Step 4}.
In light of the definition of $\underline{\widehat{\mathfrak{V}}}^{T+1}_{0}(x)$,
\[
\underline{\widehat{\mathfrak{V}}}^{T+1}_{0}(x)=\max\Big\{-K,\min\big\{K,\frac{1}{L}\sum_{l=1}^{L}\mathfrak{J}^{(l)}_{0, T+1}\big\}\Big\},
\]
we choose constant $K$ such that $K\ge|V_0(x)|$. According to Lemma \ref{lem_ineq}, we have
\begin{eqnarray*}
\Big|\underline{\widehat{\mathfrak{V}}}^{T+1}_{0}(x)-V_0(x)\Big|\le  \Big|\frac{1}{L}\sum_{l=1}^{L}\mathfrak{J}^{(l)}_{0, T+1}-V_0(x)\Big|\le \frac{1}{L}\sum_{l=1}^{L}\Big|\mathfrak{J}^{(l)}_{0,T+1}-V_0(x)\Big|.
\end{eqnarray*}
Again we use the Lemma \ref{lem_J},
\begin{eqnarray}
\label{ec_4}
\mathbb{E}\Big[1_{A(\delta,T)^c}\cdot\Big|\underline{\widehat{\mathfrak{V}}}_0^{T+1}(x)-V_0(x)\Big|\Big]
\le2 \sum_{t=1}^T\mathbb{E}\Big[ 1_{A(\delta,T)^c}\cdot\Big\|\underline{\widehat{\mathfrak{V}}}_{t}^{T}(x)-V_t(x)\Big\|_{\infty}\Big].
\end{eqnarray}
We sum the inequality (\ref{ec_ind}) from $t=1$ to $T$ in iteration $T$ and derive that
%\begin{eqnarray}
%\label{ec_3}
%\sum_{t=1}^T\mathbb{E}\Big[1_{A(\delta,T)^c}\cdot\Big\|\underline{\widehat{\mathfrak{V}}}_t^{T}(x)-V_0(x)\Big\|_{\infty}\Big]
%	\le\frac{\Big(1+2l_M+C\frac{M^{3/2}}{L^{1/4}}\Big)^{T}-1}{2l_M+C\frac{M^{3/2}}{L^{1/4}}}\Big[(1+l_M)\Delta+C\frac{M^{3/2}}{L^{1/4}}\Big].
%\end{eqnarray}
%Accordingly,
%\begin{eqnarray*}
%\mathbb{E}\Big[1_{A(\delta,T)^c}\cdot\Big|\underline{\widehat{\mathfrak{V}}}_0^{T+1}(x)-V_0(x)\Big|\Big]\le
%\frac{\Big(1+2l_M+C\frac{M^{3/2}}{L^{1/4}}\Big)^{T}-1}{l_M+C\frac{M^{3/2}}{2L^{1/4}}}\Big[(1+l_M)\Delta+C\frac{M^{3/2}}{L^{1/4}}\Big].
%\end{eqnarray*}
%Note that $l_M\ge1$. This can be obtained if we replace $f$ by $\sum_{m=1}^M\psi_m(x)$ in definition (\ref{lm}). As a consequence, the above can be further simplified as
\begin{eqnarray*}
\mathbb{E}\Big[1_{A(\delta,T)^c}\cdot\Big|\underline{\widehat{\mathfrak{V}}}_0^{T+1}(x)-V_0(x)\Big|\Big]\le
\Big(1+2l_M+C\frac{M^{3/2}}{L^{1/4}}\Big)^{T}\Big[(1+l_M)\Delta+C\frac{M^{3/2}}{L^{1/4}}\Big].
\end{eqnarray*}
\medskip
\noindent\emph{Step 5}. By combining the result of (\ref{ec_1}) and (\ref{ec_4}), we conclude that
\[
\mathbb{E}\Big[\big|\underline{\widehat{\mathfrak{V}}}^{T+1}_0(x)-V_0(x)\Big|\Big]\le CT(T+1)M \exp\Big(-\frac{L^{1/2}}{CM^{3}}\Big)+\Big(1+2l_M+C\frac{M^{3/2}}{L^{1/4}}\Big)^{T}\Big[(1+l_M)\Delta+C\frac{M^{3/2}}{L^{1/4}} \Big].
\]
For sufficient small $\alpha$, we have
\[
M \exp\Big(-\frac{L^{1/2}}{CM^{3}}\Big)\le\Big(1+2l_M+C\frac{M^{3/2}}{L^{1/4}}\Big)^{T}\frac{M^{3/2}}{L^{1/4}}.
\]
By adjusting the constant $C$ properly, we obtain the result in Theorem \ref{thm:error}.\hfill$\square$\\

\section{Supplementary Materials to Section \ref{sec:num}}
\label{app:num}

\smallskip

\subsection{Optimal Order Execution Problem:}
\begin{itemize}[leftmargin=0in]
\item[-] \textbf{The objective function}:

\noindent It is easy to see that minimizing (\ref{exm1:objective}) is equivalent to minimizing
\[
\mathbb{E}\left[\sum_{t=1}^{T}\mathbf{P}^{tr}_{t}\mathbf{S}_{t}-\tilde{\mathbf{P}}^{tr}_{0}\bar{\mathbf{R}}\right].
\]
Note that the constant $\tilde{\mathbf{P}}^{tr}_{0}\bar{\mathbf{R}}$ stands for the cost that the trader would pay for $\bar{\mathbf{R}}$ shares of assets without the price impacts.
This difference thus represents the implementation shortfall of a specific strategy, namely how much more costs the trader may incur during the course of fulfilling the execution
target. In the following lemma, we show that it equals (\ref{exm1:objective_new}).

\begin{lem}
\label{exm1:tgp}
For the trader's problem
\begin{eqnarray}
\label{app:num_t}
\min_{\{\mathbf{S}_{t}, 1 \le t \le T\}}\mathbb{E}\left[\left(\sum_{t=1}^{T}\mathbf{P}^{tr}_{t}\mathbf{S}_{t}-\tilde{\mathbf{P}}^{tr}_{0}\bar{\mathbf{R}}\right)\right],\end{eqnarray}
subject to the constraints (\ref{exm1:dyn1}-\ref{exm1:dyn3}),
it is equivalent to
\begin{eqnarray*}
\min_{\{\mathbf{S}_{t}, 1 \le t \le T\}}\mathbb{E}\left[\sum_{t=1}^{T}\mathbf{S}_{t}^{tr}h(\mathbf{S}_{t})+\sum_{t=0}^{T-1}(\tilde{\mathbf{P}}_{t+1}-\tilde{\mathbf{P}}_{t})^{tr}\mathbf{R}_{t+1}\right].
\end{eqnarray*}
\end{lem}
\textit{Proof of Lemma \ref{exm1:tgp}.}
Using the relationship (\ref{exm1:dyn2}), we observe that
\begin{eqnarray}
\label{exm1:obj_new1}
\sum_{t=1}^{T}\mathbf{P}^{tr}_{t}\mathbf{S}_{t}-\tilde{\mathbf{P}}^{tr}_{0}\bar{\mathbf{R}}&=&\sum_{t=1}^{T}(\tilde{\mathbf{P}}_{t}
+h(\mathbf{S}_t))^{tr}\mathbf{S}_{t}-\tilde{\mathbf{P}}^{tr}_{0}\bar{\mathbf{R}}\nonumber\\
&=&\sum_{t=1}^{T}\mathbf{S}_{t}^{tr}h(\mathbf{S}_{t})
+\sum_{t=1}^{T}\tilde{\mathbf{P}}^{tr}_{t}\mathbf{S}_{t}-\tilde{\mathbf{P}}^{tr}_{0}\bar{\mathbf{R}}.
\end{eqnarray}
In addition, applying Abel's summation-by-part formula to $\sum_{t=1}^{T}\tilde{\mathbf{P}}^{tr}_{t}\mathbf{S}_{t}$, we know that
\begin{eqnarray}
\label{exm1:obj_new2}
\sum_{t=1}^{T}\tilde{\mathbf{P}}^{tr}_{t}\mathbf{S}_{t}&=&\tilde{\mathbf{P}}^{tr}_{0}\left(\sum_{t=1}^{T}\mathbf{S}_{t}\right)+\sum_{t=0}^{T-1}\left((\tilde{\mathbf{P}}_{t+1}
-\tilde{\mathbf{P}}_{t})^{tr}\cdot\sum_{j=t+1}^{T}\mathbf{S}_{j}\right)\nonumber\\
&=&\tilde{\mathbf{P}}^{tr}_{0}\bar{\mathbf{R}}+\sum_{t=0}^{T-1}(\tilde{\mathbf{P}}_{t+1}-\tilde{\mathbf{P}}_{t})^{tr}\mathbf{R}_{t+1}.
\end{eqnarray}
Thus, with (\ref{exm1:obj_new1}) and (\ref{exm1:obj_new2}), we have
\begin{eqnarray*}
\sum_{t=1}^{T}\mathbf{P}^{tr}_{t}\mathbf{S}_{t}-\tilde{\mathbf{P}}^{tr}_{0}\bar{\mathbf{R}}=\sum_{t=1}^{T}\mathbf{S}_{t}^{tr}h(\mathbf{S}_{t})+\sum_{t=0}^{T-1}(\tilde{\mathbf{P}}_{t+1}-\tilde{\mathbf{P}}_{t})^{tr}\mathbf{R}_{t+1}.
\end{eqnarray*}
This verifies the equivalence of these two objective functions. Note the new value function doesn't depend on the variable $\mathbf{P}$.\hfill$\square$\\

\item[-] \textbf{The auxiliary LQC problem}:

\noindent If we ignore the temporary impact $h(\mathbf{S}_t)$ and remove the nonnegative constraint $\mathbf{S}_t\ge0$, the problem (\ref{exm1:objective_new}) with the constraints (\ref{exm1:dyn1}-\ref{exm1:dyn3}) is equivalent to LQC problem. According to the discussion in Appendix \ref{LQ_iteration}, the value function $V_t(\mathbf{X}_t, \mathbf{R}_t)$ and policy $\mathbf{S}^{*}_t$ are:
\begin{eqnarray}
	&&V_t(\mathbf{X}_t, \mathbf{R}_t)=\mathbf{X}_t^{tr}\mathbf{W}_t\mathbf{X}_t+\mathbf{R}_t^{tr}\mathbf{Q}_t\mathbf{R}_t+\mathbf{R}_t^{tr}\mathbf{K}_t\mathbf{X}_t+\mathbf{H}_t,\\
	\label{app:exm1p}
	&&\mathbf{S}^{*}_t(\mathbf{X}_t, \mathbf{R}_t)=(\mathbf{I}-\frac{1}{2}\mathbf{Q}_{t+1}^{-1}\mathbf{A}^{tr})\mathbf{R}_t+\frac{1}{2}\mathbf{Q}_{t+1}^{-1}\mathbf{K}_{t+1}\mathbf{C}\mathbf{X}_t,
\end{eqnarray}
with
\begin{eqnarray}
\label{app:exm1}
	&&\mathbf{Q}_t=-\frac{1}{4}\mathbf{A}\mathbf{Q}_{t+1}^{-1}\mathbf{A}^{tr}+\frac{1}{2}(\mathbf{A}+\mathbf{A}^{tr}),\quad \mathbf{Q}_{T}=\frac{1}{2}(\mathbf{A}+\mathbf{A}^{tr}).\\
	&&\mathbf{W}_t=\mathbf{C}^{tr}\mathbf{W}_{t+1}\mathbf{C}-\frac{1}{4}\mathbf{C}^{tr}\mathbf{K}_{t+1}^{tr}\mathbf{Q}_{t+1}^{-1}\mathbf{K}_{t+1}\mathbf{C},\quad \mathbf{W}_{T}=0.\\
	&&\mathbf{K}_t=\mathbf{B}+\frac{1}{2}\mathbf{A}\mathbf{Q}_{t+1}^{-1}\mathbf{K}_{t+1}\mathbf{C},\quad \mathbf{K}_{T}=B.\\
	\label{app:exm2}
	&&\mathbf{H}_t=\mathbf{H}_{t+1}+\mathbb{E}[\boldsymbol{\eta}^{tr}\mathbf{W}_{t+1}\boldsymbol{\eta}],\quad \mathbf{H}_{T}=0.
\end{eqnarray}
Specially if the matrix $\mathbf{A}$ is symmetric, the optimal policy (\ref{app:exm1p}) can be simplified as
\[\mathbf{S}^{*}_t(\mathbf{X}_t, \mathbf{R}_t)=\frac{1}{2}\mathbf{Q}_{t+1}^{-1}\mathbf{K}_{t+1}\mathbf{C}\mathbf{X}_t+\frac{1}{T-t+1}\mathbf{R}_t.\]

\item[-] \textbf{Parameter setting}:

\noindent To illustrate the numerical results, we consider a case with three assets and a signal vector of two variables. Assume the trader
wants to buy $1 \times 10^5$ shares for each asset within $T=20$ periods, i.e., $\bar{R}_{i}= 1 \times 10^5$ for $i=1, 2, 3$. The parameter matrices pertinent to the temporary and permanent impacts are supposed to
\begin{eqnarray*}
	&&\mathbf{A}\!=\!\left[
	\begin{smallmatrix}
		30&   7&  3\\
		7&  25&  -5\\
		3&  -5&  20
	\end{smallmatrix}
	\right]\!\times\!10^{-6}, \
	\mathbf{B}\!=\!\left[
	\begin{smallmatrix}
		5& 2 \\
		3 & 2 \\
		1 & 4
	\end{smallmatrix}
	\right], \
	\mathbf{C}\!=\!\delta\!\times\!\left[
	\begin{smallmatrix}
		0.8& 0.1 \\
		0.2 & 0.6
	\end{smallmatrix}
	\right], \
	\mathbf{D}\!=\!\left[
	\begin{smallmatrix}
		2\lambda & 0 & 0 \\
		0 & 2\lambda & 0 \\
		0 & 0 & 2\lambda
	\end{smallmatrix}
	\right]\!\times\!10^{-5},\
	\Sigma_{\eta}\!=\!\left[
	\begin{smallmatrix}
		1.0& 0.2 \\
		0.2 & 0.8
	\end{smallmatrix}
	\right].
\end{eqnarray*}
Here we parametrize matrix $\mathbf{D}$ by $\lambda$ so that we can examine the effect of the temporary price impact on the optimal execution strategies by varying $\lambda$.

\end{itemize}

\subsection{Inventory Management Problem:}
\begin{itemize}[leftmargin=0in]
\item[-]\textbf{Parameters}:
 The model parameters' values used in the experiments are given by
\[
h=1,\quad m=4\ \textrm{or}\ 9,\quad p=9\ \textrm{or}\ 19,\quad \gamma=1,\quad T=30,\quad \mathbf{x}_1=\mathbf{0}.
\]

\item[-]\textbf{Basis functions}:

\noindent For $L=4$, we choose the basis function set as
\begin{eqnarray*}
	&
  \Big\{1,\quad (x_{i,t})_{0\le i\le 3},\quad \mathbb{E}[(x_{0,t}-\tilde{d}_0)^+],\quad \mathbb{E}[((x_{0,t}-\tilde{d}_0)^++x_{1,t}-\tilde{d}_1)^+],
  \\
	&
  \resizebox{.95\hsize}{!}{$
  \mathbb{E}[(((x_{0,t}-\tilde{d}_0)^++x_{1,t}-\tilde{d}_1)^++x_{2,t}-\tilde{d}_2)^+],\quad \mathbb{E}[((((x_{0,t}-\tilde{d}_0)^++x_{1t}-\tilde{d}_1)^++x_{2,t}-\tilde{d}_2)^++x_{3,t}-\tilde{d}_3)^+],
  $}
  \\
	&
  \resizebox{.9\hsize}{!}{$
  \mathbb{E}[(((x_{1,t}-\tilde{d}_1)^++x_{2,t}-\tilde{d}_2)^++x_{3,t}-\tilde{d}_3)^+],\quad \mathbb{E}[((x_{2,t}-\tilde{d}_2)^++x_{3,t}-\tilde{d}_3)^+],\quad\mathbb{E}[(x_{3,t}-\tilde{d}_3)^+]\Big\}.
  $}
\end{eqnarray*}
The expectation is taken over $(\tilde{d}_i)_{0\le i\le 3}$, which have the same distribution with $d_t$ in the system. For $L=10$, we choose 30 basis functions in similar manner as $L=4$. That is, constant 1, one order function $(x_{i,t})_{0\le i\le 9}$, the expectation in iteration form from $\mathbb{E}[(x_{0,t}-\tilde{d}_0)^+]$ to $\mathbb{E}[((x_{0,t}-\tilde{d}_0)^+\cdots+x_{9,t}-\tilde{d}_9)^+]$, and the reverse form from  $\mathbb{E}[(x_{9,t}-\tilde{d}_9)^+]$ to $\mathbb{E}[((x_{1,t}-\tilde{d}_1)^+\cdots+x_{9,t}-\tilde{d}_9)^+]$.

\item[-]\textbf{Quasi Monte Carlo}:

\noindent As mentioned in the main body, for $L=10$, we choose low-discrepancy sequences to perform the nested simulations. To illustrate this, we note that the expectation of basis functions can be written in the form of
\[
H(\mathbf{x}_t)=\mathbb{E}[g(\mathbf{x}_t,\mathbf{d})]
\]
for some function $g(\cdot)$ and where state $\mathbf{x}=[x_{0,t},\cdots,x_{9,t}]$, geometric distribution $\mathbf{d}=[\tilde{d}_0,\cdots, \tilde{d}_{9}]$. Using the inverse transform approach we can easily rewrite $H(\mathbf{x}_t)$ as
\[
H(\mathbf{x}_t)=\mathbb{E}[g(\mathbf{x}_t,\lfloor \log(\mathbf{U})/\log(1-p)\rfloor)],
\]
where $\mathbf{U}$ is $10$-dimensional vector of independent uniform random variables in range $(0,1)$ and $\lfloor N\rfloor $ stands for the largest integer which is no bigger than $N$. We could perform this expectation with respect to $\mathbf{U}$ by using the low discrepancy sequence. Here we choose 2047 points in Sobol sequence, $U_1,\cdots,U_{2047}$, and approximate
\[
H(\mathbf{x}_t)\approx\sum_{i=1}^{2047}g(\mathbf{x}_t,\lfloor \log(U_i)/\log(1-p)\rfloor).
\]
The detailed discussion about this method, one may refer to Chapter 5 in \cite{gl}.

\end{itemize}

%%%%%%%%%%%%%%%%%
\end{document}